\documentclass{article}

\usepackage[preprint]{neurips_2020_modified}

\usepackage[utf8]{inputenc} 
\usepackage[T1]{fontenc}    
\usepackage{hyperref}       
\usepackage{url}            
\usepackage{booktabs}       
\usepackage{amsfonts}       
\usepackage{nicefrac}       
\usepackage{microtype}      

\usepackage{amssymb, amsmath, listings, mathrsfs}
\usepackage{algorithm}
\usepackage{graphicx}
\usepackage{xcolor}
\usepackage[noend]{algpseudocode}
\floatstyle{plaintop}
\usepackage[labelfont=bf]{caption}
\usepackage{thmtools,thm-restate}
\usepackage{mdframed}

\def\nr{\par \noindent}

\def\Def{\stackrel{\mathrm{def}}{=}}

\def\dom{{\rm dom \,}}
\def\beq{\begin{equation}}
\def\eeq{\end{equation}}

\def\R{\mathbb{R}}

\def\E{\mathbb{E}}
\def\O{\mathcal{O}}

\def\H{\mathcal{H}}

\def\Argmin{\mathop{\rm Argmin}}

\newcommand{\refLE}[1]{\ensuremath{\stackrel{(\ref{#1})}{\leq}}}
\newcommand{\refEQ}[1]{\ensuremath{\stackrel{(\ref{#1})}{=}}}
\newcommand{\refGE}[1]{\ensuremath{\stackrel{(\ref{#1})}{\geq}}}

\newtheorem{theorem}{Theorem}
\newmdtheoremenv{framedtheorem}{Theorem}
\newtheorem{lemma}{Lemma}
\newmdtheoremenv{framedlemma}{Lemma}
\newtheorem{corollary}{Corollary}

\newtheorem{assumption}{Assumption}
\newtheorem{definition}{Definition}

\newtheorem{example}{Example}
\newtheorem{remark}{Remark}
\newcommand{\proof}{\bf Proof: \rm \nr}
\newcommand{\qed}{\hfill $\Box$ \nr \medskip}

\def\ba{\begin{array}}
    \def\ea{\end{array}}
\def\beann{\begin{eqnarray*}}
    \def\eeann{\end{eqnarray*}}
\def\bea{\begin{eqnarray}}
\def\eea{\end{eqnarray}}

\def\BT{\begin{theorem}}
    \def\ET{\end{theorem}}
\def\BFT{\begin{framedtheorem}}
    \def\EFT{\end{framedtheorem}}
\def\BL{\begin{lemma}}
    \def\EL{\end{lemma}}
\def\BFL{\begin{framedlemma}}
    \def\EFL{\end{framedlemma}}
\def\BC{\begin{corollary}}
    \def\EC{\end{corollary}}
\def\BE{\begin{example}}
    \def\EE{\end{example}}
\def\BD{\begin{definition}}
    \def\ED{\end{definition}}
\def\BR{\begin{remark}}
    \def\ER{\end{remark}}
\def\BAS{\begin{assumption}}
    \def\EAS{\end{assumption}}
\def\BI{\begin{itemize}}
    \def\EI{\end{itemize}}

\def\BMP{\begin{minipage}{9.5cm}}
    \def\EMP{\end{minipage}}
\def\MPT{\begin{minipage}{11.5cm}}
    \def\EPT{\end{minipage}}

\def\la{\langle}
\def\ra{\rangle}

\title{Convex optimization based on global lower second-order models}

\author{%
  Nikita Doikov\thanks{Institute of Information and Communication Technologies,
    Electronics and Applied Mathematics (ICTEAM)} \\
  Catholic University of Louvain, \\
  Louvain-la-Neuve, Belgium\\
  \texttt{Nikita.Doikov@uclouvain.be}
  \And
  Yurii Nesterov\thanks{Center for Operations Research and Econometrics (CORE)} \\
  Catholic University of Louvain, \\
  Louvain-la-Neuve, Belgium\\
  \texttt{Yurii.Nesterov@uclouvain.be}
  \And
  \\
}

\begin{document}

\maketitle

\begin{abstract}
  In this paper, we present new second-order algorithms for composite convex optimization,
  called Contracting-domain Newton methods.
  These algorithms are affine-invariant and based on global
  second-order lower approximation for the smooth component of the objective.
  Our approach has an interpretation both as
  a second-order generalization of the conditional gradient method,
  or as a variant of trust-region scheme.
  Under the assumption, that the problem domain is bounded, we prove
  $\O(1/k^{2})$ global rate of convergence in functional residual,
  where $k$ is the iteration counter,
  minimizing convex functions with Lipschitz continuous Hessian.
  This significantly improves the previously known bound $\O(1/k)$ for this type of algorithms.
  Additionally, we propose a stochastic extension of our method,
  and present computational results
  for solving empirical risk minimization problem.

\end{abstract}

\section{Introduction}

Classical Newton method is one of the most popular
optimization schemes for solving ill-conditioned problems.
The method has very fast quadratic convergence, provided
that the starting point is sufficiently close to the
optimum~\cite{bennett1916newton,kantorovich1948functional,nesterov2018lectures}.
However, the questions related to its global behaviour for
a wide class of functions are still open, being in the
area of active research.

The significant progress in this direction was made
after~\cite{nesterov2006cubic}, where Cubic regularization
of Newton method with its global complexity bounds were
justified. The main idea of~\cite{nesterov2006cubic} is to
use a global \textit{upper} approximation model of the
objective, which is the second-order Taylor's polynomial
augmented by a cubic term. The next point in the iteration
process is defined as the minimizer of this model. Cubic
Newton attains global convergence for convex functions
with Lipschitz continuous Hessian. The rate of convergence
in functional residual is of the order $\O(1 / k^2)$ (here
and later on, $k$ is the iteration counter). This is much
faster than the classical $\O(1 / k)$ rate of the Gradient
Method~\cite{nesterov2018lectures}. Later on,
accelerated~\cite{nesterov2008accelerating},
adaptive~\cite{cartis2011adaptive1,cartis2011adaptive2}
and
universal~\cite{grapiglia2017regularized,doikov2019minimizing,grapiglia2019accelerated}
second-order schemes based on cubic regularization were
developed. Randomized versions of Cubic Newton, suitable
for solving high-dimensional problems were proposed
in~\cite{doikov2018randomized,hanzely2020stochastic}.

Another line of results on global convergence of Newton
method is mainly related to the framework of
self-concordant
functions~\cite{nesterov1994interior,nesterov2018lectures}.
This class is affine-invariant. From the global
perspective, it provides us with an \textit{upper}
second-order approximation of the objective, which
naturally leads to the Damped Newton Method. Several new
results are related to its analysis for generalized
self-concordant
functions~\cite{bach2010self,sun2019generalized}, and the
notion of Hessian stability~\cite{karimireddy2018global}.
However, for more refined problem classes, we can often
obtain much better complexity estimates, by using the
cubic regularization
technique~\cite{dvurechensky2018global}.

In this paper, we investigate a different approach, which
is motivated by a new global second-order \textit{lower}
model of the objective function, introduced in
Section~\ref{Sec:Lower}.

We incorporate this model into a new second-order
optimization algorithm, called Contracting-Domain Newton
Method (Section~\ref{Sec:ContrNewton}). At every
iteration, it minimizes a lower approximation of the
smooth component of the objective, augmented by a
composite term. The next point is defined as a convex
combination of the minimizer, and the previous point. By
its nature, it is similar to the scheme of Conditional
Gradient Method (or, Frank-Wolfe
algorithm,~\cite{frank1956algorithm,nesterov2018complexity}).
Under assumption of boundedness of the problem domain, for
convex functions with H\"older continuous Hessian of
degree $\nu \in [0, 1]$, we establish its $\O(1/k^{1 +
\nu})$ global rate of convergence in functional residual.
In the case $\nu = 1$, for the class of convex function
with Lipschitz continuous Hessian, this gives $\O(1/k^2)$
rate of convergence. As compared with Cubic Newton, the
new method is affine-invariant and universal, since it
does not depend on the norms and parameters of the problem
class. When the composite component is strongly convex
(with respect to arbitrary norm), we show $\O(1/k^{2 +
2\nu})$ rate for a universal scheme. If the parameters of
problem class are known, we can prove a global linear
convergence. 
We also provide different trust-region
interpretations for our algorithm.

In Section~\ref{Sec:Aggr}, we present aggregated models,
which accumulate second-order information into quadratic
Estimating Functions~\cite{nesterov2018lectures}. This
leads to another optimization process, called Aggregating
Newton Method, with the global convergence of the same
order $\O(1/k^{1 + \nu})$ as for general convex case. The
latter method can be seen as a second-order counterpart of
the dual averaging gradient
schemes~\cite{nesterov2009primal,nesterov2013gradient}.

In Section~\ref{Sec:Stoch}, we consider the problem of
finite-sum minimization. We propose stochastic extensions
of our method. During the iterations of the basic variant,
we need to increase the batch size for randomized
estimates of gradients and Hessians up to the order
$\O(k^4)$ and $\O(k^2)$ respectively. Using the
\textit{variance reduction} 
technique for the gradients,
we reduce the batch size up to the
level $\O(k^2)$ for both estimates. At the same time, the global convergence
rate of the resulting methods
is of the order $\O(1/k^{2})$,
as for general convex functions with Lipschitz continuous Hessian.

Section~\ref{Sec:Experiments} contains numerical
experiments. Section~\ref{Sec:Discussion} contains some final remarks.  
All necessary proofs are provided in
Appendix.

\section{Problem formulation and notations}
\label{Sec:Notation}

Our goal is to solve the following {\em composite} convex
minimization problem:
\beq \label{MainProblem}
\ba{rcl}
\min\limits_{x} F(x) & := & f(x) + \psi(x),
\ea
\eeq
where $\psi: \R^n \to \R \cup \{ +\infty \}$ is a
\textit{simple} proper closed convex function, and
function $f$ is convex and twice continuously
differentiable at every point $x \in \dom \psi$. Let us
fix an arbitrary (possibly non-Euclidean) norm $\| \cdot
\|$ on $\R^n$. We denote by $D$ the corresponding diameter
of $\dom \psi$:
\beq \label{DiamDef}
\ba{rcl}
D & := & \sup\limits_{x, y \in \dom \psi} \|x - y\|.
\ea
\eeq
Our main assumption on problem~\eqref{MainProblem} is that
$\dom \psi$ is bounded:
\beq \label{BoundedDom}
\ba{rcl}
D & < & +\infty.
\ea
\eeq
The most important example of $\psi$ is $\{0,
+\infty\}$-indicator of a simple compact convex set $Q =
\dom \psi$. In particular, for a ball in
$\|\cdot\|_p$-norm with $p \geq 1$, this is
\beq \label{PsiBall}
\ba{rcl}
\psi(x) & = & \begin{cases}
0, & \|x\|_p := \Bigl(\sum_{i = 1}^n |x^{(i)}|^p\Bigr)^{1/p}  \leq \frac{D}{2}, \\
+\infty, & \text{else}.
\end{cases}
\ea
\eeq
From the machine learning perspective, $D$ is usually considered as a \textit{regularization parameter} in this setting. We denote by $\la \cdot, \cdot \ra$ the standard scalar product of two vectors,
$x, y \in \R^n$: 
$$
\ba{rcl}
\la x, y \ra  & := & \sum_{i = 1}^n x^{(i)} y^{(i)}.
\ea
$$
For function $f$, we denote its gradient by $\nabla f(x)
\in \R^n$, and its Hessian matrix by $\nabla^2 f(x) \in
\R^{n \times n}$. Having fixed the norm $\| \cdot \|$ for
primal variables $x \in \R^n$, the \textit{dual} norm 
can be defined in the standard way:
$$
\ba{rcl}
\| s \|_{*}  & := &  \sup\limits_{h \in \R^n : \|h\| \leq 1} \la s, h \ra.
\ea
$$
The dual norm is necessary for measuring the size of
gradients. For a matrix $A \in \R^{n \times n}$, 
we use the corresponding induced operator norm, defined as
$$
\ba{rcl}
\| A \| & := & \sup\limits_{h \in \R^n : \|h\| \leq 1} \| Ah \|_{*}.
\ea
$$

\section{Second-order lower model of objective function}
\label{Sec:Lower}

To characterize the complexity of
problem~\eqref{MainProblem}, we need to introduce some
assumptions on the growth of derivatives. Let us assume
that the Hessian of $f$ is H\"older continuous of degree
$\nu \in [0, 1]$ on $\dom \psi$:
\beq \label{HolderHessian}
\ba{rcl}
\| \nabla^2 f(x) - \nabla^2 f(y) \| & \leq & H_{\nu} \|x - y\|^{\nu},
\qquad x, y \in \dom \psi.
\ea
\eeq
The actual parameters of this problem class may be
unknown. However, we assume that for \textit{some} $\nu
\in [0, 1]$ inequality~\eqref{HolderHessian} is satisfied
with corresponding constant $0 \leq H_{\nu} < +\infty$.
The direct consequence of~\eqref{HolderHessian} is the
following global bounds for Taylor's approximation,
for all $x, y \in \dom \psi$
\beq \label{GradHolderBound}
\ba{rcl}
\| \nabla f(y) - \nabla f(x) - \nabla^2 f(x)(y - x) \|_{*}
& \leq & \frac{H_{\nu} \|y - x\|^{1 + \nu}}{1 + \nu},
\ea
\eeq
\beq \label{FuncHolderBound}
\ba{rcl}
| f(y) - f(x) - \la \nabla f(x), y - x \ra - \frac{1}{2} \la \nabla^2 f(x)(y - x), y - x \ra |
& \leq & \frac{H_{\nu} \|y - x\|^{2 + \nu}}{(1 + \nu)(2 + \nu)}.
\ea
\eeq
Recall, that in addition to~\eqref{HolderHessian}, we assume that $f$ is \textit{convex}:
\beq \label{fConv}
\ba{rcl}
f(y) & \geq & f(x) + \la \nabla f(x), y - x \ra, \qquad x,
y \in \dom \psi.
\ea
\eeq
Employing both smoothness and convexity, we are able to enhance this global lower
bound, as follows.

\BFL \label{LemmaLB} For all $x, y \in \dom \psi$ and $t \in [0, 1]$, it holds
\beq \label{FuncLowerBound}
\ba{rcl}
f(y) & \geq & f(x) + \la \nabla f(x), y - x \ra
+ \frac{t}{2} \la \nabla^2 f(x)(y - x), y - x \ra
- \frac{t^{1 + \nu} H_{\nu} \|y - x\|^{2 + \nu}}{(1 + \nu)(2 + \nu)}.
\ea
\eeq
\EFL

\begin{figure}[h!]
    \centering
    \includegraphics[width=0.5\textwidth ]{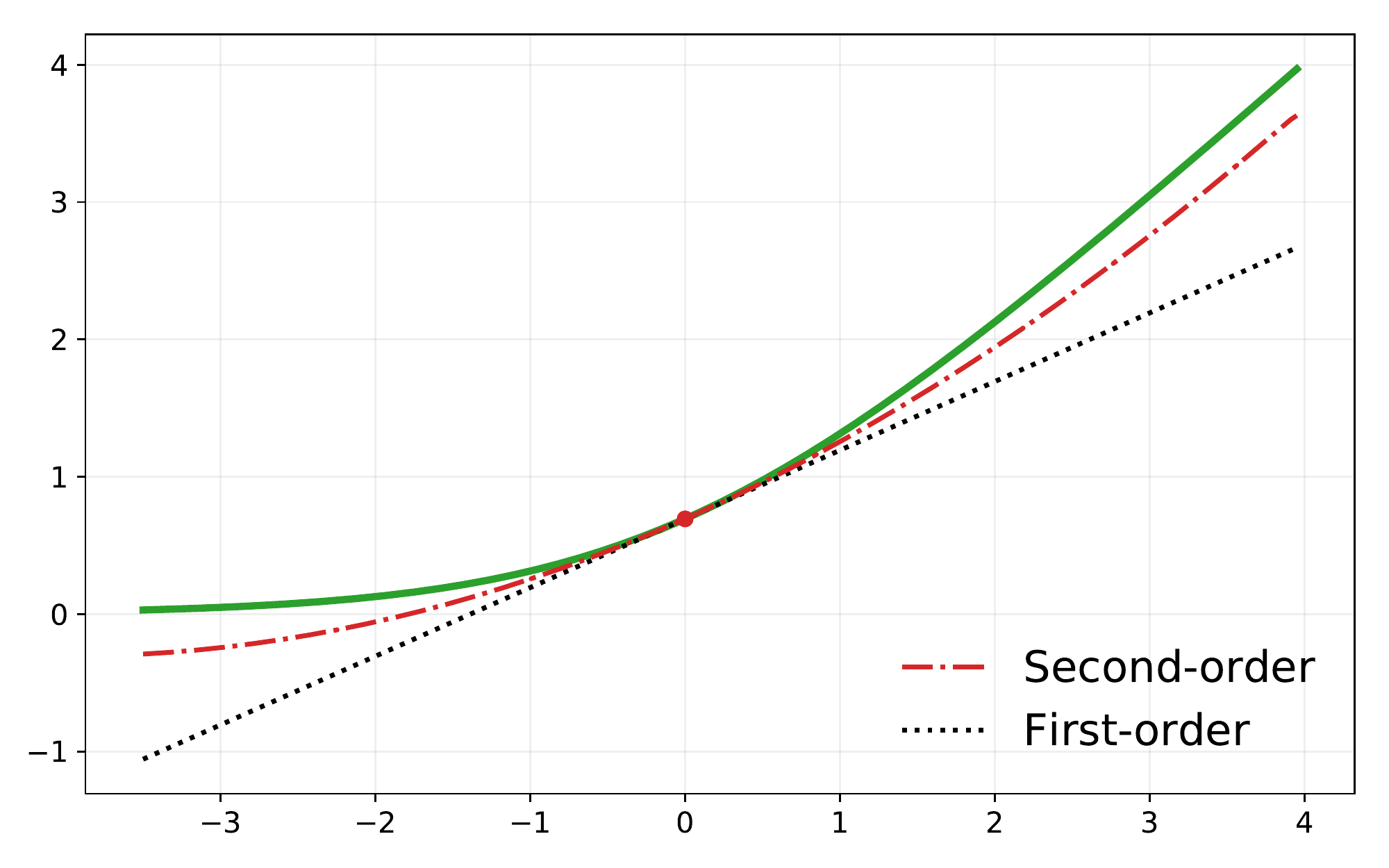}
    \caption{Global lower bounds for logistic regression loss, $f(x) = \log(1 + \exp(x))$. }
    \label{Fig:lower_bounds}
\end{figure}

Note that the right-hand side of~\eqref{FuncLowerBound} is
concave in $t \in [0, 1]$, and for $t = 0$ we obtain the
standard first-order lower bound. The maximization
of~\eqref{FuncLowerBound} over $t$ gives
\beq \label{SecondLower}
\ba{rcl}
f(y) & \geq &f(x) + \la \nabla f(x), y - x \ra
+  \frac{\bar{\gamma}_{x, y}}{2} \la \nabla^2 f(x)(y - x), y - x \ra,
\ea
\eeq
with
$$
\ba{rcl}
\bar{\gamma}_{x, y} & := &
\frac{\nu}{1 + \nu} \min\Bigl\{1,
\frac{(2 + \nu)\la \nabla^2 f(x)(y - x), y - x \ra}{2 H_{\nu}\|y - x\|^{2 + \nu}}
\Bigr\}^{\frac{1}{\nu}},
\quad x \not= y, \quad \nu \in (0, 1].
\ea
$$
Thus, \eqref{SecondLower} is always \textit{tighter}
than~\eqref{fConv}, employing additional \textit{global
second-order information}. The relationship between them
is shown on Figure~\ref{Fig:lower_bounds}. Hence, it seems
natural to incorporate the second-order lower bounds into
optimization schemes.

\section{Contracting-Domain Newton Methods}
\label{Sec:ContrNewton}

Let us introduce a general scheme of
\textit{Contracting-Domain Newton Method}, which is based
on global second-order lower bounds. Note, that the right
hand side of~\eqref{SecondLower} is nonconvex in $y$.
Hence, it can hardly be used directly in a computational
algorithm. To tackle this issue, we use a sequence of
contracting coefficients $\{ \gamma_k \}_{k \geq 0}$. Each
coefficient $\gamma_k \in (0, 1]$ can be seen as an
appropriate substitute of $\bar{\gamma}_{x, y}$
in~\eqref{SecondLower}. Then, we minimize the
corresponding global lower bound augmented by the
composite component $\psi(\cdot)$. The next point is taken
as a convex combination of the minimizer and the current
point. Let us present this method formally, as
Algorithm~\ref{ContrDomNewton}.

\begin{algorithm}[H]
    \caption{Contracting-Domain Newton Method, I}
    \label{ContrDomNewton}
    \begin{algorithmic}[1]
        \Require Choose $x_0 \in \dom \psi$.
        \Ensure {$k \geq 0$.}
        \State Pick up $\gamma_k \in (0, 1]$.
        \State Compute \vspace*{-10pt}
        $$
        \ba{rcl}
        v_{k + 1} & \in & \Argmin\limits_{y} \Bigl\{
        \la \nabla f(x_k), y - x_k \ra
        \; + \; \frac{\gamma_k}{2} \la \nabla^2 f(x_k)(y - x_k), y - x_k \ra
        \; + \; \psi(y)
        \Bigr\}.
        \ea
        $$
        \vspace*{-12pt}
        \State Set $x_{k + 1} \; := \; x_k + \gamma_k(v_{k + 1} - x_k)$.
    \end{algorithmic}
\end{algorithm}

There is a clear connection of this method with
Frank-Wolfe algorithm, \cite{frank1956algorithm}. Indeed,
instead of the standard first-order
approximation~\eqref{fConv}, we use the lower global
quadratic model. Thus, as compared with the gradient
methods, every iteration of Algorithm~\ref{ContrDomNewton}
is more expensive. However, this is a standard situation
with the second-order schemes (see the below discussion on
the iteration complexity). At the same time, our method is
\textit{affine-invariant}, since it does not depend on the
norms.

It is clear, that for $\gamma_k \equiv 1$ we obtain
iterations of the classical Newton method. Its local
quadratic convergence for composite optimization problems
was established in~\cite{lee2014proximal}. However, for
the global convergence, we need to adjust the contracting
coefficients accordingly. To state the global convergence
result, let us introduce the following linear
\textit{Estimating Functions}
(see~\cite{nesterov2018lectures}):
\beq \label{LinearEstFunc}
\ba{rcl}
\phi_k(x) & \Def & \sum_{i = 1}^k a_i
\bigl[
f(x_i) + \la \nabla f(x_i), x - x_i \ra + \psi(x)
\bigr],
\qquad
\phi_k^{*} \;\; := \;\; \min\limits_x \phi_k(x),
\ea
\eeq
for the sequence of test points $\{ x_k : x_k \in \dom \psi \}_{k \geq 1}$ and positive
scaling coefficients $\{ a_k \}_{k \geq 1}$. We relate them with contracting coefficients,
as follows
\beq \label{RelationGammaAk}
\ba{rcl}
\gamma_k & := & \frac{a_{k + 1}}{A_{k + 1}},
\qquad
A_{k} \;\; \Def \;\; \sum_{i = 1}^k a_i.
\ea
\eeq

\BFT \label{Th:ContrDomNewtonConvex}
Let $A_k := k^3$, and consequently, $\gamma_k  :=  1 - \bigl(\frac{k}{k + 1}\bigr)^{3}
 =  \O\bigl(\frac{1}{k}\bigr)$.
Then for the sequence $\{ x_k \}_{k \geq 1}$ generated by
Algorithm~\ref{ContrDomNewton}, we have
\beq \label{Th1Guarantee}
\ba{rcl}
F(x_k) - F^{*} & \leq &
\ell_k \;\; \Def \;\; F(x_k) - \frac{\phi_k^{*}}{A_k}
\;\;
\leq \;\;
\O\bigl( \frac{H_{\nu} D^{2 + \nu}}{k^{1 + \nu}} \bigr).
\ea
\eeq
\EFT

For the case $\nu = 1$ (convex functions with Lipschitz
continuous Hessian), estimate~\eqref{Th1Guarantee} gives
the convergence rate of the order  $\O(\frac{1}{k^2})$.
This is the same rate, as we can achieve on this
functional class by Cubic Regularization of Newton
Method~\cite{nesterov2006cubic}. In accordance
to~\eqref{Th1Guarantee}, in order to obtain
$\varepsilon$-accuracy in functional residual, $F(x_K) -
F^{*} \leq \varepsilon$, it is enough to perform
\beq \label{ComplEstimate}
\ba{rcl}
K & = & \O\Bigl( \bigl(
\frac{H_{\nu} D^{2 + \nu}}{\varepsilon}
\bigr)^{1 / (1 + \nu) } \Bigr)
\ea
\eeq
iterations of Algorithm~\ref{ContrDomNewton}.
In~\cite{grapiglia2017regularized}, there were proposed
first \textit{universal} second-order methods (which do
not depend on parameters $\nu$ and $H_{\nu}$ of the
problem class), having complexity guarantees of the same
order~\eqref{ComplEstimate}. These methods are based on
Cubic regularization and an adaptive search for estimating
the regularization parameter at every iteration. It is
important that Algorithm~\ref{ContrDomNewton} is both
universal and affine-invariant. Additionally, convergence
result~\eqref{Th1Guarantee} provides us with a sequence
$\{ \ell_k \}_{k \geq 1}$ of computable \textit{accuracy
certificates}, which can be used as a stopping criterion
of the method.

Now, let us assume that the composite component is \textit{strongly convex}
with parameter $\mu > 0$. Thus, for all $x, y \in \dom \psi$
and $\psi'(x) \in \partial \psi(x)$, it holds
\beq \label{PsiStronglyConvex}
\ba{rcl}
\psi(y) & \geq & \psi(x) + \la \psi'(x), y - x \ra + \frac{\mu}{2}\|y - x\|^2.
\ea
\eeq
In this situation, we are able to improve convergence
estimate~\eqref{Th1Guarantee}, as follows. 

\newpage
\BFT
\label{Th:ContrDomNewtonStronglyConvex} Let $A_k := k^5$,
and consequently, $\gamma_k := 1 - \bigl( \frac{k}{k + 1}
\bigr)^{5}  =  \O\bigl(\frac{1}{k}\bigr)$. Then for the
sequence $\{ x_k \}_{k \geq 1}$ generated by
Algorithm~\ref{ContrDomNewton}, we have
\beq \label{Th2Guarantee1}
\ba{rcl}
F(x_k) - F^{*} & \leq &
\ell_k \;\; \leq \;\;
\O\Bigl( \frac{H_{\nu} D^{\nu}}{\mu} \cdot \frac{H_{\nu} D^{2 + \nu}}{ k^{2 + 2\nu}} \Bigr).
\ea
\eeq
Moreover, if the second-order \underline{condition number}
\beq \label{CondNumber}
\ba{rcl}
\omega_{\nu} & \Def & \Bigl[  \frac{H_{\nu} D^{\nu}}{(1 +
\nu) \mu}  \Bigr]^{\frac{1}{1 + \nu}}
\ea
\eeq
is known, then, defining $A_k := (1 +
\omega_{\nu}^{-1})^k, \; k \geq 1$, $A_0 := 0$, and $\gamma_k := \frac{1}{1 +
\omega_{\nu}}, \; k \geq 1$, $\gamma_0 := 1$, we obtain
the global \underline{linear rate} of convergence
\beq \label{Th2Guarantee2}
\ba{rcl}
F(x_k) - F^{*} & \leq & \ell_k \;\; \leq \;\; \exp\bigl(
-\frac{k - 1}{1 + \omega_{\nu}}  \bigr) \cdot
\frac{H_{\nu} D^{2 + \nu}}{1 + \nu}.
\ea
\eeq
\EFT According to the estimate~\eqref{Th2Guarantee2}, in
order to get $\varepsilon$-accuracy in function value, it
is enough to perform
$$
\ba{rcl}
K & = & \O\bigl((1 +  \omega_{\nu}) \cdot \log \frac{F(x_0) - F^{*}}{\varepsilon} \bigr)
\ea
$$
iterations of the method. Hence, condition number
$\omega_{\nu}$ plays the role of the main complexity
factor. This rate corresponds to that one of Cubically
Regularized Newton Method
(see~\cite{doikov2019local,doikov2019minimizing}).
At the same time, there exists a second variant of
Contracting-Domain Newton Method, where the next point is
defined by minimization of the full second-order model for
the smooth component augmented by the composite term over
the \textit{contracted domain} (this explains the names of
our methods).

\begin{algorithm}[H]
    \caption{Contracting-Domain Newton Method, II}
    \label{ContrDomNewton2}
    \begin{algorithmic}[1]
        \Require Choose $x_0 \in \dom \psi$.
        \Ensure {$k \geq 0$.}
        \State Pick up $\gamma_k \in (0, 1]$.
        \State Denote \vspace*{-7pt}
        $$
        \ba{lcl}
        S_k(y) & := &
        \begin{cases}
        \psi(y), & \; y \in \gamma_k \dom\psi + (1 - \gamma_k) x_k, \\
        +\infty, & \; \text{else}.
        \end{cases}
        \ea
        $$ \vspace*{-10pt}
        \State Compute \vspace*{-10pt}
        $$
        \ba{rcl}
        x_{k + 1} & \in & \Argmin\limits_{y} \Bigl\{
        \la \nabla f(x_k), y - x_k \ra
        \; + \; \frac{1}{2} \la \nabla^2 f(x_k)(y - x_k), y - x_k \ra
        \; + \; S_k(y)
        \Bigr\}.
        \ea
        $$
        \vspace*{-12pt}
    \end{algorithmic}
\end{algorithm}

Note, that Algorithm~\ref{ContrDomNewton} admits similar
representation as well.~\footnote{Indeed, it is enough to
take $S_k(y) := \gamma_k \psi( x_k + \frac{1}{\gamma_k}(y
- x_k) )$.} Both methods produce the same sequences of
points when $\psi(\cdot)$ is $\{0, +\infty\}$-indicator of
a convex set. Otherwise, they are different. 
Using the same contraction technique, 
it was shown in~\cite{nesterov2018complexity} 
that the classical Frank-Wolfe algorithm can be extended
onto the case of the composite optimization problems.
Additionally, the second-order 
\textit{Contracting Trust-Region method} was proposed, which has 
the same form as Algorithm~\ref{ContrDomNewton2}.
However,
its convergence rate was established only at the level
$\O(\frac{1}{k})$. Here, we improve its rate as follows.
\BFT \label{Th:ContrDomNewton2Convex} Let $A_k := k^3$ and
$\gamma_k  :=  1 - \bigl(\frac{k}{k + 1}\bigr)^{3} =
\O\bigl(\frac{1}{k}\bigr)$. Then for the sequence $\{ x_k
\}_{k \geq 1}$ generated by
Algorithm~\ref{ContrDomNewton2}, we have
\beq \label{Th3Guarantee}
\ba{rcl}
F(x_k) - F^{*} & \leq &
\ell_k \;\; \leq \;\;
\O\bigl( \frac{H_{\nu} D^{2 + \nu}}{k^{1 + \nu}} \bigr).
\ea
\eeq
\EFT This result is very similar to
Theorem~\ref{ContrDomNewton}. However, the first algorithm
can be accelerated on the class of strongly convex
functions (see
Theorem~\ref{Th:ContrDomNewtonStronglyConvex}). Thus, it
seems that it is more preferable.

Finally, let us consider an example, when the composite
component $\psi(\cdot)$ is an $\ell_p$-ball, as
in~\eqref{PsiBall}. Then, iterations of the method can be
represented as
\beq \label{TrustRegionStep}
\ba{rcl}
x_{k + 1} & \in & x_k +  \Argmin\limits_h \Bigl\{
\la \nabla f(x_k), h \ra
+ \frac{1}{2} \la \nabla^2 f(x_k)h, h \ra
\, : \,
\| x_k + \frac{1}{\gamma_k} h \|_p \leq \frac{D}{2}
\Bigr\}.
\ea
\eeq
In this form, it looks as a variant of Trust-Region
scheme. To solve the subproblem
in~\eqref{TrustRegionStep}, we can use Interior Point
Methods~(e.g. Chapter 5 in~\cite{nesterov2018lectures}).
See also~\cite{conn2000trust}, for techniques, developed
for Trust-Region schemes. Usually, complexity of this step
can be estimated as $\O(n^3)$ arithmetic operations, which
comes from the cost of computing a suitable factorization
for the Hessian matrix. Alternatively, Hessian-free
gradient methods can be applied, for computing an inexact
step (see~\cite{carmon2020first, carderera2020second}).

\section{Aggregated second-order models}
\label{Sec:Aggr}

In this section, we propose more advanced second-order
models, based on global lower
bound~\eqref{FuncLowerBound}. Using the same notation as
before, consider a sequence of test points $\{ x_k : x_k
\in \dom \psi \}_{k \geq 0}$ and sequences of coefficients
$\{ a_k \}_{k \geq 1}$, $\{ \gamma_k \}_{k \geq 0}$,
satisfying the relations~\eqref{RelationGammaAk}. Then, we
can introduce the following \textit{Quadratic Estimating
Functions} (compare with
definition~\eqref{LinearEstFunc}):
$$
\ba{rcl}
Q_k(x) & \Def & \sum_{i = 0}^{k - 1}
a_{i + 1} \Bigl[
f(x_i) + \la \nabla f(x_i), x - x_i \ra + \frac{\gamma_i}{2}\la \nabla^2 f(x_i)(x - x_i), x - x_i\ra
+ \psi(x)
\Bigr].
\ea
$$
By~\eqref{FuncLowerBound}, we have the main property of
Estimating Functions being satisfied. Namely, for all $x
\in \dom \psi$
\beq \label{QEstUpper}
\ba{rcl}
A_k F(x) & \refGE{FuncLowerBound} &
Q_k(x) - \sum_{i = 0}^{k - 1}
\frac{a_{i + 1} \gamma_i^{1 + \nu} H_{\nu} \|x - x_i\|^{2 + \nu} }{(1 + \nu)(2 + \nu)} \\
\\
& \refGE{DiamDef} &
Q_k(x) - \frac{H_{\nu} D^{2 + \nu}}{(1 + \nu)(2 + \nu)} \sum_{i = 0}^{k - 1} a_{i + 1}\gamma_i^{1 + \nu}
\;\; =: \;\; Q_k(x) - \frac{C_k}{2}.
\ea
\eeq
Therefore, if we would be able to guarantee for our test
points the relation
\beq \label{QEstLower}
\ba{rcl}
Q_k^{*} \;\; := \;\; \min\limits_x Q_k(x)  & \geq & A_k F(x_k) - \frac{C_k}{2},
\ea
\eeq
then we could immediately obtain the global convergence in
function value. Fortunately, relation \eqref{QEstLower}
can be achieved by simple iterations.

\begin{algorithm}[H]
    \caption{Aggregating Newton Method}
    \label{AggrNewton}
    \begin{algorithmic}[1]
        \Require Choose $x_0 \in \dom \psi$. Set $A_0 := 0$, $Q_0(x) \equiv 0$.
        \Ensure {$k \geq 0$.}
        \State Pick up $a_{k + 1} > 0$. Set $A_{k + 1} := A_k + a_{k + 1}$
        and $\gamma_k := \frac{a_{k + 1}}{A_{k + 1}}$.
        \State Update Estimating Function
        \Statex $
        Q_{k + 1}(x) \; \equiv \; Q_k(x)  \, + \,  a_{k + 1} \bigl[
        f(x_k) + \la \nabla f(x_k), x - x_k \ra
        + \frac{\gamma_k}{2} \la  \nabla^2 f(x_k)(x - x_k), x - x_k \ra + \psi(x) \bigr].
        $
        \vspace*{-6pt}
        \State Compute \vspace*{-8pt}
        $$
        \ba{rcl}
        v_{k + 1} & \in & \Argmin\limits_x Q_{k + 1}(x).
        \ea
        $$
        \vspace*{-12pt}
        \State Set $x_{k + 1} \; := \; x_k + \gamma_k(v_{k + 1} - x_k)$.
    \end{algorithmic}
\end{algorithm}
Clearly, the most complicated part of this process is
Step~3, which is computation of the minimum of Estimating
Function. However, the complexity of this step remains the
same, as that one for Contracting-Domain Newton Method. We
obtain the following convergence result. \BFT
\label{Th:AggrNewtonGuarantee} For the sequence $\{ x_k
\}_{k \geq 1}$ generated by Algorithm~\ref{AggrNewton},
relation \eqref{QEstLower} is satisfied. Consequently, for
the choice $A_k := k^3$, we obtain
\beq \label{AggrNewtonGuarantee}
\ba{rcl}
F(x_k) - F^{*} & \refLE{QEstUpper} &
F(x_k) - \frac{Q_k^{*}}{A_k} + \frac{C_k}{2A_k}
\;\; \refLE{QEstLower} \;\;
\frac{C_k}{A_k}
\;\; \leq \;\; \O\bigl( \frac{H_{\nu}D^{2 + \nu}}{k^{1 + \nu}} \bigr).
\ea
\eeq
\EFT
Now, for the accuracy certificate we have new expression
$\bar{\ell}_k := F(x_k) - \frac{Q_k^{*}}{A_k} + \frac{C_k}{2A_k}$.
The value of $Q_k^{*}$ is available within the method directly.
However, in order to compute $\bar{\ell}_k$ in practice,
some estimate for $C_k$ is required.
Note, that for the given choice of coefficients $A_k := k^3$,
we have $a_k = \O(k^2)$ and $\gamma_k = \O(\frac{1}{k})$.
Therefore, new information enters into the model
with increasing weights, which seems to be natural.

\section{Stochastic finite-sum minimization}
\label{Sec:Stoch}

In this section, we consider the case when the smooth part $f$ of the objective~\eqref{MainProblem}
is represented as a sum of
$M$ convex twice-differentiable components,
\beq \label{FiniteSum}
\ba{rcl}
f(x) & := & \frac{1}{M} \sum_{i = 1}^M f_i(x).
\ea
\eeq
This setting appears in many machine learning
applications, such as \textit{empirical risk
minimization}. Often, the number $M$ is very big. Thus, it
becomes expensive to evaluate the whole gradient or the
Hessian at every iteration. Hence, \textit{stochastic} or
\textit{incremental} methods are the methods of choice in
this situation. See~\cite{bertsekas2011incremental}
for a survey of first-order incremental methods.
The Newton-type Incremental Method with superlinear
local convergence was proposed in~\cite{rodomanov2016superlinearly}.
Local linear rate of stochastic Newton methods was studied in~\cite{kovalev2019stochastic}.
Global convergence of sub-sampled Newton schemes,
based on Damped iterations, and on Cubic regularization, was established 
in~\cite{roosta2016sub,kohler2017sub,tripuraneni2018stochastic}.

The basic idea of stochastic algorithms is to substitute the true
gradients and Hessians by some random unbiased estimators
$g_k$, and $H_k$, respectively, with $\E[g_k] = \nabla
f(x_k)$ and $\E[H_k] = \nabla^2 f(x_k)$.

First, let us consider the simplest estimation strategy.
At iteration $k$, we sample uniformly and independently
two subsets of indices $S_k^g, S_k^H \subseteq \{1, \dots,
M \}$. Their sizes are $m_k^g := |S_k^g|$ and $m_k^H :=
|S_k^H|$, which are possibly different. Then, in
Algorithm~\ref{ContrDomNewton}, we can use the following
random estimators:
\beq \label{RandomEstimators}
\ba{rcl}
g_k & := & \frac{1}{m_k^{g}} \sum_{i \in S_k^g} \nabla
f_i(x_k), \qquad H_k \;\; := \;\; \frac{1}{m_k^H} \sum_{i
\in S_k^H} \nabla^2 f_i(x_k).
\ea
\eeq
Let us present for this process a result on its global
convergence. Note that in this section, we use the standard Euclidean norm
for vectors and the corresponding induced spectral norm for matrices.

\BFT \label{Th:StochNewton}
Let each component $f_i(\cdot)$ be Lipschitz continuous on $\dom \psi$ with constant $L_0$,
and have Lipschitz continuous gradients and Hessians on $\dom \psi$
with constants $L_1$ and $L_2$, respectively.
Let $\gamma_k := 1 - \bigl( \frac{k}{k + 1} \bigr)^3 = \O\bigl(\frac{1}{k} \bigr)$. Set
\beq \label{BasicBatchSize}
\ba{rcl}
m_k^g & := & 1 / \gamma_k^4,
\qquad
m_k^H \;\; := \;\; 1 / \gamma_k^2.
\ea
\eeq
Then, for the iterations $\{ x_k \}_{k \geq 1}$ of
Algorithm~\eqref{ContrDomNewton}, based on
estimators~\eqref{RandomEstimators}, it holds
\beq \label{StochConv1}
\ba{rcl}
\E[F(x_k) - F^{*}] & \leq & 
\O \Bigl( \frac{ L_2 D^3 \, + \, L_1 D^2 (1 + \log(n)) \, + \,  L_0 D  }{k^2} \Bigr).
\ea
\eeq
\EFT Therefore, in order to solve our problem with
$\varepsilon$-accuracy in expectation,
 $\E[F(x_K) - F^{*}] \leq \varepsilon$, we need to perform
$K = \O\bigl( \frac{1}{\varepsilon^{1/2}} \bigr)$
iterations of the method. In this case, the total number
of gradient and Hessian samples are $\O\bigl(
\frac{1}{\varepsilon^{5/2}} \bigr)$ and $\O\bigl(
\frac{1}{\varepsilon^{3/2}} \bigr)$, respectively. It is
interesting that we need higher accuracy for estimating
the gradients, which results in a bigger batch size.

To improve this result, we incorporate a simple
\textit{variance reduction} strategy for the gradients.
This is a popular technique in stochastic convex optimization (see~\cite{schmidt2017minimizing,johnson2013accelerating,defazio2014saga,hazan2016variance,allen2017katyusha,nguyen2017sarah,gorbunov2019unified} and references therein).
At some iterations, we recompute the full gradient. However,
during the whole optimization process this happens
logarithmic number of times in total. Let us denote by
$\pi(k)$ the maximal \textit{power of two}, which is less than
or equal to $k$: $\pi(k) := 2^{\lfloor \log_2 k
\rfloor}$, for $k > 0$, and define $\pi(0) := 0$. The
entire scheme looks as follows.

\begin{algorithm}[H]
    \caption{Stochastic Variance-Reduced Contracting-Domain Newton}
    \label{StochContrDomNewton}
    \begin{algorithmic}[1]
        \Require Choose $x_0 \in \dom \psi$.
        \Ensure $k \geq 0$.
        \State Set anchor point $z_{k} := x_{\pi(k)}$.
        \State Sample random batch $S_k \subseteq \{1, \dots, M \}$ of size $m_k$.
        \State Compute variance-reduced stochastic gradient \vspace*{-5pt}
        $$
        \ba{rcl}
        g_k & := & \frac{1}{m_k} \sum_{i \in S_k}
        \bigl( \nabla f_i(x_k)  - \nabla f_i(z_k) + \nabla f(z_k) \bigr).
        \ea
        $$
        \vspace*{-15pt}
        \State Compute stochastic Hessian \vspace*{-8pt}
        $$
        \ba{rcl}
        H_k & := & \frac{1}{m_k} \sum_{i \in S_k} \nabla^2 f_i(x_k).
        \ea
        $$
        \vspace*{-15pt}
        \State Pick up $\gamma_k \in (0, 1]$.
        \State Perform the main step \vspace*{-5pt}
        $$
        \ba{rcl}
        x_{k + 1} & \in & \Argmin\limits_{y} \Bigl\{
        \la g_k, y - x_k \ra
        \; + \; \frac{1}{2} \la H_k(y - x_k), y - x_k \ra
        \; + \; \gamma_k \psi(x_k + \frac{1}{\gamma_k}(y - x_k))
        \Bigr\}.
        \ea
        $$
        \vspace*{-12pt}
    \end{algorithmic}
\end{algorithm}

Note that this is just Algorithm~\ref{ContrDomNewton} with random estimators $g_k$ and $H_k$
instead ot the true gradient and Hessian. The following global convergence result holds.

\BFT \label{Th:StochVRNewton}
Let each component $f_i(\cdot)$ have Lipschitz continuous gradients and Hessians
on $\dom \psi$ with constants $L_1$ and $L_2$, respectively.
Let $\gamma_k := 1 - \bigl( \frac{k}{k + 1} \bigr)^3 = \O(\frac{1}{k})$. Set batch size
\beq \label{VRBatchSize}
\ba{rcl}
m_k & := & 1 / \gamma_k^2.
\ea
\eeq
Then, for all iterations $\{ x_k \}_{k \geq 1}$ of
Algorithm~\ref{StochContrDomNewton}, we have
\beq \label{StochConv2}
\ba{rcl}
\E[ F(x_k) - F^{*} ] & \leq & 
\O \Bigl( \frac{ L_2 D^3 \, + \, L_1 D^2 (1 + \log(n)) \, + \, L_1^{1/2} D (F(x_0) - F^{*}) }{k^2} \Bigr).
\ea
\eeq
\EFT

It is thanks to the variance reduction that
we can use the same batch size for both estimators now.
To solve the problem with $\varepsilon$-accuracy in expectation,
we need $K = \O\bigl( \frac{1}{\varepsilon^{1/2}} \bigr)$
iterations of the method. And the total number 
of gradient and Hessian samples during these iterations is
$\O\bigl( \frac{1}{\varepsilon^{3/2}} \bigr)$.

\section{Experiments}
\label{Sec:Experiments}

\begin{figure}[h!]
	\centering
	\includegraphics[width=0.32\textwidth ]{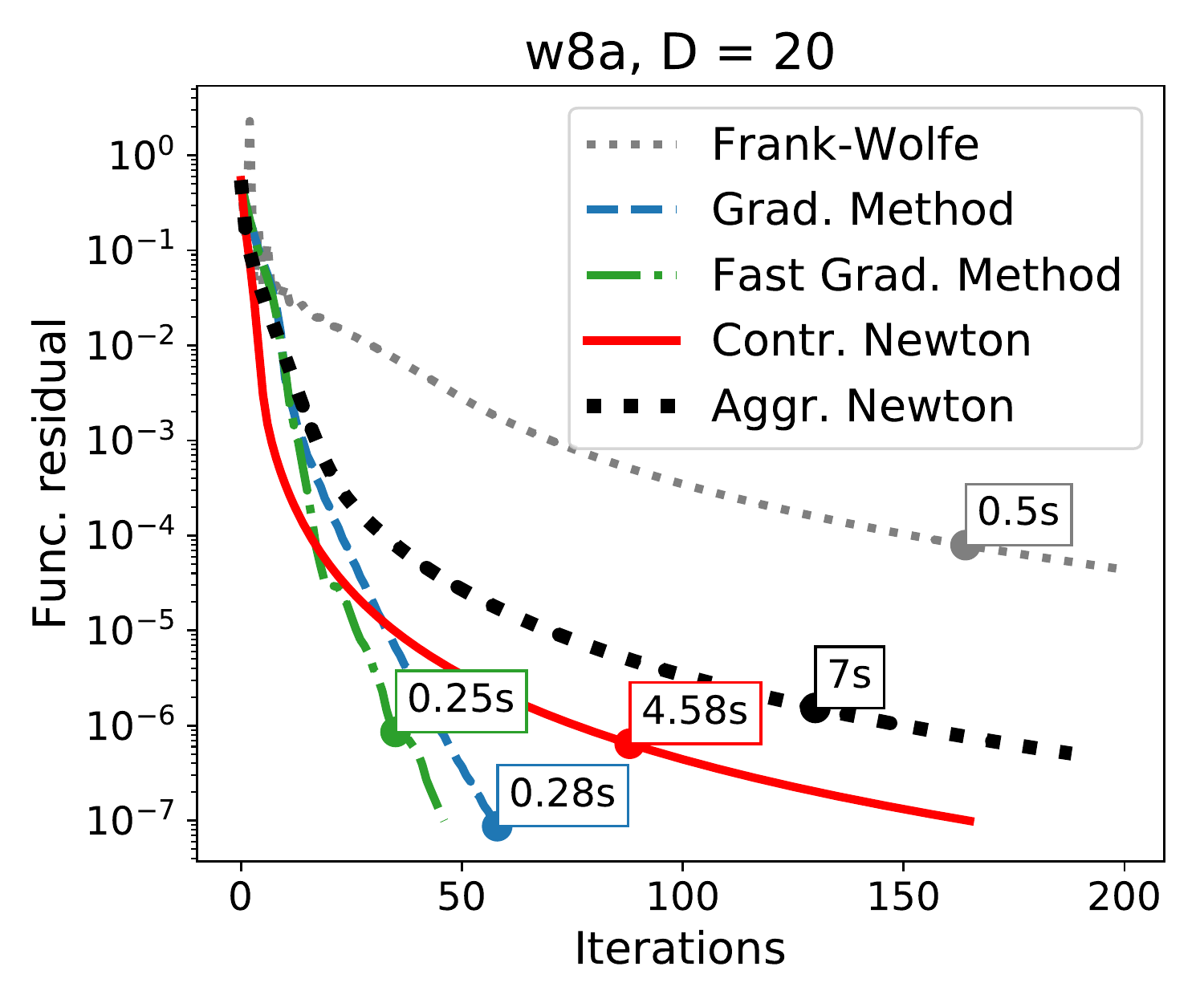}
	\includegraphics[width=0.32\textwidth ]{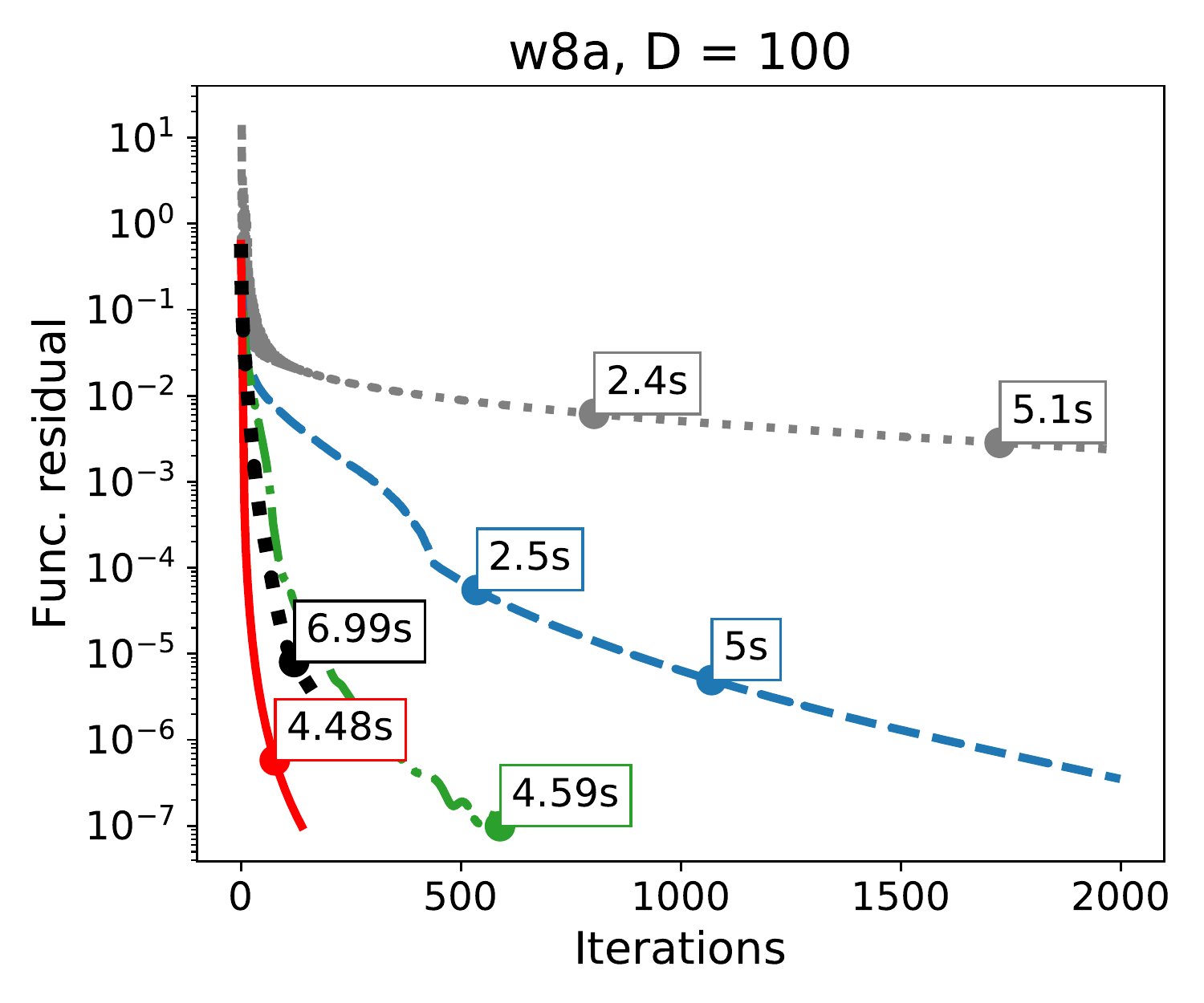}
	\includegraphics[width=0.32\textwidth ]{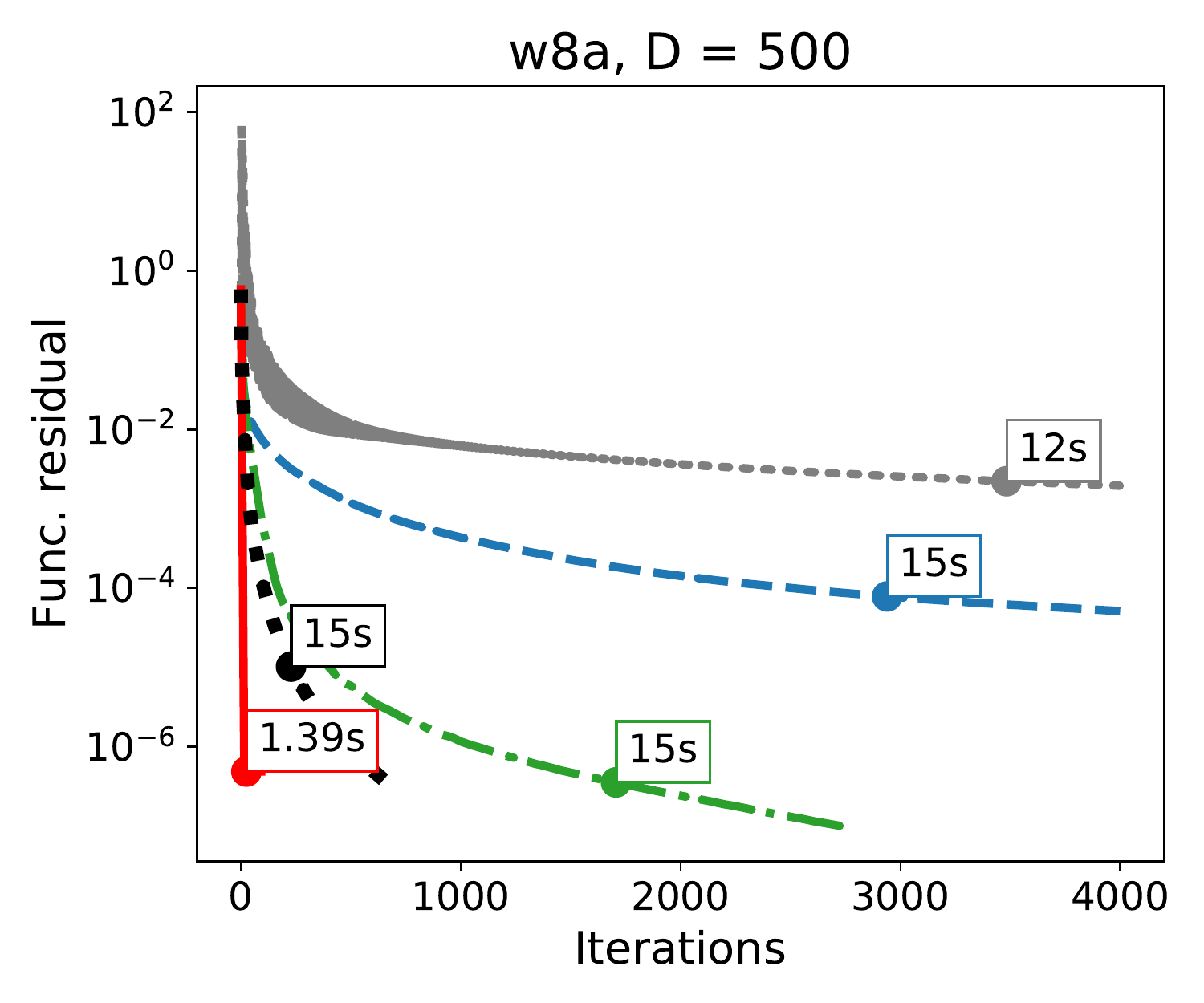}
	\caption{Training logistic regression, \textit{w8a} $(M = 49749, n = 300)$.  }
	\label{Fig:logreg_w8a}
\end{figure}

Let us demonstrate computational results for the problem
of training Logistic Regression model, regularized by $\ell_2$-ball constraints.
Thus, the smooth part of the objective has the finite-sum
representation~\eqref{FiniteSum}, each component is
$ f_i(x) := \log(1 + \exp(\la a_i, x \ra)) $.
The composite part is given by~\eqref{PsiBall}, with $p = 2$.
Diameter $D$ plays the role of regularization parameter, while
vectors $\{ a_i : a_i \in \R^n \}_{i = 1}^M$
are determined by the dataset\footnote{\url{https://www.csie.ntu.edu.tw/~cjlin/libsvmtools/datasets/}}.
First, we compare the performance of Contracting-Domain 
Newton Method
(Algorithm~\ref{ContrDomNewton}) 
and Aggregating Newton Method (Algorithm~\ref{AggrNewton})
 with first-order optimization schemes:
 Frank-Wolfe algorithm~\cite{frank1956algorithm},
 the classical Gradient Method, and the Fast Gradient Method~\cite{nesterov2013gradient}.
 For the latter two we use a line-search at each iteration, to estimate the Lipschitz constant. 
 The results are shown on Figure~\ref{Fig:logreg_w8a}.

 We see, that for bigger $D$, it becomes harder to solve the optimization problem.
 Second-order methods demonstrate good performance both in terms of 
 the iterations, and the total computational time.
 \footnote{Clock time was evaluated using the machine with Intel Core i5 CPU, 1.6GHz; 8 GB RAM.
 All methods were implemented in C++.
The source code can be found at \url{https://github.com/doikov/contracting-newton/}}
 
\begin{figure}[h!]
	\centering
	\includegraphics[width=0.32\textwidth ]{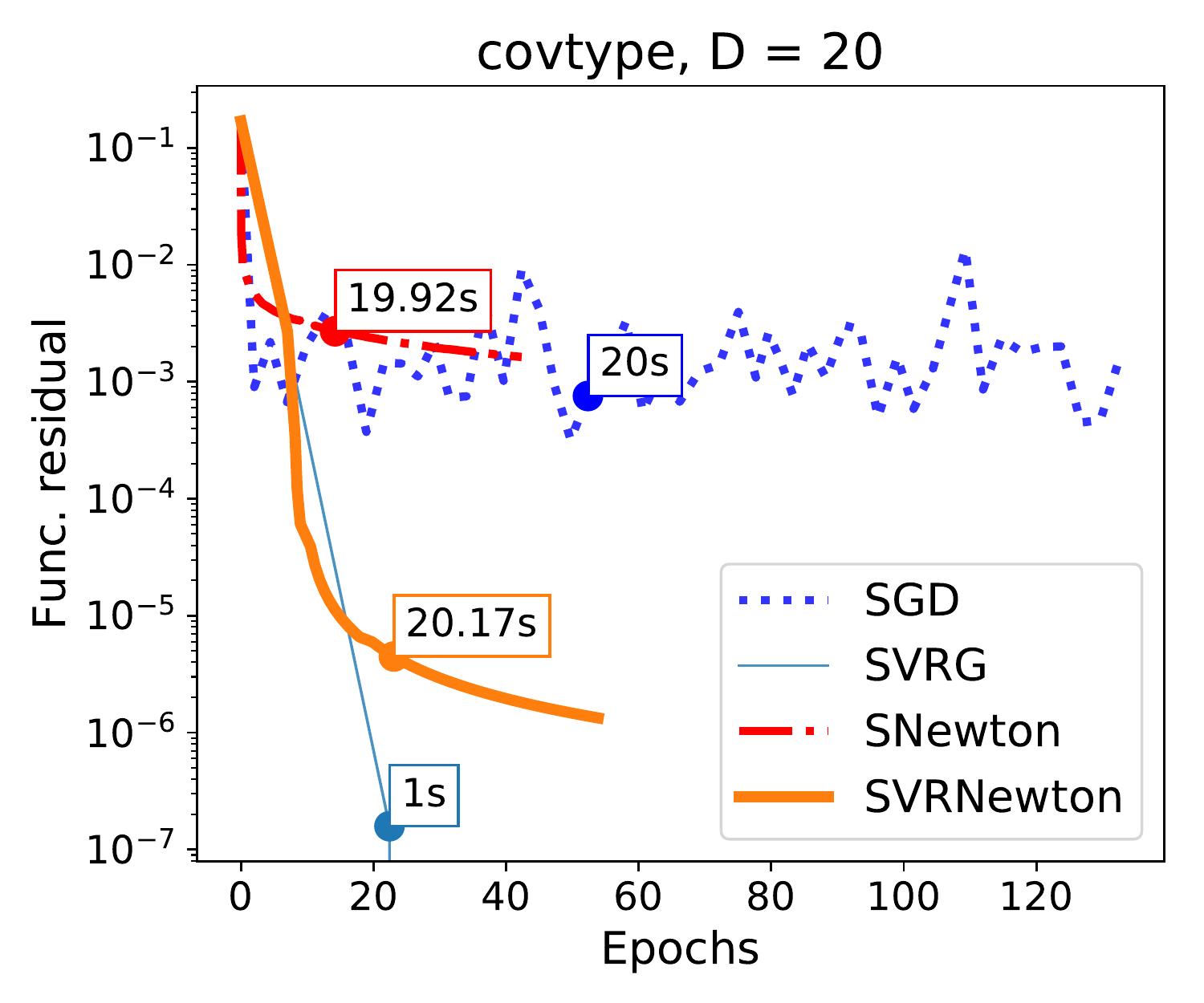}
	\includegraphics[width=0.32\textwidth ]{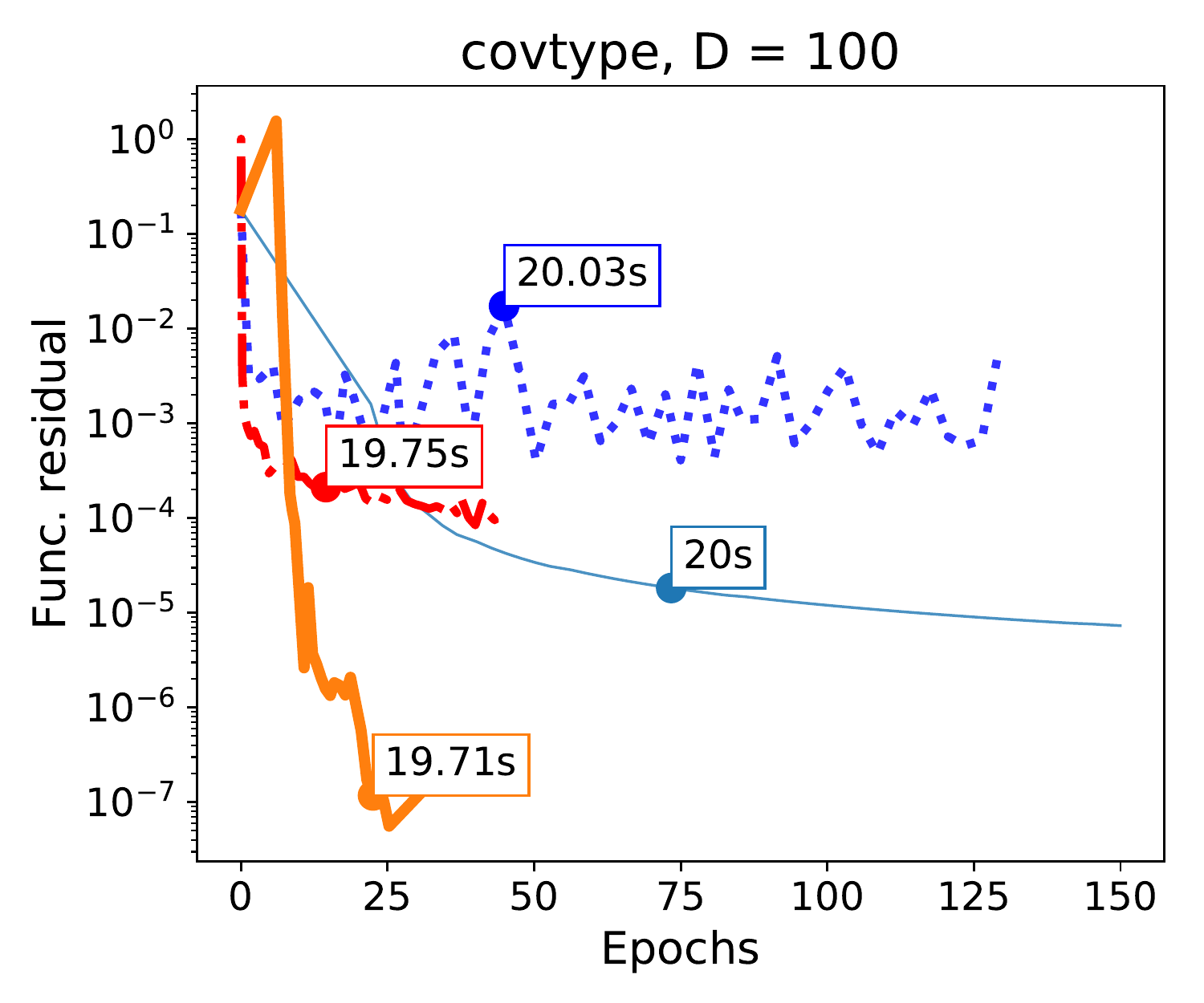}
	\includegraphics[width=0.32\textwidth ]{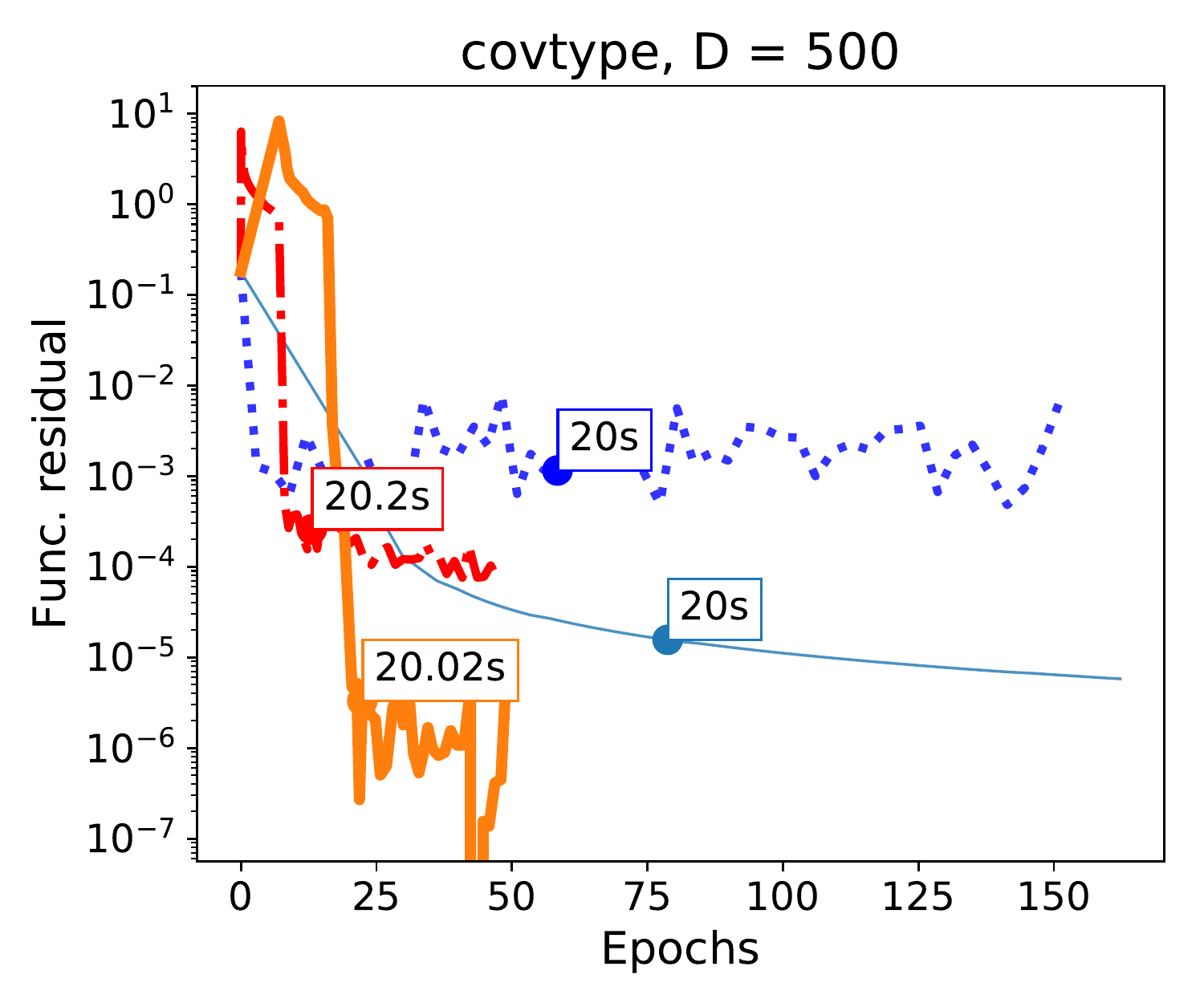}
	\caption{Stochastic methods for training logistic regression, 
		\textit{covtype} $(M = 581012, n = 54)$. }
	\label{Fig:logreg_stoch}
\end{figure}

In the next set of experiments, we compare the basic stochastic version of our method,
using estimators~\eqref{RandomEstimators} --- SNewton,
the method with the variance reduction
(Algorithm~\ref{StochContrDomNewton}) --- SVRNewton,
and first-order algorithms (with constant step-size, tuned for each problem): SGD and SVRG~\cite{johnson2013accelerating}.
We see (Figure~\ref{Fig:logreg_stoch}) that using the variance reduction strategy
significantly improve the convergence for both first-order and second-order 
stochastic optimization methods. 

According to these graphs, our second-order algorithms can be 
more efficient when solving ill-conditioned problems,
producing the better solution within a given computational time.
See also Section~\ref{SectionEE} in Appendix 
for extra experiments.

\section{Discussion}
\label{Sec:Discussion}

Let us discuss complexity estimates,
which we established in our work.  
For the basic versions of our method
we have the global convergence in the functional residual
of the form
$$
\ba{rcl}
F(x_k) - F^{*} & \leq & \O\bigl(\frac{H_{\nu} D^{2 + \nu}}{k^{1 + \nu}}\bigr).
\ea
$$ 
Note that the complexity parameter $H_{\nu}$ depends
only on the variation of the Hessian (in arbitrary norm).
It can be much smaller than the maximal eigenvalue 
of the Hessian, which typically appears in the rates
of first-order methods. 
It is important that our algorithms 
are free from using the norms or any other particular
parameters of the problem class.

At the same time, the arithmetic complexity of one step of our methods
for simple sets can be estimated as the sum of the cost of computing the Hessian, 
and $\O(n^3)$ additional operations (to compute a suitable factorization of the matrix).
For example, the cost of computing the gradient of Logistic Regression 
is $\O(Mn)$, and the Hessian is $\O(M n^2)$, where $M$ is the dataset size.
Hence, it is preferable to use our algorithms with exact steps
in the situation when $M$ is much bigger than $n$.

\begin{ack}
The research results of this paper were obtained in the framework of ERC Advanced Grant 788368.
\end{ack}

\bibliographystyle{plain}
\bibliography{bibliography}


\appendix
\clearpage 
\onecolumn 
\newcommand{\toptitlebar}{
	\hrule height 1pt
	\vskip 0.25in
	\vskip -\parskip%
}
\newcommand{\bottomtitlebar}{
	\vskip 0.29in
	\vskip -\parskip
	\hrule height 1pt
	\vskip 0.09in%
}
\vbox{%
	\hsize\textwidth
	\linewidth\hsize
	\vskip 0.1in
	\toptitlebar
	\centering
	{\LARGE\bf Appendix \par}
	\bottomtitlebar
	\vskip 0.2in
}

\section{Proof of Lemma~\ref{LemmaLB}}

First, let us note that inequality~\eqref{GradHolderBound} follows from the
following simple observation,
using Newton-Leibniz formula and H\"older continuity of the Hessian,
for all $x, y \in \dom \psi$
$$
\ba{rcl}
\| \nabla f(y) - \nabla f(x) - \nabla^2 f(x)(y - x) \|_{*}
& = &
\| \int\limits_{0}^1 ( \nabla^2 f(x + \tau(y - x)) - \nabla^2 f(x) ) (y - x) d\tau \|_{*} \\
\\
& \overset{\eqref{HolderHessian}}{\leq} &
\frac{H_{\nu} \|y - x\|^{1 + \nu}}{1 + \nu}.
\ea
$$
We are ready to prove the lemma.

\begin{lemma}
For all $x, y \in \dom \psi$ and $t \in [0, 1]$, it holds
$$
\ba{rcl}
f(y) & \geq & f(x) + \la \nabla f(x), y - x \ra
+ \frac{t}{2} \la \nabla^2 f(x)(y - x), y - x \ra
- \frac{t^{1 + \nu} H_{\nu} \|y - x\|^{2 + \nu}}{(1 + \nu)(2 + \nu)}.
\ea
$$
\end{lemma}
\proof
Let us prove the following bound,
for all $x, y \in \dom \psi$ and $t \in [0, 1]$
\beq \label{GradLowerBound}
\ba{rcl}
\la \nabla f(y) - \nabla f(x), y - x \ra 
& \geq & 
t \la \nabla^2 f(x)(y - x), y - x \ra
- \frac{t^{1 + \nu} H_{\nu} \|y - x\|^{2 + \nu}}{1 + \nu}.
\ea
\eeq
For $t = 1$ it follows from~\eqref{GradHolderBound}. 
Therefore,
we may assume that $t < 1$. 
Let us take $z_t := x + t(y - x)$. 
Then, by convexity of $f$, we have
$$
\ba{rcl}
\la \nabla f(y), y - x \ra 
& = & 
\frac{1}{1 - t} \la \nabla f(y), y - z_t \ra \\
\\
& \geq & \frac{1}{1 - t} \la \nabla f(z_t), y - z_t \ra
\;\; = \;\;
\la \nabla f(z_t), y - x \ra.
\ea
$$
Now, from H\"older continuity of the Hessian, we get
$$
\ba{rcl}
\la \nabla f(z_t), y - x \ra
& \refGE{GradHolderBound} &
\la \nabla f(x), y - x \ra 
+ \la \nabla^2 f(x)(z_t - x), y - x \ra
- \frac{H_{\nu} \|z_t - x\|^{1 + \nu} \|y - x\|}{1 + \nu} \\
\\
& = & 
\la \nabla f(x), y - x \ra + t \la \nabla^2 f(x)(y - x), y - x \ra
- \frac{t^{1 + \nu} H_{\nu }\|y - x\|^{2 + \nu}}{1 + \nu}.
\ea
$$
Thus we prove~\eqref{GradLowerBound}.
Then, the claim of the lemma can be obtained by
simple integration:
$$
\ba{rcl}
f(y) - f(x) - \la \nabla f(x), y - x \ra
& = &
\int\limits_0^1 
\la \nabla f(z_{\tau}) - \nabla f(x), y - x \ra d\tau \\
\\
& \refGE{GradLowerBound} &
\int\limits_0^1
t \tau \la \nabla^2 f(x)(y - x), y - x \ra
- \frac{(t \tau)^{1 + \nu} H_{\nu} \|y - x\|^{2 + \nu}}{1 + \nu} d\tau \\
\\
& = &
\frac{t}{2}\la \nabla^2 f(x)(y - x), y - x \ra
- \frac{t^{1 + \nu} H_{\nu} \|y - x\|^{2 + \nu}}{(1 + \nu)(2 + \nu)}.
\ea
$$
\qed

\section{Convergence of Contracting-Domain Newton Method}

In this section, we prove the global convergence of
Algorithms~\ref{ContrDomNewton} and~\ref{ContrDomNewton2}. 
We use the same notation as in the main part. There is 
a sequence of controlling coefficients~$\{ a_k \}_{k \geq 1}$~(see relations~\eqref{RelationGammaAk}),
and a sequence of linear Estimating Functions $\{ \phi_k(x) \}_{k \geq 0}$. We denote by $\mu \geq 0$ the constant of strong convexity of $\psi(\cdot)$. 
We allow $\mu = 0$ in the following auxiliary lemma,
in order to cover both the general convex and the strongly convex cases.

\setcounter{lemma}{1}
\BL \label{Lemma:ContrDomNewtonConv}

For the sequences $\{ x_k \}_{k \geq 1}$ and $\{ v_k \}_{k \geq 1}$,
produced by Algorithm~\ref{ContrDomNewton}, we have
\beq \label{ContrDomNewtonConv1}
\ba{rcl}
A_k F(x_k)
& \leq & \phi_k(x) \; + \; B_k(x), \qquad x \in \dom \psi,
\ea
\eeq
with
\beq \label{BkDef}
\ba{rcl}
B_k(x) & \equiv &
\sum\limits_{i = 1}^k \Bigl[
\frac{H_{\nu} a_i^{2 + \nu}  \|x - v_i\| \cdot \|x_{i - 1} - v_i\|^{1 + \nu} }{(1 + \nu) A_i^{1 + \nu}}
- \frac{\mu a_i  \|x - v_i\|^2}{2}
- \frac{\mu a_i A_{i - 1}  \|x_{i - 1} - v_i\|^2 }{2A_i}
\Bigr].
\ea
\eeq
\EL
\proof
Let us prove~\eqref{ContrDomNewtonConv1} by induction.
It obviously holds for $k = 0$, since $A_0 := 0$, 
$\phi_0(x) \equiv 0$, and $B_0(x) \equiv  0$
by definition.
Assume that it holds for the current $k \geq 0$, and consider the next iterate.
Stationary condition for the method step is
\beq \label{StatCondition}
\ba{rcl}
\la \nabla f(x_k) + \nabla^2 f(x_{k})(x_{k + 1} - x_k), x - v_{k + 1} \ra
+ \psi(x)
& \geq & \psi(v_{k + 1}) + \frac{\mu}{2}\|x - v_{k + 1}\|^2,
\ea
\eeq
for all $x \in \dom \psi$.  Then, we have
$$
\ba{rcl}
\phi_{k + 1}(x)
& \equiv &
a_{k + 1} \bigl[  f(x_{k + 1}) + \la \nabla f(x_{k + 1}), x - x_{k + 1} \ra + \psi(x) \bigr]
+ \phi_k(x) \\
\\
& \refGE{ContrDomNewtonConv1} &
a_{k + 1} \bigl[  f(x_{k + 1}) + \la \nabla f(x_{k + 1}), x - x_{k + 1} \ra + \psi(x)  \bigr]
+ A_k F(x_k) - B_k(x) \\
\\
& \overset{(*)}{\geq} & 
A_{k + 1} \bigl[ f(x_{k + 1}) + \la \nabla f(x_{k + 1}), 
\frac{a_{k + 1} x + A_k x_k}{A_{k + 1}} - x_{k + 1} \ra  \bigr]
+ a_{k + 1} \psi(x) \\
\\
& \quad & \qquad + \; A_k \psi(x_k) - B_k(x) \\
\\
& = &
A_{k + 1} f(x_{k + 1}) 
+ a_{k + 1}\la \nabla f(x_{k + 1}), x - v_{k + 1} \ra
+ a_{k + 1} \psi(x) \\
\\
& \quad & \qquad + \; A_k \psi(x_k) - B_k(x) \\
\\
& = &
A_{k + 1} f(x_{k + 1})
+ a_{k + 1} \bigl[ 
\la \nabla f(x_k) + \nabla^2 f(x_k)(x_{k + 1} - x_k), x - v_{k + 1} \ra
+ \psi(x)
\bigr] \\
\\
& \quad & \qquad + \;
a_{k + 1} \la \nabla f(x_{k + 1}) - \nabla f(x_k) - \nabla^2 f(x_k)(x_{k + 1} - x_k),
x - v_{k + 1} \ra  \\
\\
& \quad & \qquad + \; 
A_k \psi(x_k) - B_k(x) \\
\\
& \overset{\eqref{StatCondition},\eqref{GradHolderBound}}{\geq} &
A_{k + 1} f(x_{k + 1})
+ a_{k + 1} \bigl[ \psi(v_{k + 1}) + \frac{\mu}{2}\|x - v_{k + 1}\|^2 \bigr] \\
\\
& \quad & \qquad
- \; 
\frac{H_{\nu}  a_{k + 1}^{2 + \nu}  \|x - v_{k + 1}\| 
\cdot \|v_{k + 1} - x_k\|^{1 + \nu}  }{(1 + \nu) A_{k + 1}^{1 + \nu}}
+ A_k \psi(x_k) - B_k(x) \\
\\
& \overset{(**)}{\geq} &
A_{k + 1} F(x_{k + 1})
+ \frac{\mu a_{k + 1} \|x - v_{k + 1}\|^2}{2} 
+ \frac{\mu a_{k + 1} A_k }{2 A_{k + 1}} \|x_k - v_{k + 1}\|^2  \\
\\
& \quad & \qquad
- \; 
\frac{H_{\nu}  a_{k + 1}^{2 + \nu}  \|x - v_{k + 1}\| 
	\cdot \|v_{k + 1} - x_k\|^{1 + \nu}  }{(1 + \nu) A_{k + 1}^{1 + \nu}}
+ A_k \psi(x_k) - B_k(x) \\
\\
& \equiv &
A_{k + 1} F(x_{k + 1}) - B_{k + 1}(x),
\ea
$$
where $(*)$ and $(**)$ stand for convexity of $f$, and strong convexity of $\psi$,
correspondingly. Thus we have~\eqref{ContrDomNewtonConv1} established
for all $k \geq 0$.
\qed

\subsection{Proof of Theorem~\ref{Th:ContrDomNewtonConvex}}
\BT
Let $A_k := k^3$, and consequently, $\gamma_k  :=  1 - \bigl(\frac{k}{k + 1}\bigr)^{3}
=  \O\bigl(\frac{1}{k}\bigr)$.
Then for the sequence $\{ x_k \}_{k \geq 1}$ generated by
Algorithm~\ref{ContrDomNewton}, we have
$$
\ba{rcl}
F(x_k) - F^{*} & \leq &
\ell_k \;\; \Def \;\; F(x_k) - \frac{\phi_k^{*}}{A_k}
\;\;
\leq \;\;
\O\bigl( \frac{H_{\nu} D^{2 + \nu}}{k^{1 + \nu}} \bigr).
\ea
$$
\ET
\proof
First, by convexity of $f$ we have, for all $x \in \dom \psi$
$$
\ba{rcl}
\phi_k(x) & \leq & A_k F(x).
\ea
$$
Therefore, for the solution $x^{*}$ of our problem: $F^{*} = F(x^{*})$, it holds
$$
\ba{rcl}
F(x_k) - F^{*} & \leq & 
F(x_k) - \frac{\phi_k(x^{*})}{A_k}
\;\;
\leq 
\;\; \ell_k \;\; \Def \;\; F(x_k) - \frac{\phi_k^{*}}{A_k}, 
\ea
$$
and this is the first part of~\eqref{Th1Guarantee}.

At the same time, by Lemma~\ref{Lemma:ContrDomNewtonConv},
and using boundness of the domain, we have
$$
\ba{rcl}
\phi_k^{*} & := & \min\limits_{x \in \dom \psi} \Bigl\{ \phi_k(x) \Bigr\}
\;\;\; \refGE{ContrDomNewtonConv1} \;\;
\min\limits_{x \in \dom \psi} \Bigl\{  A_k F(x_k) - B_k(x)   \Bigr\} \\
\\
& \geq & 
A_k F(x_k) - \frac{H_{\nu} D^{2 + \nu} }{1 + \nu} 
\sum\limits_{i = 1}^k \frac{a_i^{2 + \nu}}{A_i^{1 + \nu}}
\ea
$$
Therefore, for the choice $A_k := k^3$, we finally obtain
$$
\ba{rcl}
\ell_k & \leq & \frac{H_{\nu} D^{2 + \nu}}{(1 + \nu) A_k}
\sum\limits_{i = 1}^k \frac{a_i^{2 + \nu}}{A_i^{1 + \nu}}
\;\; = \;\; 
\frac{H_{\nu} D^{2 + \nu}}{(1 + \nu) k^3}
\sum\limits_{i = 1}^k \frac{ (i^3 - (i - 1)^3)^{2 + \nu} }{i^{3(1 + \nu)}}
\\
\\
& \leq &
\frac{H_{\nu} D^{2 + \nu}}{(1 + \nu) k^3}
\sum\limits_{i = 1}^k \frac{ 3^{2 + \nu} i^{2(2 + \nu)} }{i^{3(1 + \nu)}}
\;\; = \;\;  
\frac{3^{2 + \nu} H_{\nu} D^{2 + \nu}}{(1 + \nu) k^3}
\sum\limits_{i = 1}^k i^{1 - \nu}  \\
\\
& = & \O\bigl( \frac{H_{\nu} D^{2 + \nu}}{k^{1 + \nu}} \bigr).
\ea
$$

\qed

\subsection{Proof of Theorem~\ref{Th:ContrDomNewtonStronglyConvex}}
\BT Let $A_k := k^5$,
and consequently, $\gamma_k := 1 - \bigl( \frac{k}{k + 1}
\bigr)^{5}  =  \O\bigl(\frac{1}{k}\bigr)$. Then for the
sequence $\{ x_k \}_{k \geq 1}$ generated by
Algorithm~\ref{ContrDomNewton}, we have
$$
\ba{rcl}
F(x_k) - F^{*} & \leq &
\ell_k \;\; \leq \;\;
\O\Bigl( \frac{H_{\nu} D^{\nu}}{\mu} \cdot \frac{H_{\nu} D^{2 + \nu}}{ k^{2 + 2\nu}} \Bigr).
\ea
$$
Moreover, if the second-order \underline{condition number}
$$
\ba{rcl}
\omega_{\nu} & \Def & \Bigl[  \frac{H_{\nu} D^{\nu}}{(1 +
	\nu) \mu}  \Bigr]^{\frac{1}{1 + \nu}}
\ea
$$
is known, then, defining $A_k := (1 +
\omega_{\nu}^{-1})^k, \; k \geq 1$, $A_0 := 0$, and $\gamma_k := \frac{1}{1 +
	\omega_{\nu}}, \; k \geq 1$, $\gamma_0 := 1$, we obtain
the global \underline{linear rate} of convergence
$$
\ba{rcl}
F(x_k) - F^{*} & \leq & \ell_k \;\; \leq \;\; \exp\bigl(
-\frac{k - 1}{1 + \omega_{\nu}}  \bigr) \cdot
\frac{H_{\nu} D^{2 + \nu}}{1 + \nu}.
\ea
$$
\ET
\proof
Starting from the same reasoning, as in the proof of
Theorem~\ref{Th:ContrDomNewtonConvex}, we get
$$
\ba{rcl}
F(x_k) - F^{*} & \leq & \ell_k \;\; \Def \;\; F(x_k) - \frac{\phi_k^{*}}{A_k}.
\ea
$$
Let us denote by $u_k$ the minimum of the Estimating Function $\phi_k$.
Thus,
$$
\ba{rcl}
\ell_k & = & F(x_k) - \frac{\phi_k(u_k)}{A_k}
\;\; \refLE{ContrDomNewtonConv1} \;\; \frac{1}{A_k}B_k(u_k)
\;\; \equiv \;\; \frac{1}{A_k}  \sum\limits_{i = 1}^k B_k^{(i)},
\ea
$$
with 
\beq \label{BkiBound}
\ba{rcl}
B_k^{(i)} & \Def &
a_i \Bigl[  
\frac{H_{\nu} a_i^{1 + \nu} \|u_k - v_i \| \cdot \|x_{i - 1} - v_i \|^{1 + \nu}}{(1 + \nu) A_i^{1 + \nu}}
- \frac{\mu \| u_k - v_i \|^2 }{2}
\Bigr] - \frac{\mu a_i A_{i - 1} \| x_{i - 1} - v_i \|^2}{2 A_i} \\
\\
& \leq & 
a_i \max\limits_{t \geq 0} \Bigl\{   
\frac{H_{\nu} a_i^{1 + \nu} \|x_{i - 1} - v_i\|^{1 + \nu} t}{(1 + \nu) A_i^{1 + \nu}}
- \frac{\mu t^2}{2}
\Bigr\}
- \frac{\mu a_i A_{i - 1} \| x_{i - 1} - v_i \|^2}{2 A_i} \\
\\
& = & 
\frac{a_i}{2 \mu} \Bigl( 
\frac{H_{\nu} a_i^{1 + \nu} \|x_{i - 1} - v_i\|^{1 + \nu}}{(1 + \nu) A_i^{1 + \nu}}
\Bigr)^2
- \frac{\mu a_i A_{i - 1} \| x_{i - 1} - v_i \|^2}{2 A_i}.
\ea
\eeq
Therefore, for the choice $A_k := k^5$, we have
$$
\ba{rcl}
\ell_k & \leq & 
\frac{1}{A_k} \sum\limits_{i = 1}^k
\frac{a_i}{2\mu} \Bigl(  \frac{H_{\nu} a_i^{1 + \nu} 
	\|x_{i - 1} - v_i\|^{1 + \nu}}{(1 + \nu) A_i^{1 + \nu}}  \Bigr)^2
\;\; \leq \;\;
\frac{H_{\nu}^2 D^{2(1 + \nu)}}{2 \mu (1 + \nu)^2 A_k}
\sum\limits_{i = 1}^k \frac{a_i^{2(1 + \nu) + 1}}{A_i^{2(1 + \nu)}} \\
\\
& = &
\frac{H_{\nu}^2 D^{2(1 + \nu)}}{2 \mu (1 + \nu)^2 k^5}
\sum\limits_{i = 1}^k \frac{ (i^5 - (i - 1)^5)^{2(1 + \nu) + 1} }{i^{10(1 + \nu)}}
\;\; \leq \;\;
\frac{5^{2(1 + \nu) + 1} H_{\nu}^2 D^{2(1 + \nu)}}{2 \mu (1 + \nu)^2 k^5}
\sum\limits_{i = 1}^k i^{2 - 2\nu} \\
\\
& = & 
\O\bigl( \frac{H_{\nu} D^{\nu}}{\mu} \cdot \frac{H_{\nu} D^{2+\nu}}{k^{2 + 2\nu}} \bigr).
\ea
$$
Thus we have justified~\eqref{Th2Guarantee1}. 
To obtain the linear rate~\eqref{Th2Guarantee2}, we set
$$
\ba{rcl}
A_k & := & (1 + \omega_{\nu}^{-1} )^{k},
\qquad k \geq 1,
\ea
$$
and $A_0 := 0$. So, $a_1 = A_1$ and
$$
\ba{rcl}
a_i & = & A_i - A_{i - 1} \;\; = \;\; \omega_{\nu}^{-1}A_{i - 1},
\qquad i \geq 2.
\ea
$$
Therefore, for the values $\{ B_k^{(i)} \}_{i = 1}^k$, we have
$$
\ba{rcl}
B_k^{(1)} & \leq & a_1 \frac{H_{\nu}D^{2 + \nu}}{1 + \nu}
\;\; = \;\; A_1 \frac{H_{\nu} D^{2 + \nu}}{1 + \nu},
\ea
$$
and
$$
\ba{rcl}
B_k^{(i)} & \refLE{BkiBound} & 
\frac{ H_{\nu}^{2} D^{2\nu} \|x_{i - 1} - v_i\|^2 a_i^{3 + 2\nu}  }{2 \mu (1 + \nu)^2 A_i^{2+ 2\nu}} 
- \frac{\mu a_i A_{i - 1} \|x_{i - 1} - v_i\|^2}{2 A_i} \\
\\
& = &
\frac{ \mu a_i A_{i - 1} \|x_{i - 1} - v_i\|^2 }{2 A_i}
\Bigl(  \Bigl[ \frac{H_{\nu} D^{\nu}}{(1 + \nu) \mu} \Bigr]^2 
\frac{ a_i^{2 + 2\nu} }{A_i^{1 + 2 \nu} A_{i - 1} }  - 1 \Bigr) \\
\\
& \leq &
\frac{ \mu a_i A_{i - 1} \|x_{i - 1} - v_i\|^2 }{2 A_i}
\Bigl(  \Bigl[ \frac{H_{\nu} D^{\nu}}{(1 + \nu) \mu} \Bigr]^2 
\Bigl[ \frac{ a_i }{ A_{i - 1} } \Bigr]^{2(1 + \nu)}  - 1 \Bigr) \\
\\
& = & 0, \qquad 2 \leq i \leq k,
\ea
$$
since by our choice
$$
\ba{rcl}
\frac{a_i}{A_{i - 1}} & = & \omega_{\nu}^{-1}
\;\; \refEQ{CondNumber} \;\;
\Bigl[ \frac{(1 + \nu) \mu}{H_{\nu} D^{\nu}}  \Bigr]^{\frac{1}{1 + \nu}}.
\ea
$$
Finally, we obtain
$$
\ba{rcl}
\ell_k & \leq & \frac{1}{A_k} B_k^{(1)} \;\; \leq \;\;
\frac{A_1}{A_k} \cdot \frac{H_{\nu} D^{2 + \nu}}{1 + \nu}
\;\; = \;\; \frac{1}{(1 + \omega_{\nu}^{-1})^{k - 1}} \cdot
\frac{H_{\nu} D^{2 + \nu}}{1 + \nu} \\
\\
& \leq & 
\exp\bigl( - \frac{k - 1}{1 + \omega_{\nu}} \bigr)
\cdot \frac{H_{\nu} D^{2 + \nu}}{1 + \nu}.
\ea
$$

\qed

\subsection{Proof of Theorem~\ref{Th:ContrDomNewton2Convex}}
\BT Let $A_k := k^3$ and
$\gamma_k  :=  1 - \bigl(\frac{k}{k + 1}\bigr)^{3} =
\O\bigl(\frac{1}{k}\bigr)$. Then for the sequence $\{ x_k
\}_{k \geq 1}$ generated by
Algorithm~\ref{ContrDomNewton2}, we have
$$
\ba{rcl}
F(x_k) - F^{*} & \leq &
\ell_k \;\; \leq \;\;
\O\bigl( \frac{H_{\nu} D^{2 + \nu}}{k^{1 + \nu}} \bigr).
\ea
$$
\ET
\proof
The proof is very similar to that one for Algorithm~\ref{ContrDomNewton}.
First, stationary condition for one iteration of 
Algorithm~\ref{ContrDomNewton2} is
\beq \label{Alg2StatCondition}
\ba{cl}
& \la \nabla f(x_k) + \nabla^2 f(x_k)(x_{k + 1} - x_k), x - v_{k + 1} \ra
+ \frac{1}{\gamma_k} \psi\bigl( \gamma_k x +(1 - \gamma_k) x_k \bigr) \\
\\
& \quad  \geq \quad 
\frac{1}{\gamma_k} \psi(x_{k + 1}),
\ea
\eeq
for all $x \in \dom \psi$ and $k \geq 0$ (compare with~\eqref{StatCondition}),
where
$$
\ba{rcl}
v_{k + 1} & := & x_k + \frac{1}{\gamma_k}(x_{k + 1} - x_k)
\;\; \in \;\; \dom \psi.
\ea
$$

Now, let us prove by induction the following bound
\beq \label{Th3Induction}
\ba{rcl}
\phi_k(x) & \geq & A_k F(x_k) - {B}_k, \qquad x \in \dom \psi,
\ea
\eeq
with 
${B}_k := \frac{H_{\nu} D^{2 + \nu}}{1 + \nu} 
\sum_{i = 1}^k \frac{a_i^{2 + \nu}}{A_i^{1 + \nu}}$.
It obviously holds for $k = 0$, since both sides are zero.
Assume that it holds for the current $k \geq 0$. 
Then, we have for the next iterate
$$
\ba{rcl}
\phi_{k + 1}(x) & \equiv & 
a_{k + 1} \bigl[ f(x_{k + 1}) + \la \nabla f(x_{k + 1}), x - x_{k + 1} \ra + \psi(x) \bigr]
+ \phi_k(x) \\
\\
& \refGE{Th3Induction} &
a_{k + 1} \bigl[ f(x_{k + 1}) + \la \nabla f(x_{k + 1}), x - x_{k + 1} \ra + \psi(x) \bigr]
+ A_k F(x_k) - {B}_k \\
\\
& \overset{(*)}{\geq} &
A_{k + 1} \bigl[ f(x_{k + 1}) 
+ \la \nabla f(x_{k + 1}), \frac{a_{k + 1} x + A_k x_k}{A_{k + 1}} - x_{k + 1} \ra \bigr]
+ a_{k + 1} \psi(x) + A_k \psi(x_k) \\
\\
& \quad & \qquad - \; {B}_k \\
\\
& \overset{(**)}{\geq} &
A_{k + 1} \bigl[ f(x_{k + 1}) 
+ \la \nabla f(x_{k + 1}), \frac{a_{k + 1} x + A_k x_k}{A_{k + 1}} - x_{k + 1} \ra 
+ \psi\bigl( \frac{a_{k + 1} x + A_k x_k}{A_{k + 1}}  \bigr)
\bigr] -  {B}_k,
\ea
$$
where $(*)$ and $(**)$ stand for convexity of $f$ and $\psi$, correspondingly.
Using both stationary condition and smoothness, we obtain,
for all $x \in \dom \psi$
$$
\ba{cl}
& \la \nabla f(x_{k + 1}), \frac{a_{k + 1} x + A_k x_k}{A_{k + 1}} - x_{k + 1} \ra
+ \psi \bigl( \frac{a_{k + 1} x + A_k x_k}{A_{k + 1}}  \bigr) \\
\\
& \quad = \quad
\gamma_k \la \nabla f(x_{k + 1}), x - v_{k + 1} \ra
+ \psi\bigl( \gamma_k x + (1 - \gamma_k) x_k\bigr) \\
\\
& \quad = \quad
\gamma_k \la \nabla f(x_k) + \nabla^2 f(x_k)(x_{k + 1} - x_k), x - v_{k + 1} \ra
+ \psi\bigl(  \gamma_k x + (1 - \gamma_k) x_k \bigr) \\
\\
& \quad \qquad \qquad + \;
\gamma_k \la \nabla f(x_{k + 1}) - \nabla f(x_k) 
- \nabla^2 f(x_k)(x_{k + 1} - x_k), x - v_{k + 1} \ra \\
\\
& \;\overset{\eqref{Alg2StatCondition},\eqref{GradHolderBound}}{\geq} \;
\psi(x_{k + 1}) - \frac{\gamma_k H_{\nu} \|x_{k + 1} - x_k \|^{1 + \nu} \|x - v_{k + 1}\|}{1 + \nu}
\;\; = \;\;
\psi(x_{k + 1}) - \frac{\gamma_k^{2 + \nu} H_{\nu}\|v_{k + 1} - x_k \|^{1 + \nu}  
	\|x - v_{k + 1}\|}{1 + \nu} \\
\\
& \quad \geq \quad
\psi(x_{k + 1}) - \frac{\gamma_k^{2 + \nu} H_{\nu} D^{2 + \nu}}{1 + \nu}.
\ea
$$
Therefore, we have
$$
\ba{rcl}
\phi_{k + 1}(x) & \geq & 
A_{k + 1} \bigl[ f(x_{k + 1}) + \psi(x_{k + 1}) 
- \frac{\gamma_k^{2 + \nu} H_{\nu} D^{2 + \nu}}{1 + \nu} \bigr] - {B}_k \\
\\
& = & A_{k + 1} F(x_{k + 1}) - {B}_{k + 1},
\ea
$$
and~\eqref{Th3Induction} is justified for all $k \geq 0$. Finally, by convexity of $f$, we get
$$
\ba{rcl}
F(x_k) - F^{*} & \leq & \ell_k \;\; \Def \;\; F(x_k) - \frac{\phi_k^{*}}{A_k} \\
\\
& \refLE{Th3Induction} & 
\frac{{B}_k}{A_k} 
\;\; = \;\; \frac{H_{\nu} D^{2 + \nu}}{(1 + \nu)A_k} 
\sum \limits_{i = 1}^k \frac{a_i^{2 + \nu}}{A_i^{1 + \nu}} \\
\\
& = & \O\bigl( \frac{H_{\nu} D^{2 + \nu}}{k^{1 + \nu}} \bigr),
\ea
$$
where the last equation holds from the choice $A_k := k^3$ (see the end of the proof of Theorem~\ref{Th:ContrDomNewtonConvex}).
\qed

\section{Convergence of Aggregating Newton Method}

In this section, we establish the convergence result for Algorithm~\ref{AggrNewton}.

\subsection{Proof of Theorem~\ref{Th:AggrNewtonGuarantee}}
\BT
For the sequence $\{ x_k
\}_{k \geq 1}$ generated by Algorithm~\ref{AggrNewton},
relation \eqref{QEstLower} is satisfied. 
\ET
\proof
Let us establish the relation~\eqref{QEstLower} by induction. It obviously
holds for $k = 0$. Assume that it is proven for the current iterate $k \geq 0$, 
and consider the next step:
$$
\ba{cl}
& Q_{k + 1}(v_{k + 1}) \\
\\
& \quad \equiv \quad
a_{k + 1} \bigl[ 
f(x_k) + \la \nabla f(x_k), v_{k + 1} - x_k \ra 
+ \frac{\gamma_k}{2} \la \nabla^2 f(x_k)(v_{k + 1} - x_k), v_{k + 1} - x_k \ra \\
\\
& \qquad \qquad \qquad \; + \; \psi(v_{k + 1}) \bigr] \; + \; Q_k(v_{k + 1}) \\
\\
& \;\,\,\, \refGE{QEstLower} \;\;\,
a_{k + 1} \bigl[ 
f(x_k) + \la \nabla f(x_k), v_{k + 1} - x_k \ra 
+ \frac{\gamma_k}{2} \la \nabla^2 f(x_k)(v_{k + 1} - x_k), v_{k + 1} - x_k \ra \\
\\
& \qquad \qquad \qquad \; + \; \psi(v_{k + 1}) \bigr] \; + \; A_k F(x_k) - \frac{C_k}{2} \\
\\
& \quad = \quad
A_{k + 1}\bigl[ f(x_{k}) +  \gamma_k \la \nabla f(x_k), v_{k + 1} - x_k \ra
+ \frac{\gamma_k^2}{2} \la \nabla^2 f(x_k)(v_{k + 1} - x_k), v_{k + 1} - x_k \ra ] \\
\\
& \qquad \qquad \qquad \; + \; a_{k + 1} \psi(v_{k + 1}) + A_k \psi(x_k) - \frac{C_k}{2} \\
\\
& \quad = \quad 
A_{k + 1}\bigl[ f(x_{k}) +  \la \nabla f(x_k), x_{k + 1} - x_k \ra
+ \frac{1}{2} \la \nabla^2 f(x_k)(x_{k + 1} - x_k), x_{k + 1} - x_k \ra ] \\
\\
& \qquad \qquad \qquad \; + \; a_{k + 1} \psi(v_{k + 1}) + A_k \psi(x_k) - \frac{C_k}{2} \\
\\
& \;\,\,\, \refGE{FuncHolderBound} \;\;\, 
A_{k + 1} \bigl[ f(x_{k + 1}) 
- \frac{H_{\nu}\|x_{k + 1} - x_k \|^{2 + \nu}}{(1 + \nu)(2 + \nu)}  \bigr]
+ a_{k + 1} \psi(v_{k + 1}) + A_k \psi(x_k) - \frac{C_k}{2} \\
\\
& \quad = \quad 
A_{k + 1} f(x_{k + 1}) 
- \frac{ A_{k + 1} \gamma_k^{2 + \nu} H_{\nu} \|v_{k + 1} - x_k \|^{2 + \nu}}{(1 + \nu)(2 + \nu)} 
+ a_{k + 1} \psi(v_{k + 1}) + A_k \psi(x_k) - \frac{C_k}{2} \\
\\
& \quad \geq \quad
A_{k + 1} f(x_{k + 1})
- \frac{a_{k + 1} \gamma_k^{1 + \nu} \H_{\nu} D^{2 + \nu}}{(1 + \nu)(2 + \nu)}
+ A_{k + 1} \psi(x_{k + 1}) - \frac{C_k}{2} \\
\\
& \quad = \quad A_{k + 1} F(x_{k + 1}) - \frac{C_{k + 1}}{2}.
\ea 
$$
Thus, we have~\eqref{QEstLower} justified for all $k \geq 0$.
\qed

\section{Convergence of stochastic methods}

Let us consider the following general iterations, 
for solving optimization problem~\eqref{MainProblem}:
\beq \label{GeneralContrDomNewton}
\ba{rcl}
x_{k + 1} & \in & \Argmin\limits_{y}
\Bigl\{  \la g_k, y - x_k \ra + \frac{1}{2}\la H_k(y - x_k), y - x_k \ra
+ S_k(y) \Bigr\}, \quad k \geq 0
\ea
\eeq
with $S_k(y) := \gamma_k \psi(x_k + \frac{1}{\gamma_k}(y - x_k))$.
This is Algorithm~\ref{ContrDomNewton} with substituted vector $g_k$
and matrix $H_k$ instead of the true gradient and the Hessian.
First, we need to study the convergence of this process. For simplicity,
let us study the case $\nu = 1$ only
(convex functions with Lipschitz continuous Hessian, 
we denote the corresponding Lipschitz constant by $L_2$).
Recall, that in this section we use the standard Euclidean norm 
for vectors and induced spectral norm for matrices.

As before, we use the sequence of positive numbers
 $\{ a_k \}_{k \geq 1}$, and set
 $$
\ba{rcl}
\gamma_k & := & \frac{a_{k + 1}}{A_{k + 1}}, 
\qquad A_k \;\; \Def \;\; \sum\limits_{i = 1}^k a_i.
\ea
$$

\BL \label{Lemma:GeneralContrDomNewton}
For iterations~\eqref{GeneralContrDomNewton}, we have
for all $k \geq 1$
\beq \label{GeneralItersConv}
\ba{rcl}
F(x_k) - F^{*} & \leq & \frac{B_k}{A_k},
\ea
\eeq
with
$$
\ba{rcl}
B_k & := &
\frac{L_2 D^3}{2} \sum\limits_{i = 0}^{k - 1} \frac{a_{i + 1}^3}{A_{i + 1}^2}
+ D \sum\limits_{i = 0}^{k - 1} a_{i + 1} \| \nabla f(x_i) - g_i \|
+ D^2 \sum\limits_{i = 0}^{k - 1} \frac{a_{i + 1}^2}{A_{i + 1}} \| \nabla^2 f(x_i) - H_i \|.
 \ea
$$
\EL
\proof
Let us prove by induction the following inequality
\beq \label{GeneralInduction}
\ba{rcl}
A_k F(x) & \geq & A_k F(x_k) - B_k, \qquad x \in \dom \psi.
\ea
\eeq
It obviously holds for $k = 0$, and for $k \geq 1$ it is equivalent to~\eqref{GeneralItersConv}.

Assume that~\eqref{GeneralInduction} is satisfied for some $k \geq 0$, and consider the next step:
\beq \label{GenProofFirst}
\ba{rcl}
A_{k + 1} F(x) & = & a_{k + 1} F(x) + A_k F(x) \\
\\
& \refGE{GeneralInduction} &
a_{k + 1} F(x)+ A_k F(x_k) - B_k \\
\\
& \overset{(*)}{\geq} & A_{k + 1} f\bigl( \frac{a_{k + 1} x + A_k x_k}{A_{k + 1}} \bigr)
+ a_{k + 1} \psi(x) + A_k \psi(x_k) - B_k \\
\\
& \overset{(*)}{\geq} &
A_{k + 1} \bigl[ f(x_{k + 1}) 
+ \la \nabla f(x_{k + 1}), \frac{a_{k + 1} x + A_k x_k}{A_{k + 1}} - x_{k + 1} \ra \bigr]
+ a_{k + 1} \psi(x)  \\
\\
& \quad & \qquad + \; A_k \psi(x_k) - B_k,
\ea
\eeq
where $(*)$ stands for convexity of $f$. Now, let us denote
the point
$$
\ba{rcl}
v_{k + 1} & := & x_k + \frac{1}{\gamma_k} (x_{k + 1} - x_k) \;\; \in \;\; \dom \psi.
\ea
$$
Then, stationary condition for the method step~\eqref{GeneralContrDomNewton} can
be written as
\beq \label{GenStatCond}
\ba{rcl}
\la g_k + H_k(x_{k + 1} - x_k), x - v_{k + 1} \ra + \psi(x) & \geq & \psi(v_{k + 1}),
\ea
\eeq
for all $x \in \dom \psi$. Therefore,
\beq \label{GenProofSecond}
\ba{cl}
& A_{k + 1} \la \nabla f(x_{k + 1}), \frac{a_{k + 1} x + A_k x_k}{A_{k + 1}} - x_{k + 1} \ra
+ a_{k + 1} \psi(x) \\
\\
& \quad = \quad 
a_{k + 1} \bigl[ \la \nabla f(x_{k + 1}), x - v_{k + 1} \ra + \psi(x) \bigr] \\
\\
& \quad = \quad 
a_{k + 1} \bigl[   
\la g_k + H_k (x_{k + 1} - x_k), x - v_{k + 1} \ra + \psi(x) \\
\\
& \qquad \qquad \qquad 
+ \; \la \nabla f(x_k) - g_k, x - v_{k + 1} \ra \\
\\
& \qquad \qquad \qquad
+ \; \la (\nabla^2 f(x_k) - H_k)(x_{k + 1} - x_k), x - v_{k + 1} \ra \\
\\
& \qquad \qquad \qquad 
+ \; \la \nabla f(x_{k + 1}) - \nabla f(x_k) - \nabla^2 f(x_k)(x_{k + 1} - x_k), 
x - v_{k + 1} \ra
\bigr] \\
\\
& \; \overset{\eqref{GenStatCond},\eqref{GradHolderBound}}{\geq} \;
a_{k + 1}\bigl[ \psi(v_{k + 1}) - \| \nabla f(x_k) - g_k \| \cdot \|x - v_{k + 1}\| \\
\\
& \qquad \qquad \qquad 
- \; \gamma_k \| \nabla^2 f(x_k) - H_k\| \cdot \|v_{k + 1} - x_k \| \cdot \| x - v_{k + 1}\| \\
\\
& \qquad \qquad \qquad 
- \; \frac{L_2 \gamma_k^2 \|v_{k + 1} - x_k \|^2 \cdot \|x - v_{k + 1} \| }{2}
\bigr] \\
\\
& \quad \geq \quad 
a_{k + 1} \psi(v_{k + 1}) - a_{k + 1} D \| \nabla f(x_k) - g_k \|_{*}
- \frac{ a_{k + 1}^2 D^2 \| \nabla^2 f(x_k) - H_k \|}{A_{k + 1}}
- \frac{a_{k + 1}^3 L_2 D^3}{A_{k + 1}^2}.
\ea
\eeq
Thus, combining all together,
and using convexity of $\psi$, we obtain
$$
\ba{rcl}
A_{k + 1} F(x) & \overset{\eqref{GenProofFirst},\eqref{GenProofSecond}}{\geq} &
A_{k + 1} f(x_{k + 1}) + a_{k + 1} \psi(v_{k + 1}) + A_k \psi(x_k) - B_k \\
\\
& \quad & \qquad 
- \; a_{k + 1} D \| \nabla f(x_k) - g_k \|
- \frac{a_{k + 1}^2 D^2 \| \nabla^2 f(x_k) - H_k \|}{A_{k + 1}}
- \frac{a_{k + 1}^3 L_2 D^3}{A_{k + 1}^2} \\
\\
& \geq & A_{k + 1} F(x_{k + 1}) - B_{k + 1}.
\ea
$$
So, we have \eqref{GeneralInduction} justified for all $k \geq 0$.
\qed

Now, we are ready to prove convergence results
for the process~\eqref{GeneralContrDomNewton}
with the basic variant of stochastic estimators~\eqref{RandomEstimators},
and with the variance reduction strategy for the gradients,
incorporated into Algorithm~\ref{StochContrDomNewton}.

\subsection{Proof of Theorem~\ref{Th:StochNewton}}
\BT
Let each component $f_i(\cdot)$ be Lipschitz continuous on $\dom \psi$ with constant $L_0$,
and have Lipschitz continuous gradients and Hessians on $\dom \psi$
with constants $L_1$ and $L_2$, respectively.
Let $\gamma_k := 1 - \bigl( \frac{k}{k + 1} \bigr)^3 = \O\bigl(\frac{1}{k} \bigr)$. Set
$$
\ba{rcl}
m_k^g & := & 1 / \gamma_k^4,
\qquad
m_k^H \;\; := \;\; 1 / \gamma_k^2.
\ea
$$
Then, for the iterations $\{ x_k \}_{k \geq 1}$ of
Algorithm~\eqref{ContrDomNewton}, based on
estimators~\eqref{RandomEstimators}, it holds
$$
\ba{rcl}
\E[F(x_k) - F^{*}] & \leq & 
\O \Bigl( \frac{ L_2 D^3 \, + \, L_1 D^2 (1 + \log(n)) \, + \,  L_0 D  }{k^2} \Bigr).
\ea
$$
\ET
\proof
Let us fix iteration $k \geq 0$. For one uniform random sample $i \in \{1, \dots, M\}$, we have
\beq \label{GradOneSample}
\ba{rcl}
\E \| \nabla f(x_k) - \nabla f_i(x_k) \|^2
& = & \E \| \nabla f_i(x_k) \|^2  - \| \nabla f(x_k) \|^2
\;\; \leq \;\; L_0^2.
\ea
\eeq
Therefore, for the random batch of size $m_k^g$, we obtain
\beq \label{GradBatchSample}
\ba{rcl}
\E \| \nabla f(x_k) - g_k \| & \leq & \sqrt{ \E \| \nabla f(x_k) - g_k \|^2 } \\
\\
& = & 
\sqrt{  \frac{1}{(m_k^g)^2} 
\E \|  \sum_{i \in S_k^g} ( \nabla f(x_k) - \nabla f_i(x_k))  \|^2 } \\
\\
& = &
\sqrt{ 
\frac{1}{(m_k^g)^2} \sum_{i \in S_k^g} \E \| \nabla f(x_k) - \nabla f_i(x_k) \|^2} \\
\\
& \refLE{GradOneSample} &
\frac{L_0}{\sqrt{m_k^g}}.
\ea
\eeq
More advanced reasoning for matrices (Matrix Bernstein Inequality; see Chapter 6 in~\cite{tropp2015introduction}) gives
\beq \label{HessBatchSample}
\ba{rcl}
\E \| \nabla^2 f(x_k) - H_k \| & \leq & 
L_1 \Bigl(  \sqrt{ \frac{2 \log(2n)}{m_k^H} } + \frac{2 \log(2n)}{3m_k^H}  \Bigr) \\
\\
& \leq & \frac{L_1( 3 \sqrt{2 \log(2n)} + 2 \log(2n)  )}{3 \sqrt{m_k^H}}
\; \leq \;
\frac{L_1(  6 + 7 \log(2n)   )}{6\sqrt{m_k^H}}.
\ea
\eeq
So, using these estimates together, we have, for every $k \geq 1$
$$
\ba{rcl}
\E[ F(x_k) - F^{*} ]
& \refLE{GeneralItersConv} & 
\frac{1}{A_k} \Bigl( 
\frac{L_2 D^3}{2} \sum\limits_{i = 0}^{k - 1} \frac{a_{i + 1}^3}{A_{i + 1}^2}
+ D \sum\limits_{i = 0}^{k - 1} a_{i + 1} \E \| \nabla f(x_i) - g_i \| \\
\\
& \quad & \qquad \qquad 
+ \; D^2 \sum\limits_{i = 0}^{k - 1} \frac{a_{i + 1}^2}{A_{i + 1}} \E\| \nabla^2 f(x_i) - H_i \|
\Bigr) \\
\\
& \overset{\eqref{GradBatchSample},\eqref{HessBatchSample}}{\leq} &
\frac{1}{A_k} \Bigl(
\frac{L_2 D^3}{2} \sum\limits_{i = 0}^{k - 1} \frac{a_{i + 1}^3}{A_{i + 1}^2}
+ L_0 D \sum\limits_{i = 0}^{k - 1} \frac{a_{i + 1}}{\sqrt{m_i^g}} \\
\\
& \quad & \qquad \qquad
+ \; \frac{L_ 1 D^2 (6 + 7 \log(2n))}{6} 
\sum\limits_{i = 0}^{k - 1} \frac{a_{i + 1}^2}{A_{i + 1} \sqrt{m_i^H}}
\Bigr) \\
\\
& \refEQ{BasicBatchSize} &
\frac{1}{A_k} \Bigl(
\frac{L_2 D^3}{2} + L_0 D + \frac{L_1 D^2 (6 + 7 \log(2n))}{6}
\Bigr)
\sum\limits_{i = 0}^{k - 1} \frac{a_{i + 1}^3}{A_{i + 1}^2}. \\
\ea
$$
Thus, for the choice $A_k := k^3$, we get
$$
\ba{rcl}
\E[ F(x_k) - F^{*} ] & \leq & 
\O\Bigl( \frac{ L_2 D^3 + L_1 D^2 (1 + \log(n)) +  L_0 D }{k^2} \Bigr).
\ea
$$
\qed

\subsection{Proof of Theorem~\ref{Th:StochVRNewton}}
\BT 
Let each component $f_i(\cdot)$ have Lipschitz continuous gradients and Hessians
on $\dom \psi$ with constants $L_1$ and $L_2$, respectively.
Let $\gamma_k := 1 - \bigl( \frac{k}{k + 1} \bigr)^3 = \O(\frac{1}{k})$. Set batch size
$$
\ba{rcl}
m_k & := & 1 / \gamma_k^2.
\ea
$$
Then, for all iterations $\{ x_k \}_{k \geq 1}$ of
Algorithm~\ref{StochContrDomNewton}, we have
$$
\ba{rcl}
\E[ F(x_k) - F^{*} ] & \leq & 
\O \Bigl( \frac{ L_2 D^3 \, + \, L_1 D^2 (1 + \log(n)) \, + \, L_1^{1/2} D (F(x_0) - F^{*}) }{k^2} \Bigr).
\ea
$$
\ET
\proof

Let us consider the following stochastic estimate
$$
\ba{rcl}
g_k^i & := & \nabla f_i(x_k) - \nabla f_i(z_k) + \nabla f(z_k),
\ea
$$
for a uniform random sample $i \in \{1, \dots, M\}$,
and a current iterate $k \geq 0$.
We denote by $x^{*}$ the solution of our problem: $F^{*} = F(x^{*})$,
stationary condition for which is
\beq \label{MainProblemStatCond}
\ba{rcl}
\la \nabla f(x^{*}), x - x^{*} \ra + \psi(x) & \geq & \psi(x^{*}),
\qquad x \in \dom \psi.
\ea
\eeq
Then, it holds
$$
\ba{rcl}
\E \| \nabla f(x_k) - g_k^{i} \|^2
& = &
\E \| (\nabla f(x_k) - \nabla f(x^{*}))  \\
\\
& \quad & \quad 
+ \quad ( \nabla f_i(z_k) - \nabla f_i (x^{*}) - \nabla f(z_k) + \nabla f(x^{*}) ) \\
\\
& \quad & \quad 
+ \quad (\nabla f_i(x^{*}) - \nabla f_i(x_k)) \|^2 \\
\\
& \leq &
3  \E \| \nabla f(x_k) - \nabla f(x^{*}) \|^2 \\
\\
& \quad & \quad 
+ \quad 3\E \| (\nabla f_i(z_k) - \nabla f_i(x^{*})) - (\nabla f(z_k) - \nabla f(x^{*})) \|^2 \\
\\
& \quad & \quad
+ \quad 3\E \| \nabla f_i(x_k) - \nabla f_i(x^{*}) \|^2 \\
\\
& \leq &
3\Bigl( \E\| \nabla f(x_k) - \nabla f(x^{*}) \|^2
+ \E \| \nabla f_i(z_k) - \nabla f_i(x^{*}) \|^2 \\
\\
& \quad & \quad
+ \quad \E \| \nabla f_i(x_k) - \nabla f_i(x^{*}) \|^2 \Bigr),
\ea
$$
where we used the following simple bounds:
$$
\ba{rcl}
\| a + b + c\|^2 & \leq & 3\|a\|^2 + 3\|b\|^2 + 3 \|c\|^2, \\
\\
\E \| \xi - \E \xi \|^2 & \leq & \E \| \xi \|^2,
\ea
$$
which are valid for any $a, b, c \in \R^n$ and arbitrary random vector $\xi \in \R^n$.

Now, by Lipschitz continuity of the gradients, we have 
(see Theorem 2.1.5 in~\cite{nesterov2018lectures})
$$
\ba{rcl}
\| \nabla f(x_k) - \nabla f(x^{*}) \|^2 
& \leq &
2L_1 \bigl( f(x_k) - f(x^{*}) - \la \nabla f(x^{*}), x_k - x^{*} \ra \bigr) \\
\\
& \refLE{MainProblemStatCond} & 
2L_1 \bigl( F(x_k) - F^{*} \bigr).
\ea
$$
The same holds for the random sample $i$, for arbitrary fixed $x \in \dom \psi$
$$
\ba{rcl}
\E_i \| \nabla f_i(x) - \nabla f_i(x^{*}) \|^2
& \leq & 2L_1 \E_i 
\bigl[ 
f_i(x) - f_i(x^{*}) - \la \nabla f_i(x^{*}), x - x^{*} \ra
\bigr] \\
\\
& = & 
2L_1 \bigl( 
f(x) - f(x^{*}) - \la \nabla f(x^{*}), x - x^{*} \ra
\bigr) \\
\\
& \refLE{MainProblemStatCond} & 
2L_1 \bigl( F(x) - F^{*} ).
\ea
$$
Thus, we obtain
\beq \label{GradOneVREstimate}
\ba{rcl}
\E \| \nabla f(x_k) - g_k^{i} \|^2 
& \leq & 
12 L_1 \E[ F(x_k) - F^{*} ] + 6 L_1 \E[ F(z_k) - F^{*} ].
\ea
\eeq
Consequently, for the random batch
$$
\ba{rcl}
g_k & := & \frac{1}{m_k} \sum_{i \in S_k} g_k^i,
\ea
$$
we have (compare with~\eqref{GradBatchSample})
\beq \label{GradBatchVREstimate}
\ba{rcl}
\E \| \nabla f(x_k) - g_k \| & \leq & 
\sqrt{ \frac{1}{(m_k)^2}  
\sum_{i \in S_k} \E \| \nabla f(x_k) - g_k^{i} \|^2} \\
\\
& \refLE{GradOneVREstimate} &
\sqrt{\frac{6L_1}{m_k}\bigl(  2\E[ F(x_k) - F^{*} ] + \E [F(z_k) - F^{*}]  \bigr)  } \\
\\
& \leq &
\sqrt{\frac{12 L_1}{m_k} \E[ F(x_k) - F^{*}]}
+ \sqrt{\frac{6 L_1}{m_k} \E[ F(z_k) - F^{*} ] }.
\ea
\eeq
So, using the variance reduction for the gradients, 
and the basic estimate for the Hessians, we have, for every $k \geq 1$
$$
\ba{rcl}
\E [ F(x_k) - F^{*} ]
& \overset{\eqref{GeneralItersConv},\eqref{GradBatchVREstimate},\eqref{HessBatchSample}}{\leq} &
\frac{1}{A_k} \Bigl(
\frac{L_2 D^3}{2} \sum\limits_{i = 0}^{k - 1} \frac{a_{i + 1}^3}{A_{i + 1}^2} \\
\\
& \quad & \quad
+ \; D \sqrt{6L_1} \sum\limits_{i = 0}^{k - 1} \frac{a_{i + 1}}{\sqrt{m_i}}
\bigl(  \sqrt{2\E[F(x_i) - F^{*}]} + \sqrt{\E[F(z_i) - F^{*}]}   \bigr) \\
\\
& \quad & \quad
+ \; \frac{L_1 D^2(6 + 7 \log(2n))}{6} 
\sum\limits_{i = 0}^{k - 1} \frac{a_{i + 1}^2}{A_{i + 1}\sqrt{m_i}}
\Bigr) \\
\\
& \overset{\eqref{VRBatchSize}}{=} &
\frac{1}{A_k} \Bigl(
\Bigl[
\frac{3L_2 D^3 + L_1 D^2 (6 + 7 \log(2n))}{6}
\Bigr] \sum\limits_{i = 0}^{k - 1} \frac{a_{i + 1}^3}{A_{i + 1}^2} \\
\\
& \quad & \quad
+ \; D \sqrt{6L_1} \sum\limits_{i = 0}^{k - 1}
\frac{a_{i + 1}^2}{A_{i + 1}}
\bigl(  \sqrt{2\E[F(x_i) - F^{*}]} + \sqrt{\E[F(z_i) - F^{*}]}   \bigr)
\Bigr).
\ea
$$
Now, let us set $A_{i + 1} := (i + 1)^3$, and thus $a_{i + 1} := (i + 1)^3 - i^3 \leq 3(i + 1)^2$, so
we have
\beq \label{StochVRResBound}
\ba{rcl}
\E [ F(x_k) - F^{*} ]
& \leq & 
\frac{\alpha \, + \, \beta (\sqrt{2} + 1) (F(x_0) - F^{*})}{k^2}  \\
\\
& \quad & \quad 
+ \; \frac{\beta}{k^3} \sum\limits_{i = 1}^{k - 1}
\Bigl( (i + 1) \bigl( 
\sqrt{2\E[F(x_i) - F^{*}]} + \sqrt{\E[F(z_i) - F^{*}]}
  \bigr)
  \Bigr),
\ea
\eeq
where
$$
\ba{rcl}
\alpha & := & 27 \cdot \Bigl[
\frac{3L_2 D^3 + L_1 D^2 (6 + 7 \log(2n))}{6}
\Bigr], \qquad 
\beta \;\; := \;\; 9 \cdot D \sqrt{6 L_1}.
\ea
$$
We are going to prove by induction, for every $k \geq 1$
\beq \label{StochVRInduction}
\ba{rcl}
\E[F(x_k) - F^{*}] & \leq & \frac{c}{k^2},
\ea
\eeq
with
\beq \label{cDef}
\ba{rcl}
c & := & \bigl(  4 \beta + \sqrt{\alpha + 3 \beta (F(x_0) - F^{*}) + 16 \beta^2}  \bigr)^2 
\;\; \leq \;\;
74 \beta^2 + 2 \alpha + 6 \beta(F(x_0) - F^{*}) \\
\\
& = &
\O\bigl(  L_2 D^3  + L_1 D^2 (1 + \log(n)) + L_1^{1/2} D(F(x_0) - F^{*}) \bigr).
\ea
\eeq
Hence, if \eqref{StochVRInduction} is true, then we essentially obtain the claim of the theorem.
For $k = 1$, \eqref{StochVRInduction} follows directly from~\eqref{StochVRResBound}.
Assume that~\eqref{StochVRInduction} holds for all $1 \leq i \leq k$, and consider
iteration $k + 1$:
$$
\ba{rcl}
\E[F(x_{k + 1}) - F^{*}] 
& \overset{\eqref{StochVRResBound},\eqref{StochVRInduction}}{\leq} &
\frac{\alpha \, + \, \beta (\sqrt{2} + 1) (F(x_0) - F^{*})}{k^2} 
+ 
\frac{\beta}{k^3} \sum\limits_{i = 1}^{k}
\Bigl( (i + 1) \Bigl( 
\frac{\sqrt{2c}}{i} + \frac{\sqrt{c}}{\pi(i)}
\Bigr)
\Bigr) \\
\\
& \overset{(*)}{\leq} &
\frac{\alpha \, + \, \beta (\sqrt{2} + 1) (F(x_0) - F^{*})}{k^2} 
+ 
\frac{\beta \sqrt{c}}{k^3} \sum\limits_{i = 1}^{k}
\Bigl( (i + 1) \Bigl( 
\frac{2\sqrt{2} \, + \, 4}{i + 1}
\Bigr)
\Bigr) \\
\\
& = &
\frac{\alpha \, + \, (\sqrt{2} + 1)\beta(F(x_0) - F^{*}) 
\, + \, (2\sqrt{2} + 4) \beta \sqrt{c}   }{k^2} \\
\\
& \leq & 
\frac{ \alpha \, + \, 3 \beta (F(x_0) - F^{*}) \, + \, 8 \beta \sqrt{c} }{k^2}
\;\; \refEQ{cDef} \;\; \frac{c}{k^2},
\ea
$$
where in $(*)$ we have used two simple bounds: $i \leq 2\pi(i)$,
and $i + 1 \leq 2i$, valid for all $i \geq 1$.
\qed

\newpage

\section{Extra experiments}
\label{SectionEE}

\begin{figure}[h!]
	\centering
	\includegraphics[width=0.32\textwidth ]{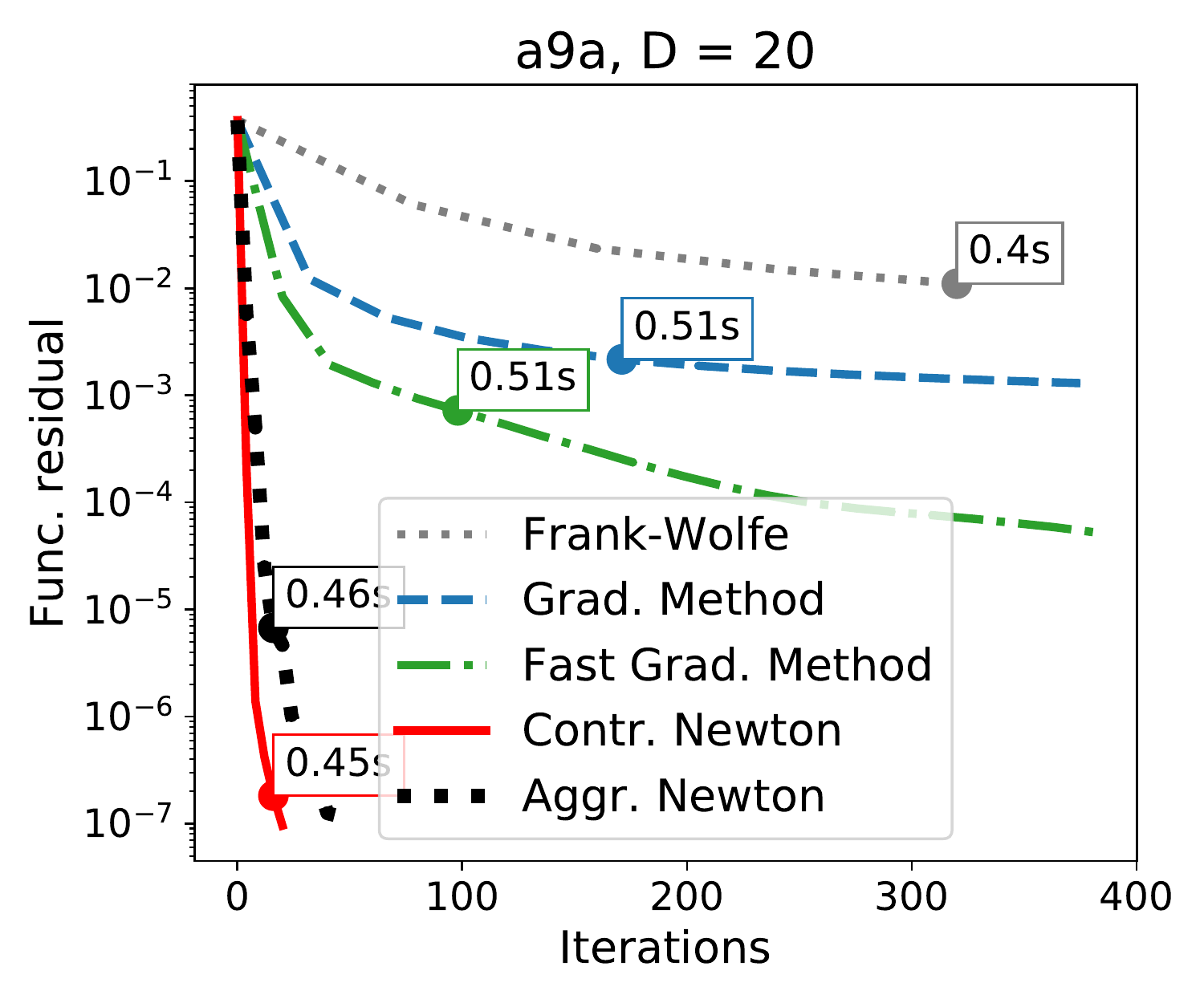}
	\includegraphics[width=0.32\textwidth ]{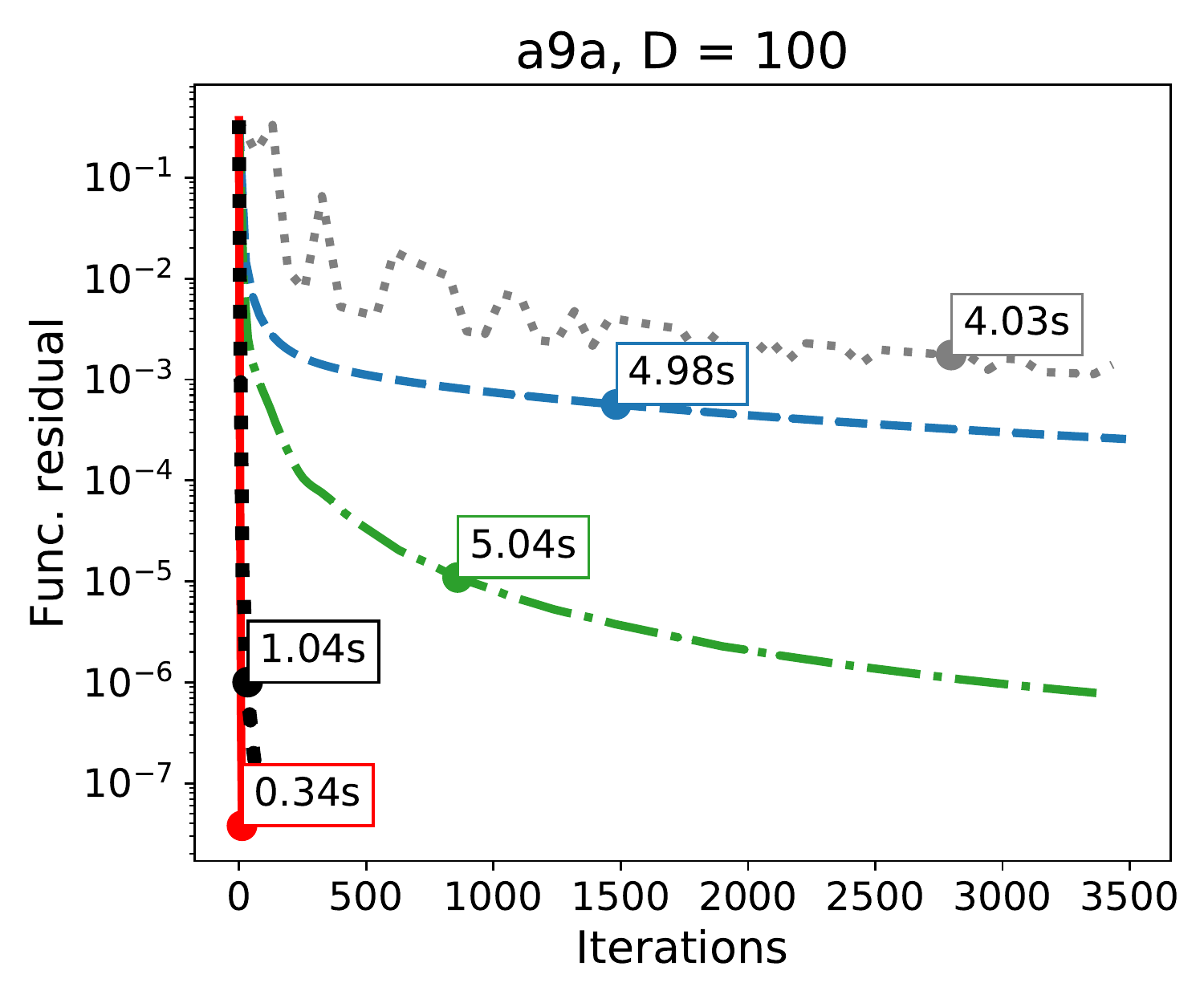}
	\includegraphics[width=0.32\textwidth ]{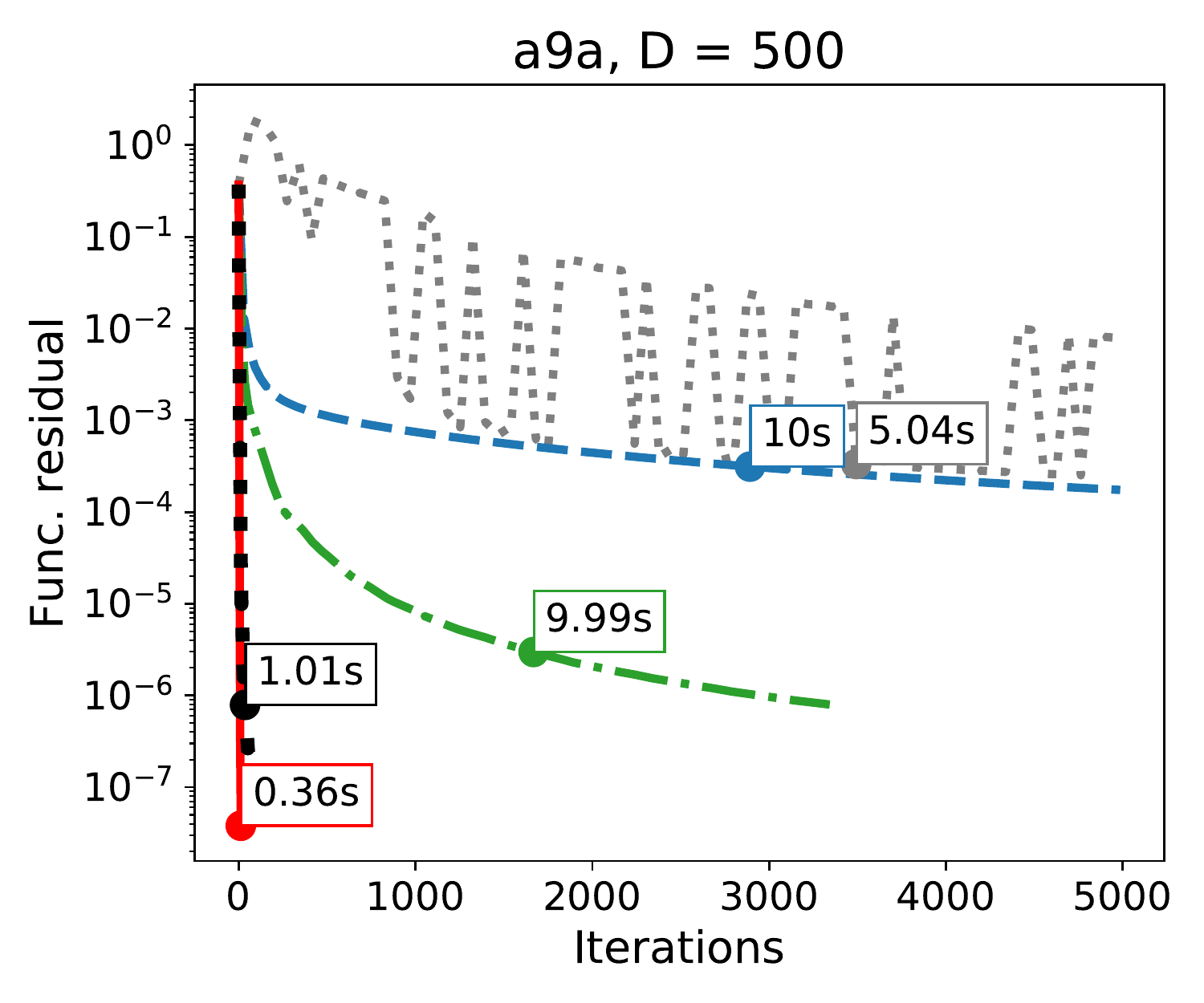}
	
	\includegraphics[width=0.32\textwidth ]{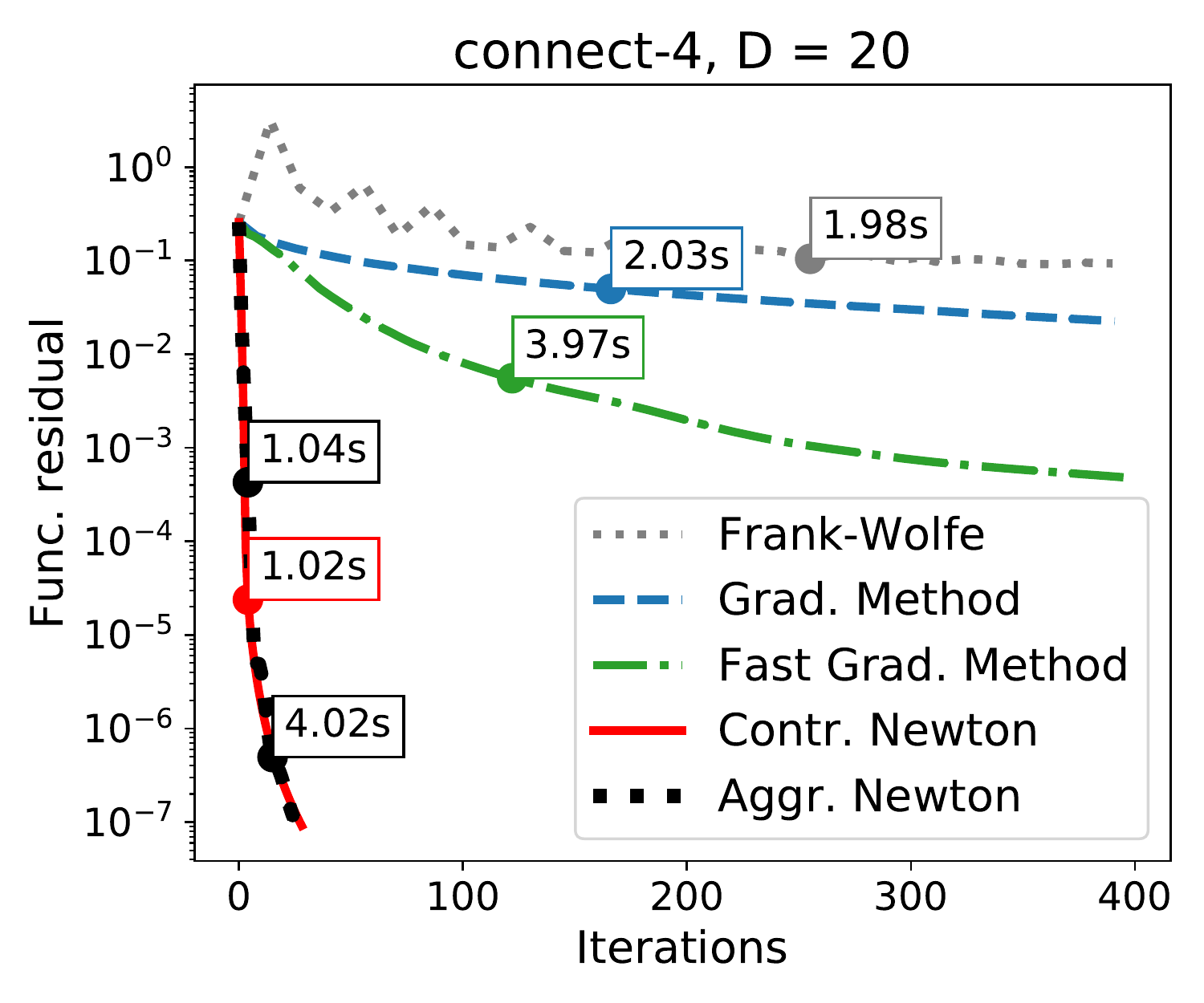}
	\includegraphics[width=0.32\textwidth ]{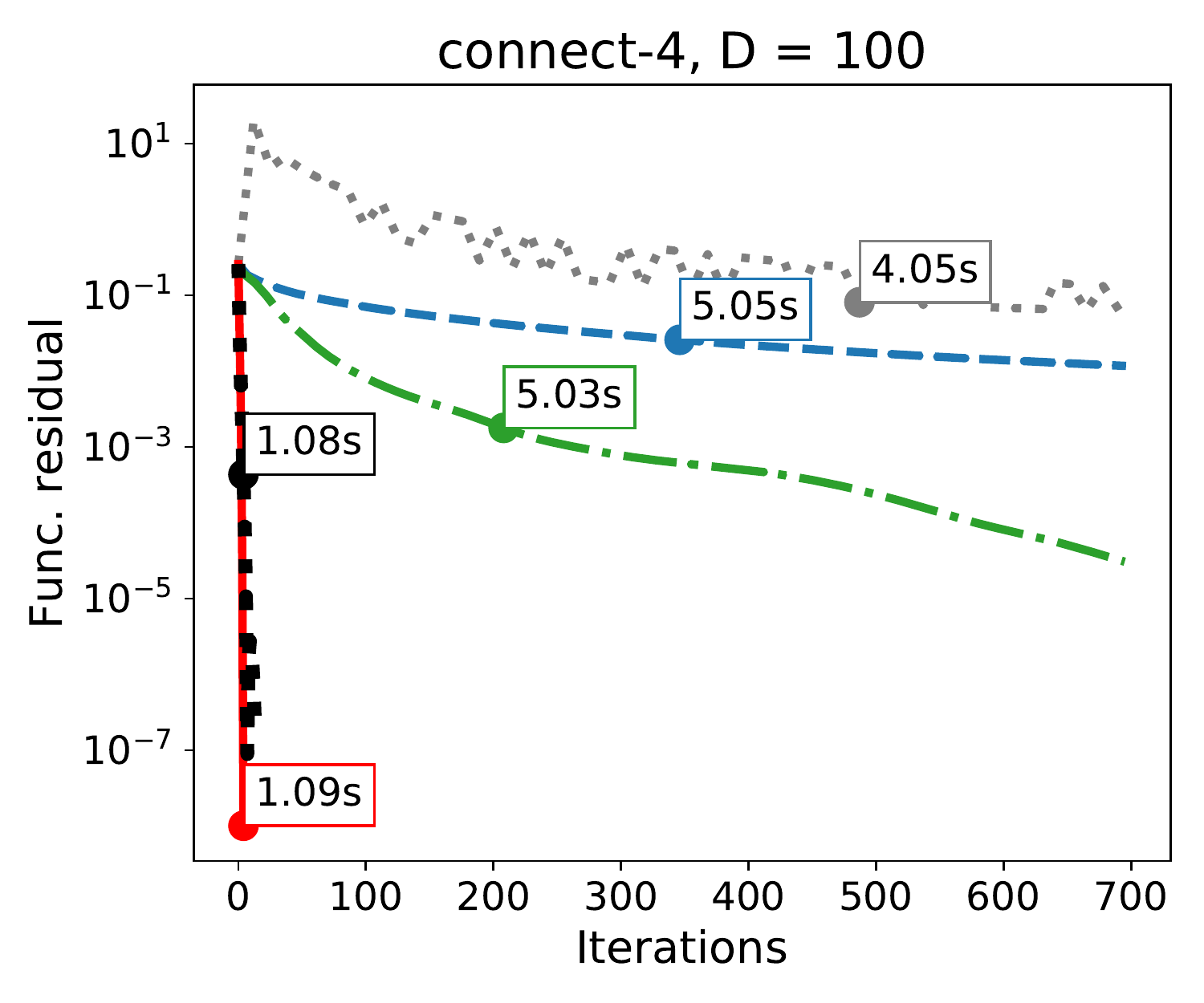}
	\includegraphics[width=0.32\textwidth ]{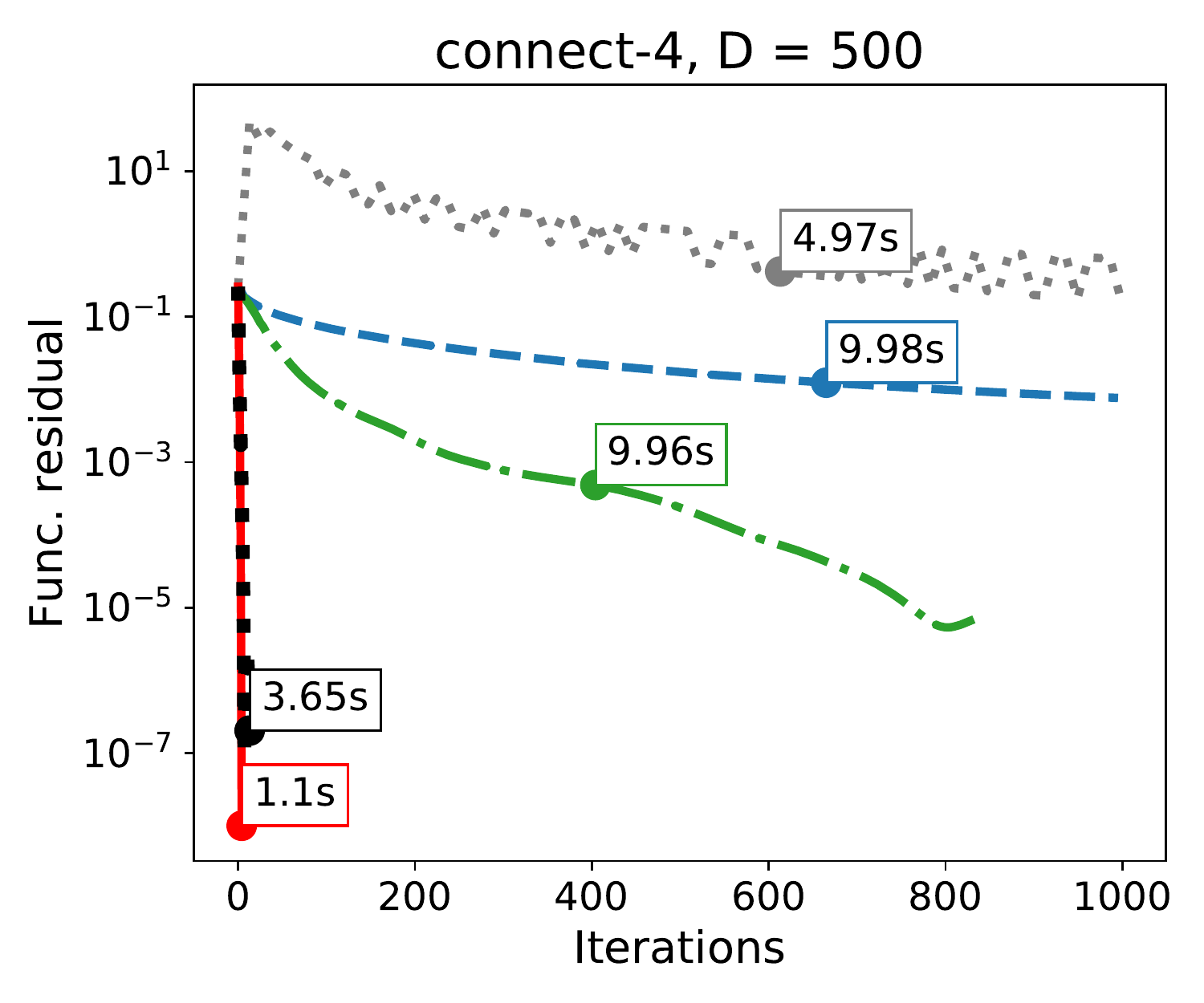}
	
	\includegraphics[width=0.32\textwidth ]{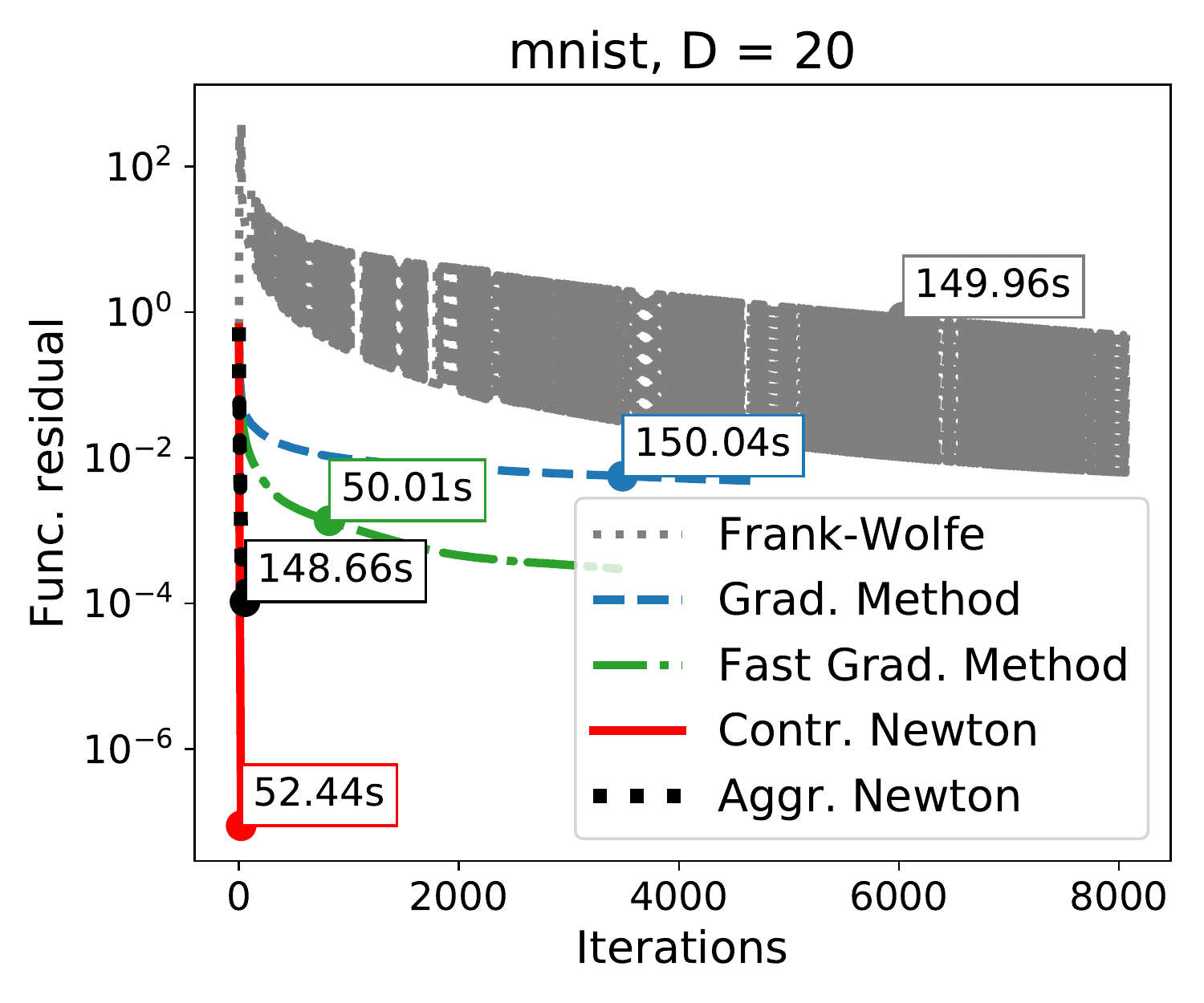}
	\includegraphics[width=0.32\textwidth ]{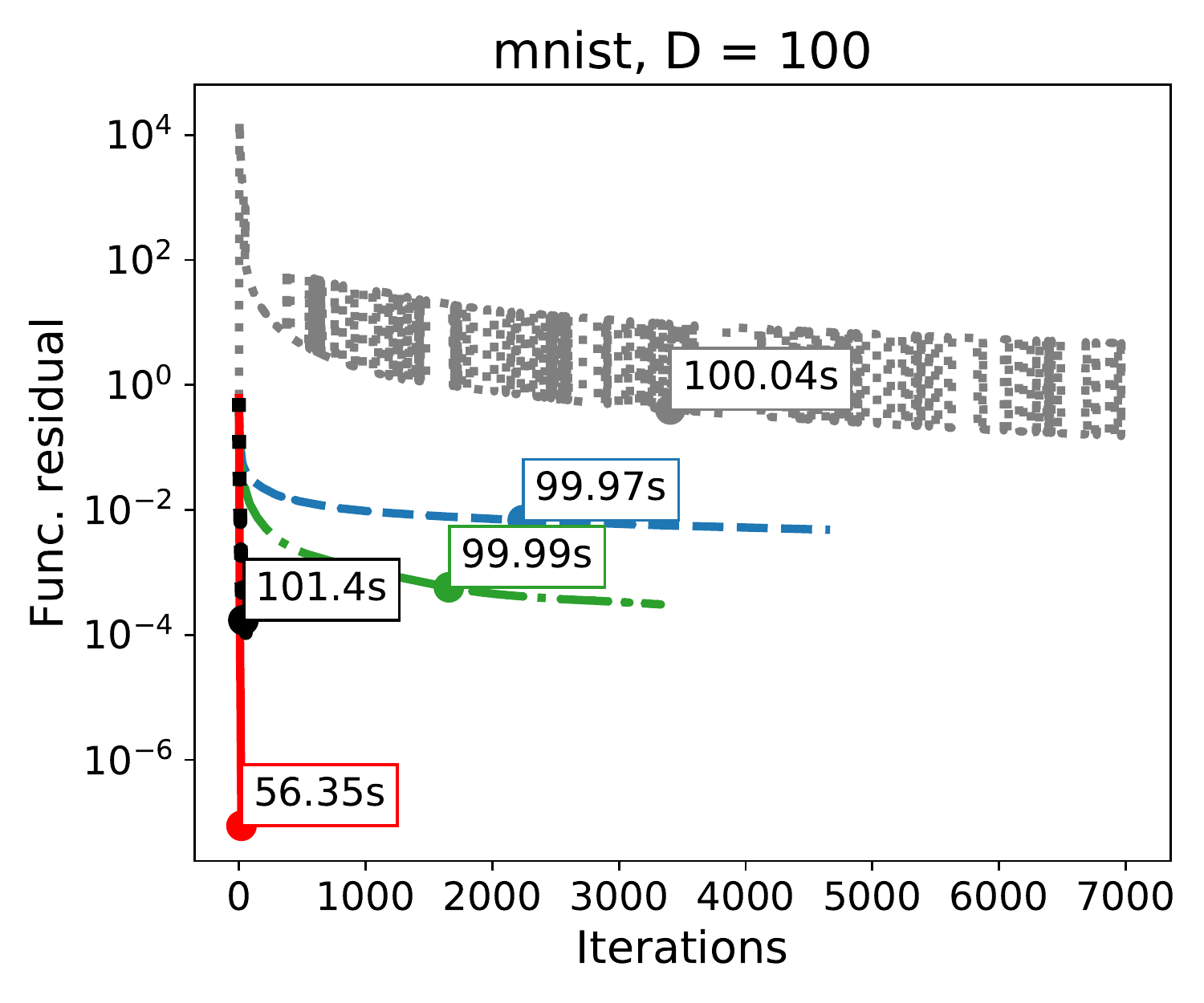}
	\includegraphics[width=0.32\textwidth ]{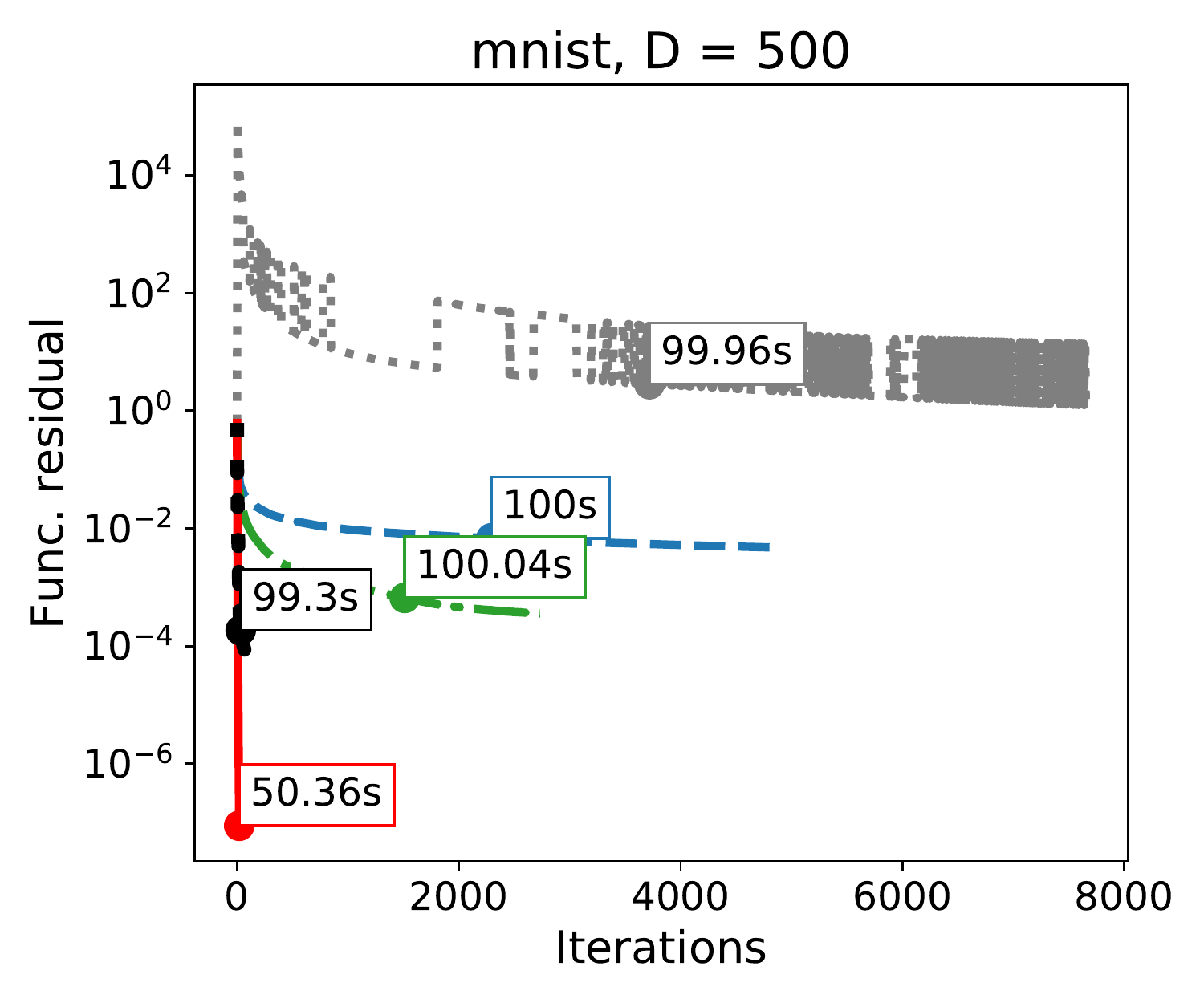}
	\caption{Training logistic regression, datasets: \textit{a9a} $(M = 32561, n = 123)$, 
	 \textit{connect-4} $(M = 67557, n = 126)$, \textit{mnist} $(M = 60000, n = 780)$.  }
	\label{Fig:logreg_extra_basic}
\end{figure}

\begin{figure}[h!]
	\centering
	
	\includegraphics[width=0.32\textwidth ]{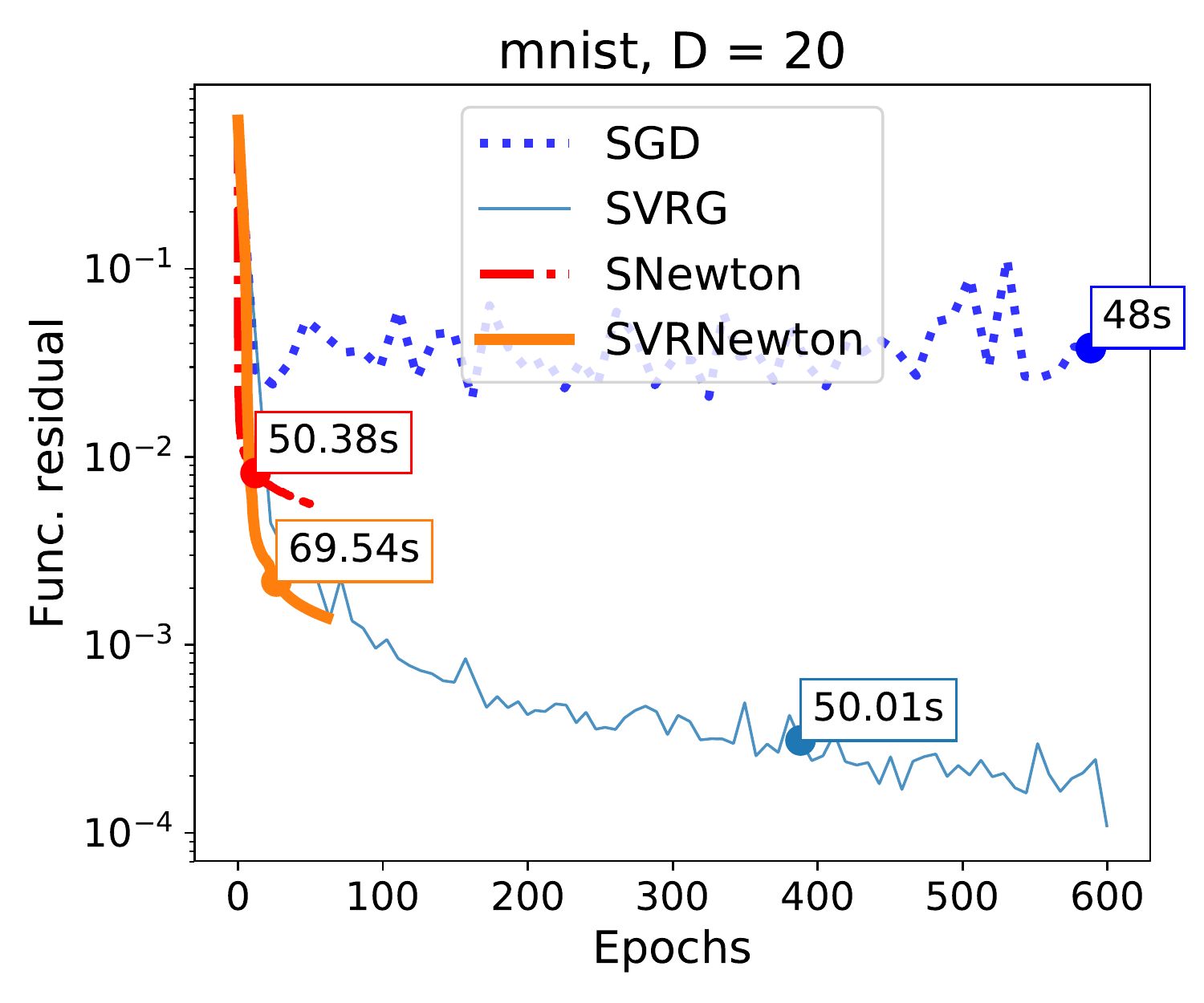}
	\includegraphics[width=0.32\textwidth ]{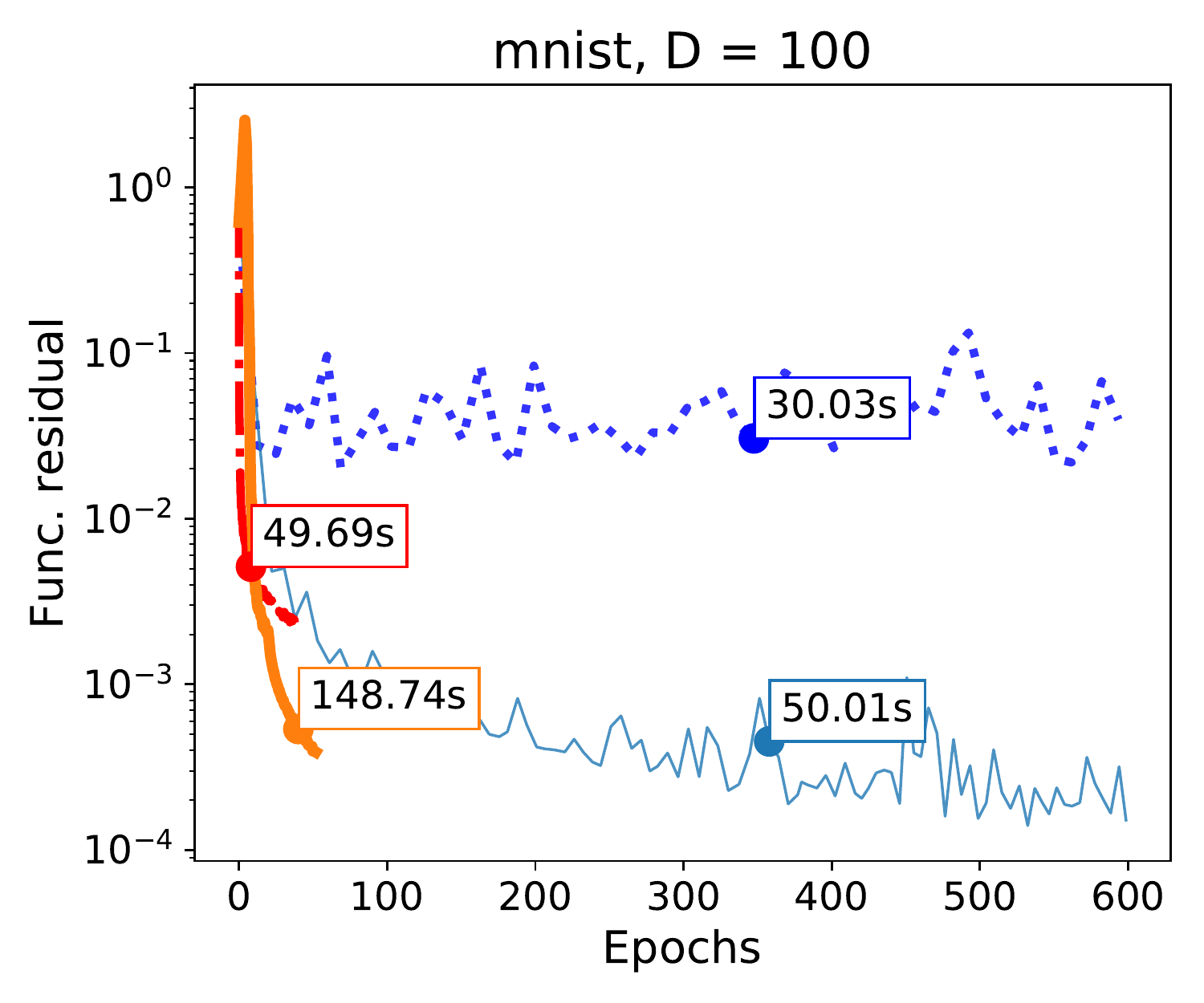}
	\includegraphics[width=0.32\textwidth ]{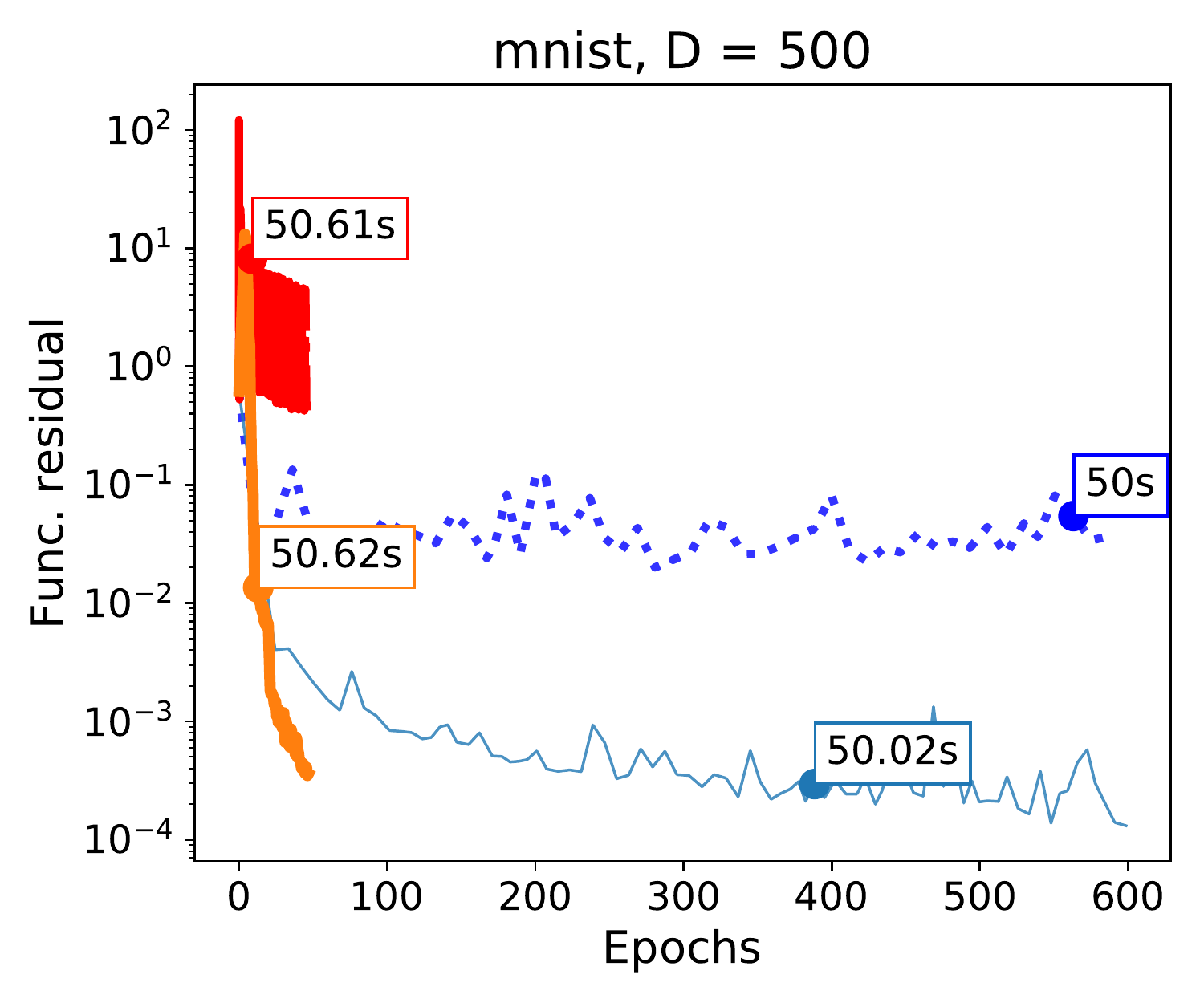}
	
	\includegraphics[width=0.32\textwidth ]{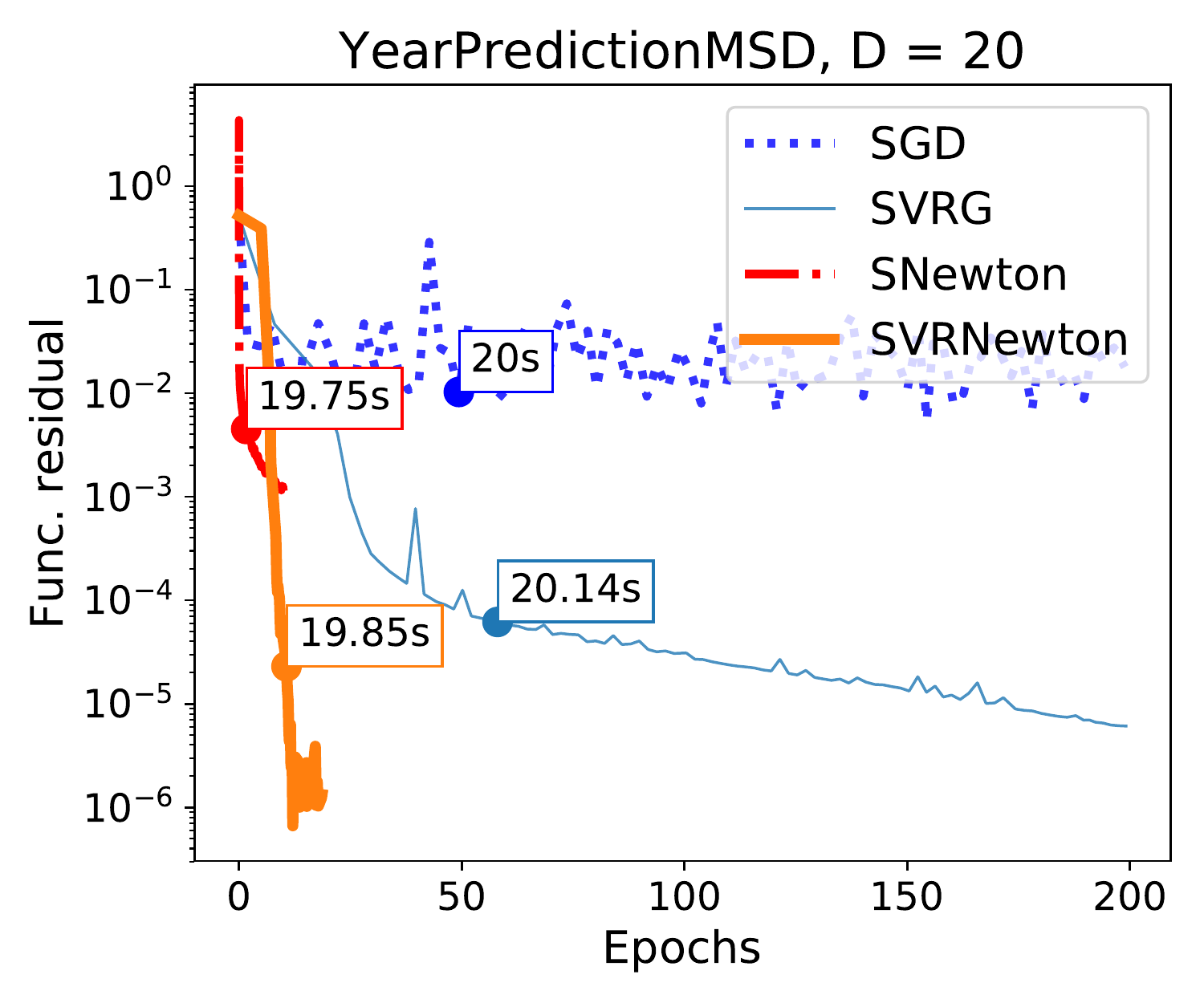}
	\includegraphics[width=0.32\textwidth ]{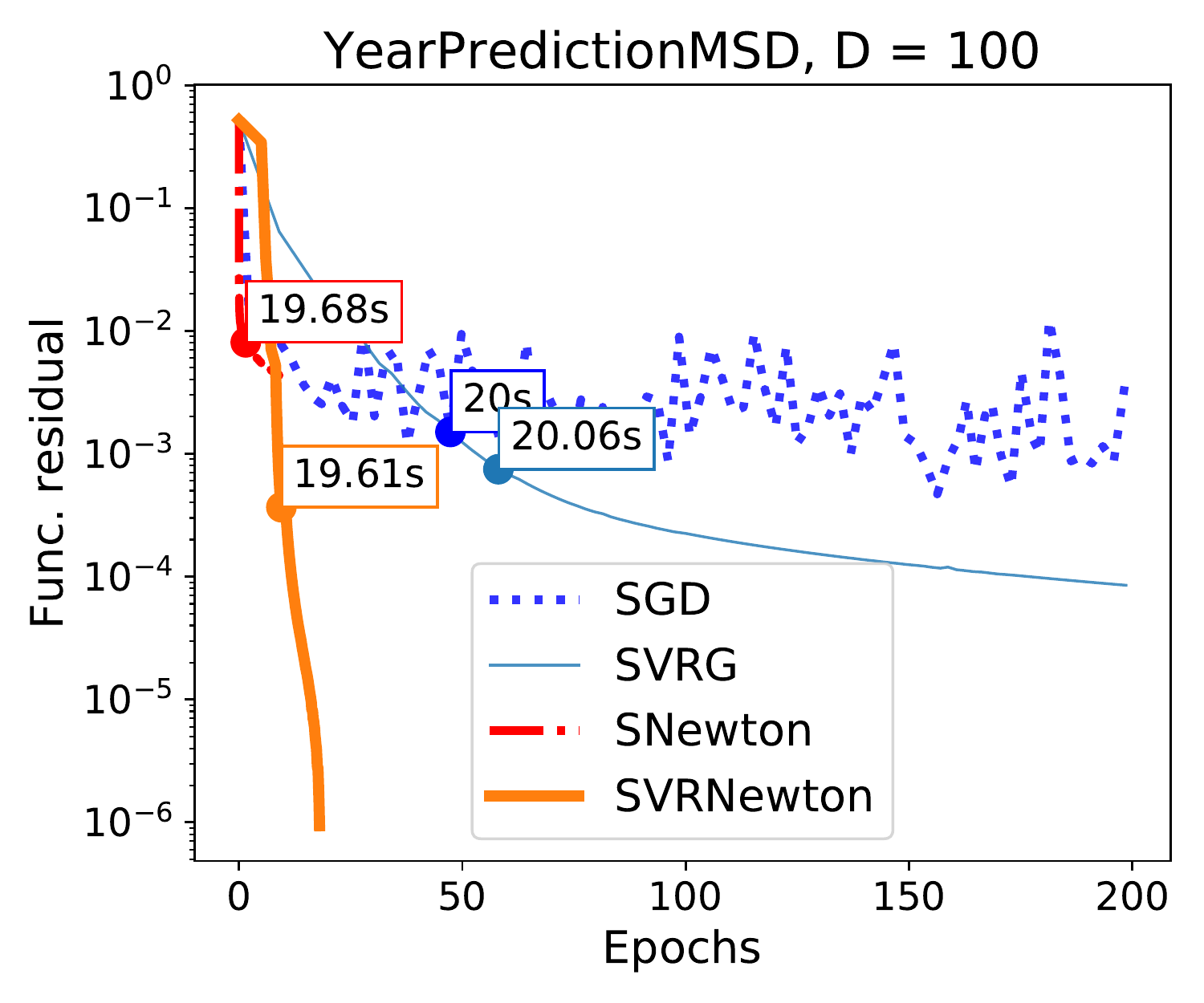}
	\includegraphics[width=0.32\textwidth ]{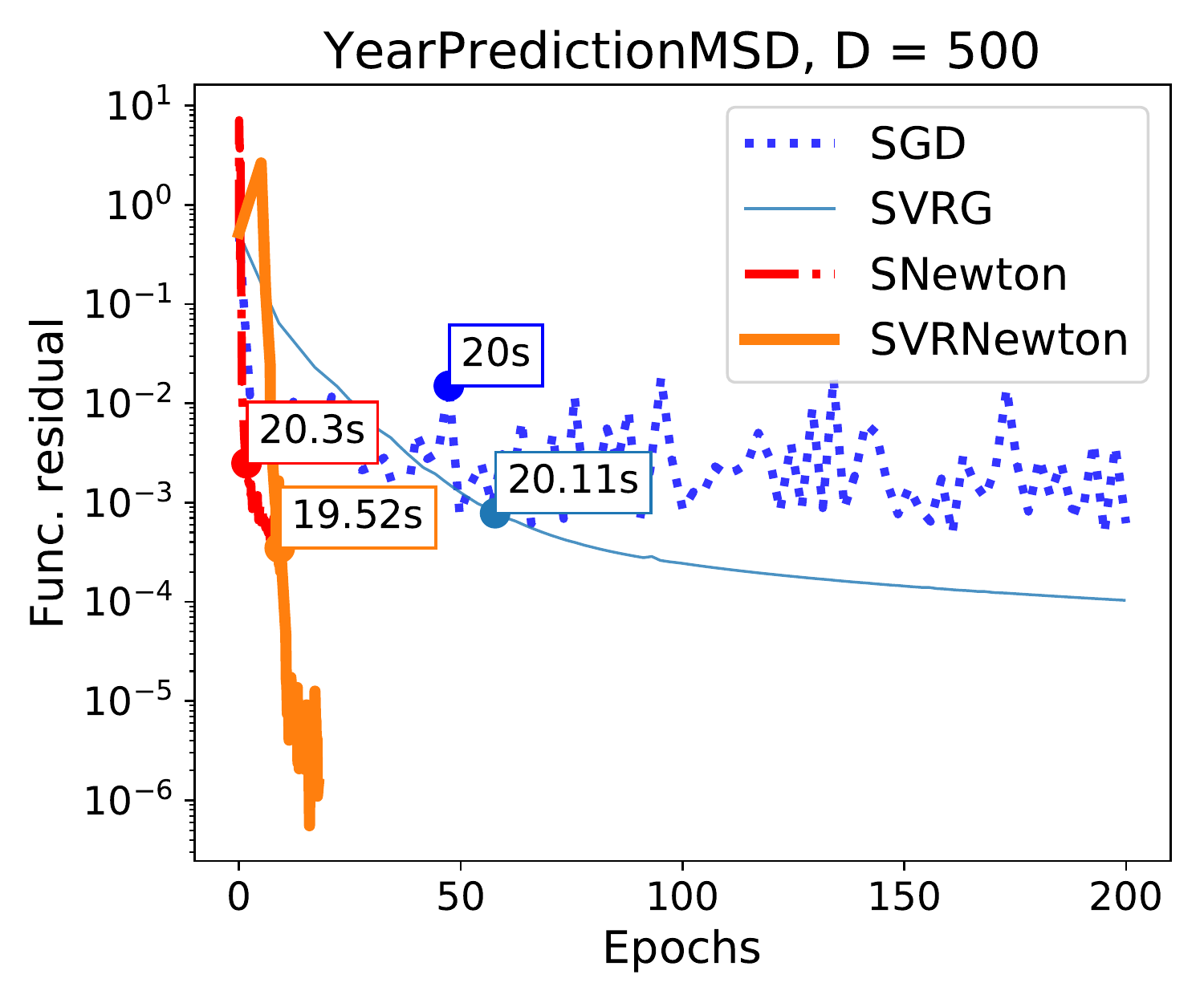}
	
	\includegraphics[width=0.32\textwidth ]{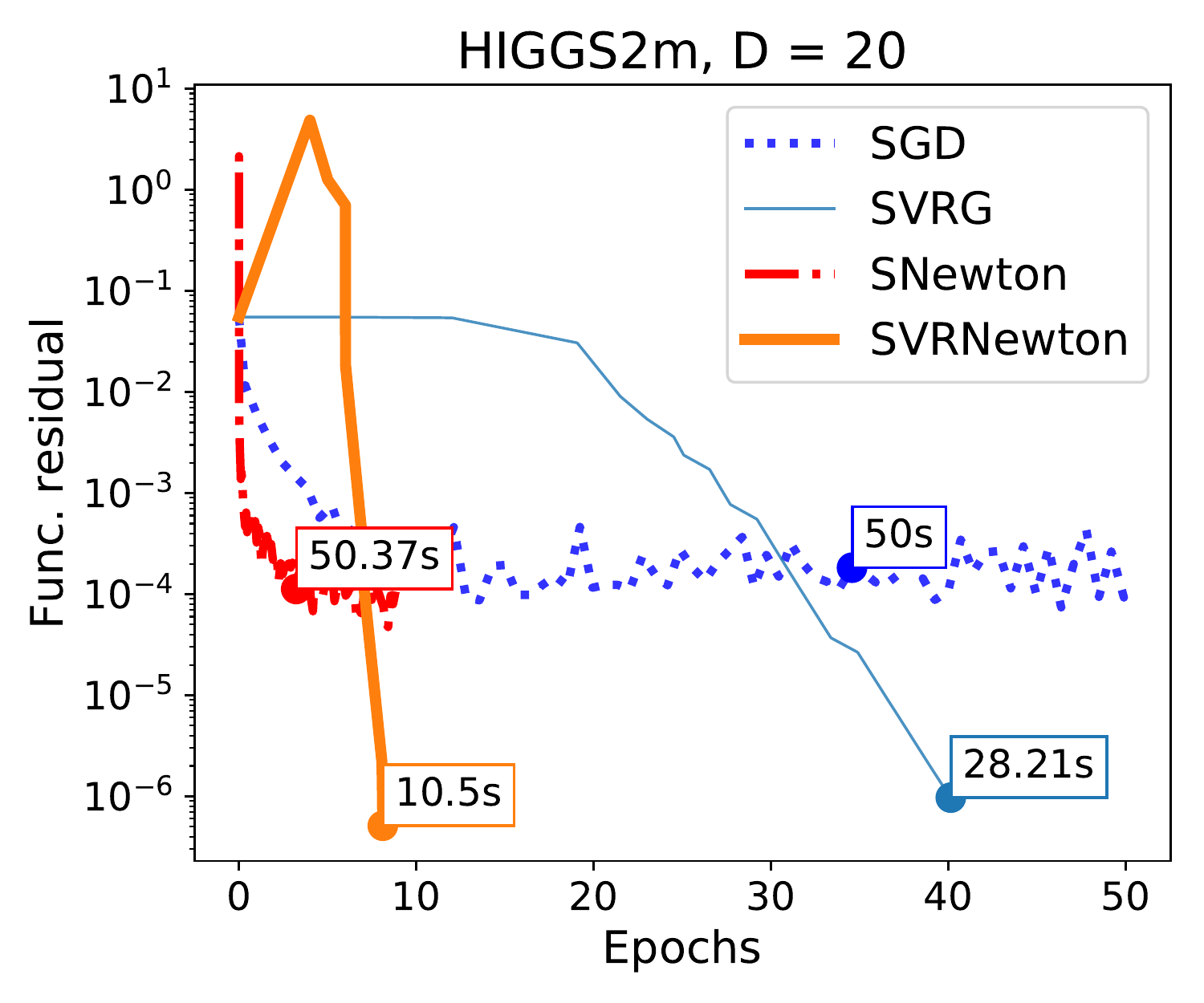}
	\includegraphics[width=0.32\textwidth ]{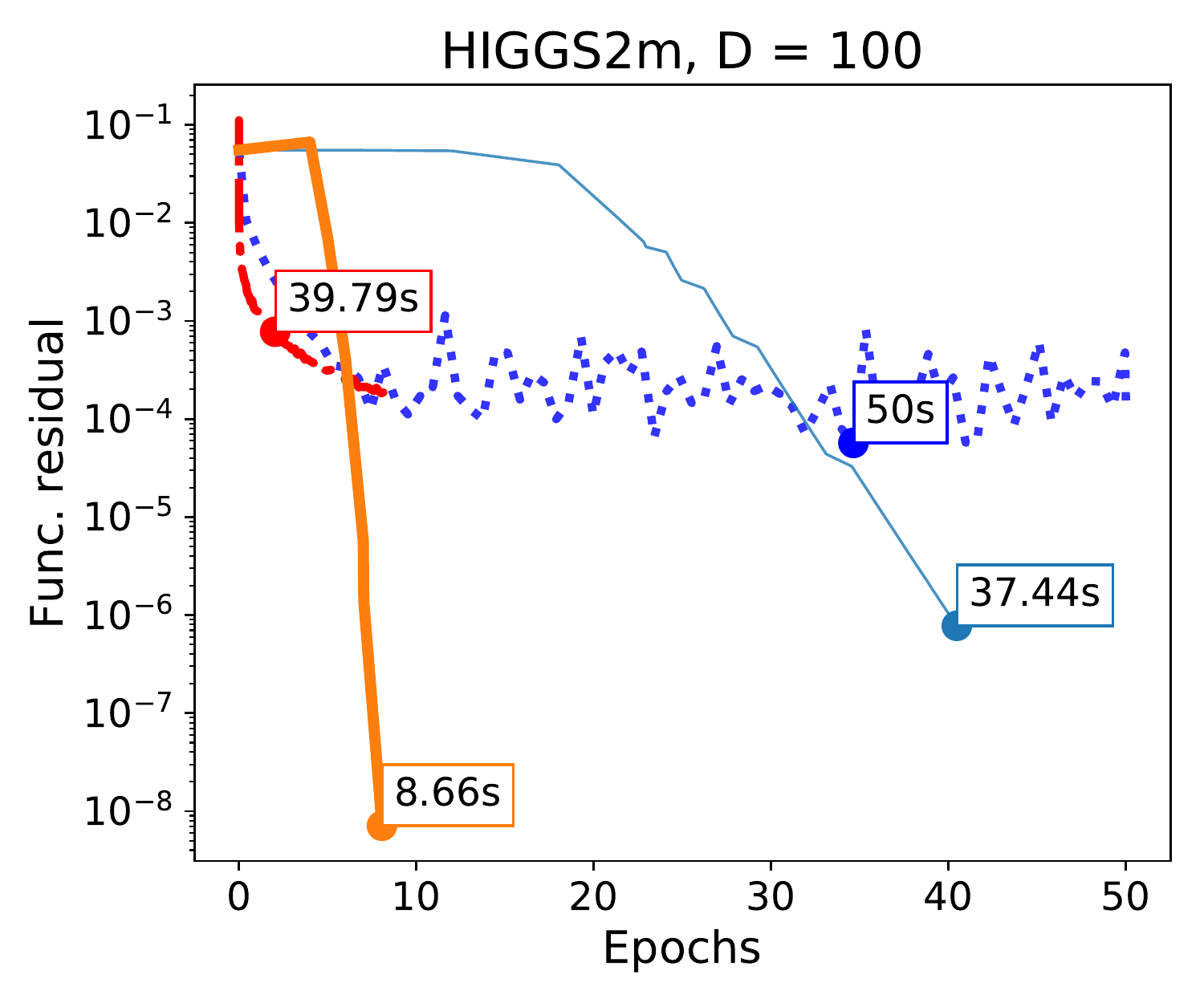}
	\includegraphics[width=0.32\textwidth ]{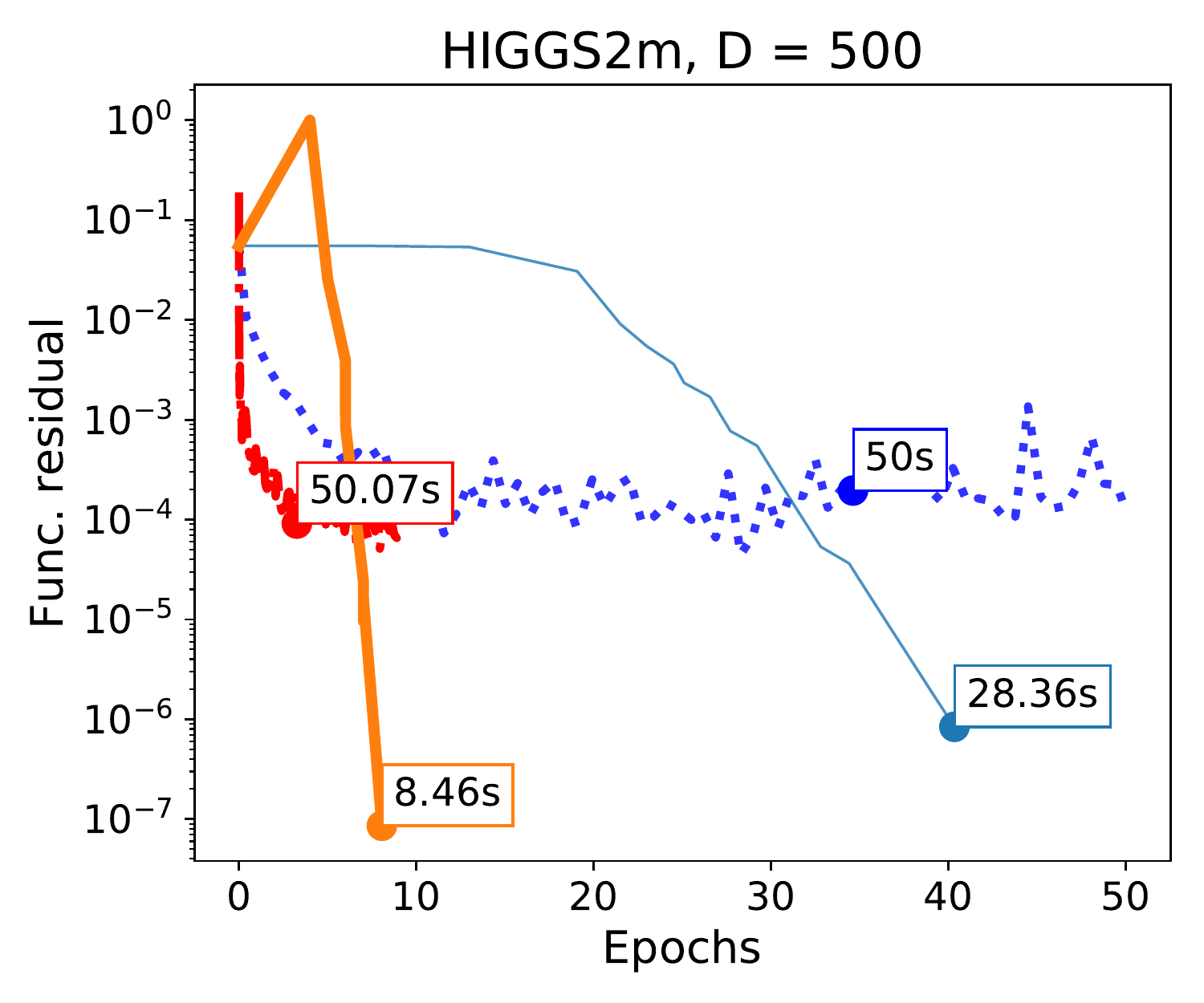}
	
	\caption{Stochastic methods for training logistic regression, datasets: 
	\textit{mnist} $(M = 60000, n = 780)$, 
	\textit{YearPredictionMSD} $(M = 463715, n = 90)$, 
 	\textit{HIGGS2m} $(M = 2 \cdot 10^6, n = 28)$. }
	\label{Fig:logreg_extra_stoch}
\end{figure}

In this section, we provide additional experimental results for the problem
of training Logistic Regression model, regularized by $\ell_2$-ball constraints:
Figure~\ref{Fig:logreg_extra_basic} for the exact methods, and
Figure~\ref{Fig:logreg_extra_stoch} for the stochastic algorithms.

We see, that the second-order schemes usually outperforms first-order methods,
in terms of the number of iterations, and the number of epochs.
Despite the fact, that the Newton step is more expensive,
in many situations we see superiority of the second-order schemes
in terms of the total computational time as well.

\begin{figure}[h!]
	\centering
	\includegraphics[width=0.32\textwidth ]{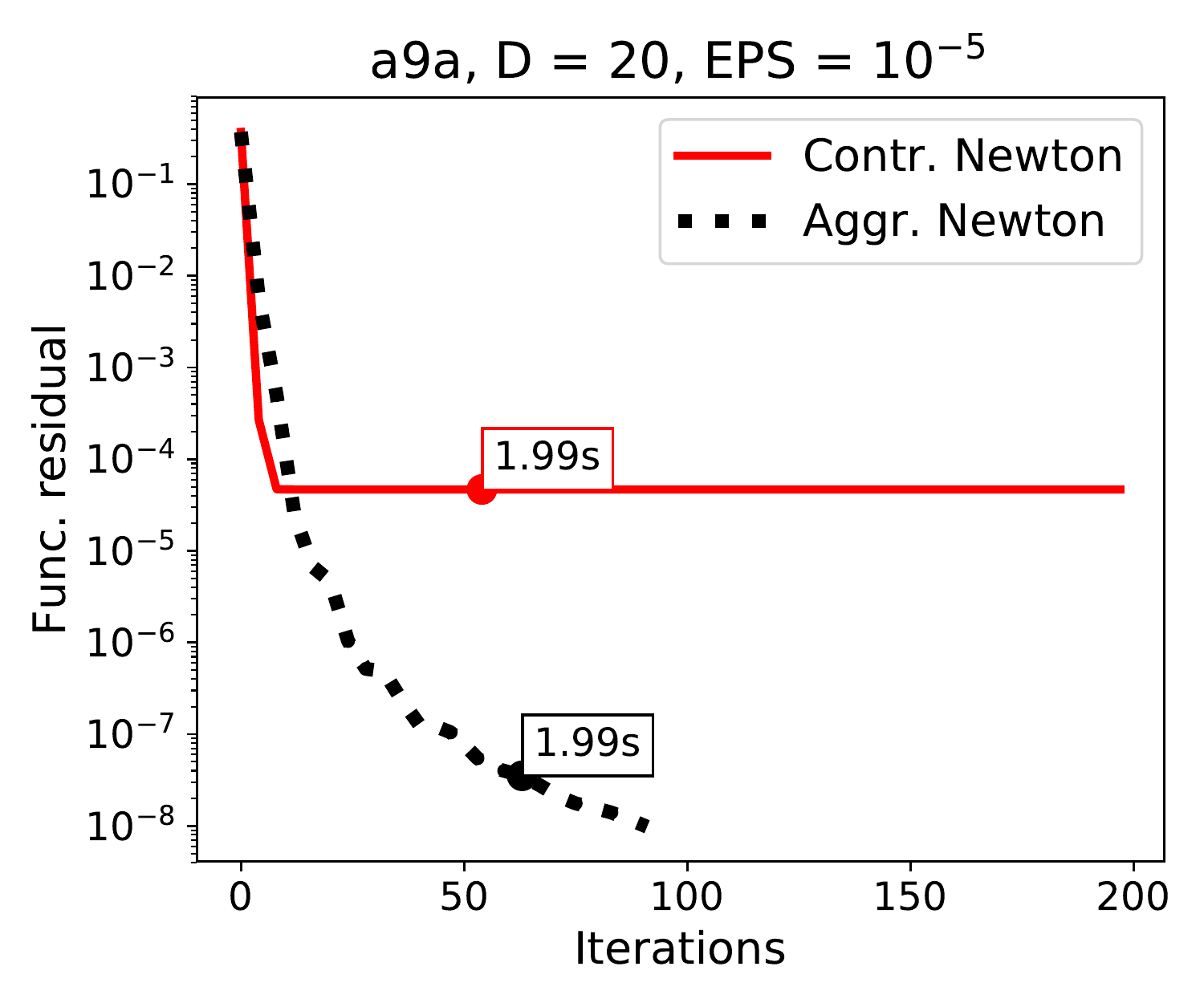}
	\includegraphics[width=0.32\textwidth ]{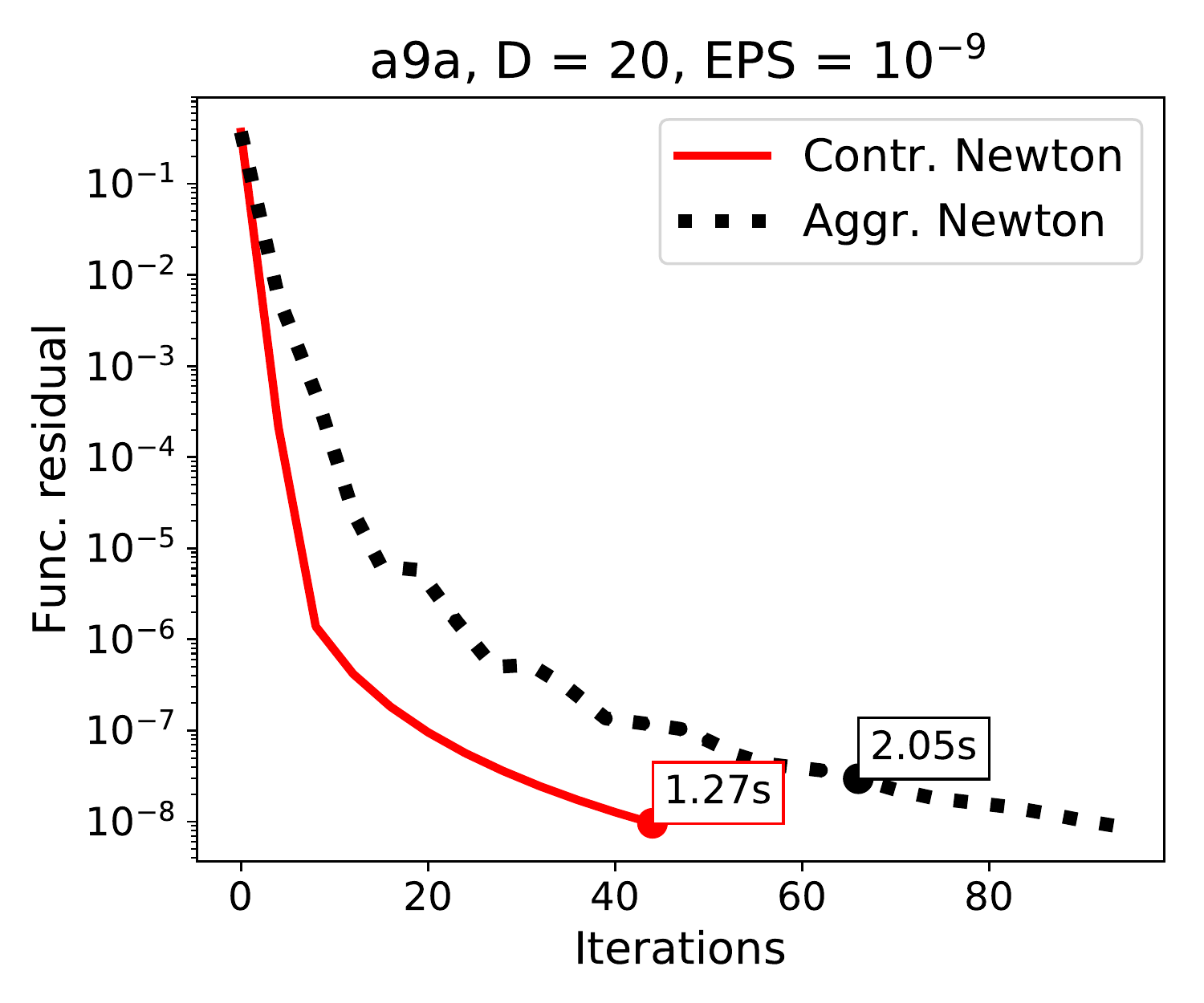}
	
	\includegraphics[width=0.32\textwidth ]{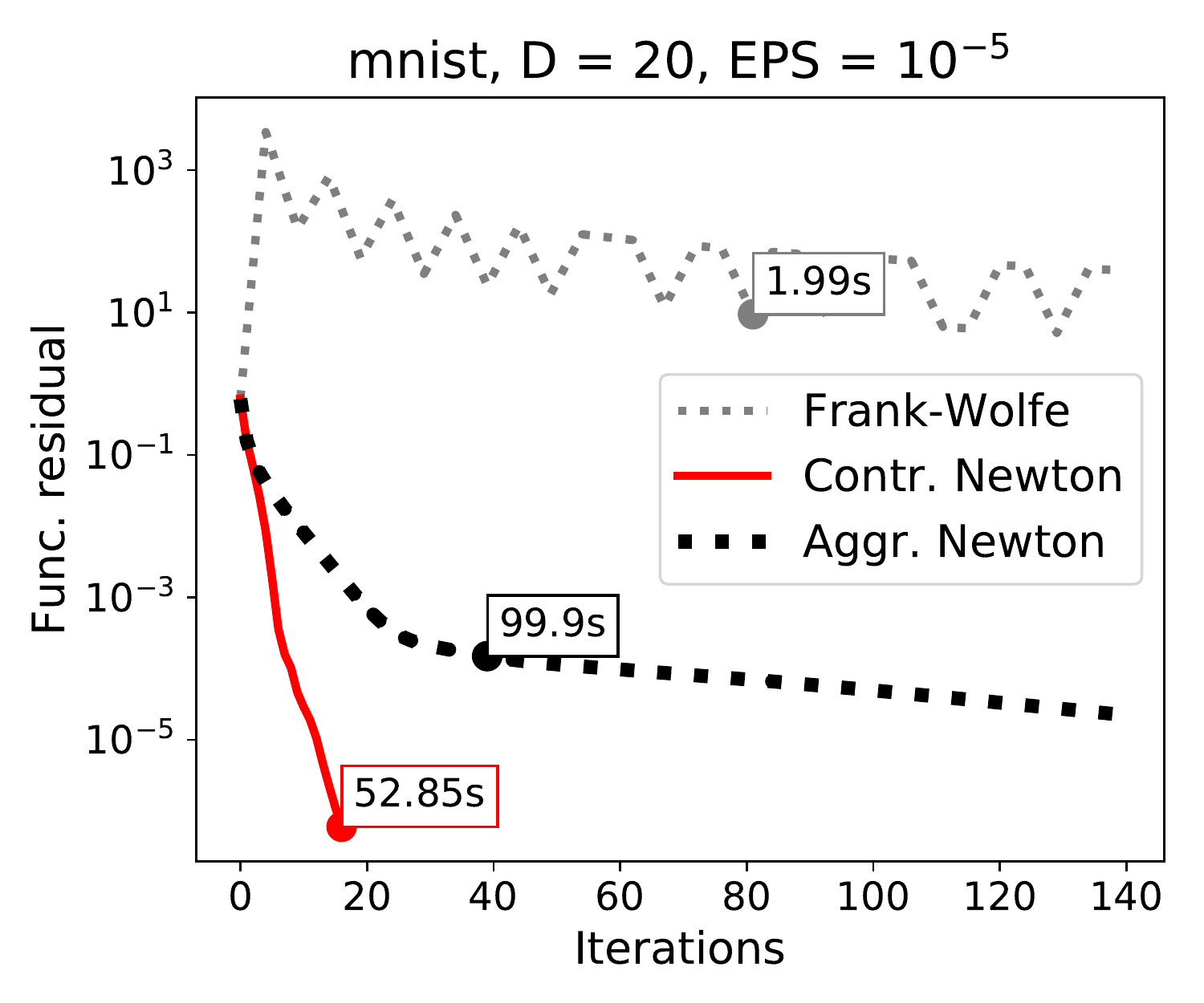}
	\includegraphics[width=0.32\textwidth ]{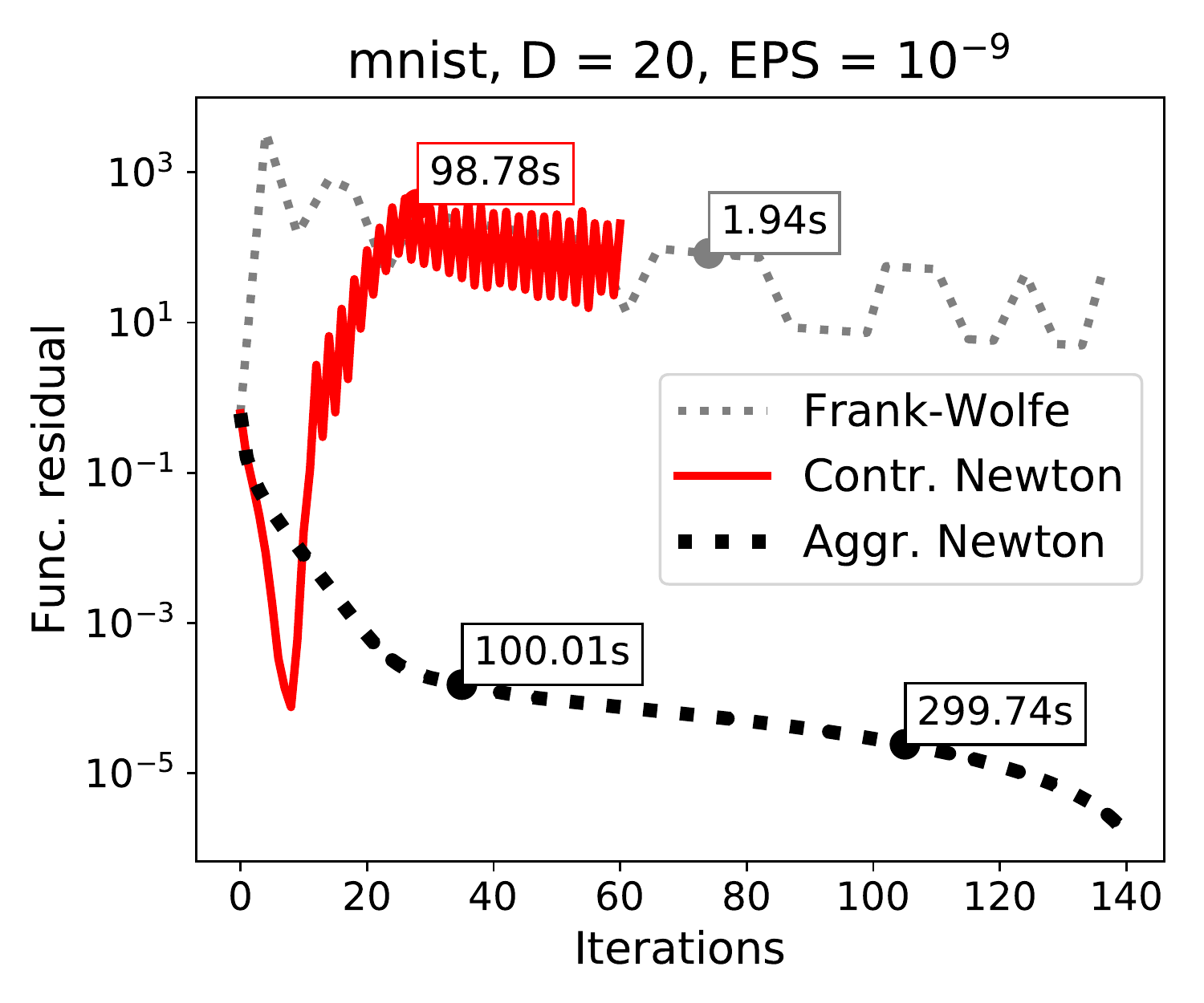}
	\caption{Influence of the parameter of inner accuracy.}
	\label{Fig:accuracy_effect}
\end{figure}

Comparing Contracting-Domain Newton Method (Algorithm~\ref{ContrDomNewton}),
and Aggregating Newton Method (Algorithm~\ref{AggrNewton}), we conclude that both
of the algorithms show reasonably good performance in practice.
The latter one works a bit slower. However,
the aggregation of the Hessians helps to improve numerical stability.
On Figure~\ref{Fig:accuracy_effect}, we demonstrate
influence of the parameter of inner accuracy (EPS),
which we use in our subsolver, on the convergence of the algorithms.
We see much more robust behaviour for Aggregating Newton Method,
while the first algorithm can potentially stop, or even start to diverge,
if the parameter is chosen in a wrong way.

To compute one step of our second-order methods for this task, 
we need to solve subproblem~\eqref{TrustRegionStep}
for $p = 2$. This is minimization of quadratic function over the standard Euclidean ball.
First, we compute \textit{tridiagonal} decomposition of the Hessian 
(it requires $\O(n^3)$ arithmetical operations). 
Then, we solve the dual to our subproblem (which is maximization of one-dimensional concave function) 
by classical Newton iterations (the cost of each iteration is $\O(n)$).
For more details, see Chapter 7 in~\cite{conn2000trust}.

\end{document}